# High-Resolution Weighted Essentially Non-Oscillatory Compact Least-Squares Schemes with Implicit Time Integration for Compressible Navier-Stokes Equations on Curvilinear Grids

**Yongzhi Luo[a,b,1], Huiheng Fan[a,b,1], Wei-Gang Zeng[c], Yu-Xin Ren[d] , Jianhua Pan[a,b*]**

[a] *Zhejiang Provincial Engineering Research Center for the Safety of Pressure Vessel and Pipeline, Ningbo University, Ningbo, 315211, China*

[b] *Key Laboratory of Impact and Safety Engineering, Ministry of Education, Ningbo University, Ningbo 315211, China*

[c] *Academy for Advanced Interdisciplinary Studies, Northeast Normal University, Changchun, 130024, China*

[d] *Department of Engineering Mechanics, Tsinghua University, Beijing 100084, China*

* *Corresponding author: Jianhua Pan. Email: panjianhua@nbu.edu.cn*

[1.] *These authors contributed equally to this work.*

**Abstract**

This paper presents a family of high-resolution weighted essentially non-oscillatory compact least-squares schemes with implicit time integration for the compressible Navier-Stokes equations on curvilinear grids. Compared with the original compact least-squares schemes, the proposed method introduces two main improvements. First, instead of enforcing the accuracy constraints over the entire computational domain, including discontinuous regions, it constructs the reconstruction matrix only along smooth reconstruction lines using an accuracy-preserving weighting strategy. This treatment effectively suppresses the persistent high-frequency oscillations observed in the original compact least-squares schemes and yields sharper profiles near discontinuities. Second, the method simplifies the shock-capturing procedure and improves efficiency by using a single set of polynomials, whereas the original compact least-squares schemes require both unlimited and limited polynomials. Combined with spectral optimization, the proposed method exhibits more favorable spectral properties than conventional weighted essentially non-oscillatory schemes. The smoothness indicators and penalty matrices are constructed through an efficient iterative procedure with modest additional cost. Numerical results for inviscid and viscous one-, two-, and three-dimensional flows demonstrate that the proposed method provides robust shock-capturing capability while maintaining high resolution in smooth regions and across contact discontinuities.

## 1. Introduction

Compressible-flow simulations are ubiquitous in aeronautics and aeroacoustics, including turbulent combustion and mixing in scramjet engines and aerodynamic noise around supersonic vehicles. Over the past two decades, high-order methods have attracted increasing attention for such simulations. Compressible flows involve complex nonlinear phenomena in which shocks, contact discontinuities, and multiscale vortex structures coexist within the computational domain. This coexistence poses two main challenges for shock-capturing methods: resolving smooth flow structures with high fidelity and capturing discontinuities sharply and robustly.

The high resolution of smooth flow structures can be achieved by extending the bandwidth of resolved wavenumber region, either by increasing the order of accuracy or by optimizing the spectral properties. The optimization of spectral properties is an appealing approach since it does not introduce additional computational overhead when compared to the original schemes. The schemes that have superior spectral properties can be roughly categorized into two classes. The first class is the dispersion relation preserving (DRP) schemes pioneered by Tam and Webb [1] and followed by Cheong and Lee [2], Sun et al. [3,4], Li and Ren [5], Zeng et al. [6], etc. In the DRP schemes, one or more free interpolation parameters are optimized to reach a broader resolved wavenumber region. The second class is the compact schemes pioneered by Lele [7] and followed by Deng and Maekawa [8], Guo et al. [9,10], etc. The reconstruction stencil of compact schemes is all the control volumes along the reconstruction axes in essence, and the bandwidth of resolved wavenumber region is greatly increased. Besides that, Wang et al. [11–13] proposed a series of compact least-squares (CLS) reconstruction based on the minimization of variational functions. The CLS schemes are compact while possessing free parameters which are beneficial to achieve optimized spectral properties without sacrificing the order of accuracy.

This work focuses on the development of CLS schemes. Although the compact schemes can reach a broader range of resolved wavenumbers, a shock-capturing technique is essential to handle the discontinuities, which is the second challenge in compressible flow simulations. In the original CLS schemes, a gradient limiter [11–15] is utilized as a post-processing technique to smooth out the high-frequency oscillations produced by Gibbs phenomenon near discontinuities. Thus, two sets of polynomials are maintained in the CLS framework: one set is the polynomials obtained through the accuracy enforcement which we refer to as the linearly reconstructed polynomials, the other set consists of polynomials obtained by applying a slope limiter on the linear reconstructed polynomials and is called post-processed polynomials. If we define the axis along which the accuracy condition in CLS method is enforced as the reconstruction lines, the reconstruction lines in the CLS method extend across the entire computational domain, regardless of shocks or contact discontinuities. This is not mathematically justified since the accuracy condition becomes invalid across discontinuities. Consequently, the high-frequency oscillations may spread over to and deteriorate the smooth regions, which is a major challenge in a compact scheme. In the CLS schemes, a damping mechanism by the WBAP limiter [14,15] is imposed on the linear-reconstructed polynomials. However, this approach has a limitation, since it heavily relies on the performance of the limiter. As the order of reconstructed polynomials increases, it becomes harder for a gradient limiter to smooth out all the high-frequency components, as shown by Fig. 4 in the following section.

Therefore, a fundamental modification of the compact reconstruction procedure is needed to mitigate Gibbs-type oscillations and to remove the need for two separate sets of polynomials in the CLS framework, i.e., the reconstruction matrix of the compact scheme can only be constructed along smooth reconstruction lines which are separated by discontinuities. In detail, the accuracy relation can only be enforced along smooth lines, and an appropriate boundary condition needs to be imposed on the two end points of a smooth reconstruction line. In the shock capturing framework, the exact location of discontinuities is not explicitly defined. Alternatively, a weighting technique can be utilized to define such location of discontinuities implicitly. If a face of a control volume is assumed to locate inside a discontinuity, then the importance of accuracy condition associated with this face is decreased asymptotically to zero and boundary conditions which enforce low gradients at the ends of a smooth reconstruction lines are applied. Similar ideas can be found in the literature, e.g., Guo et al. [9,10] proposed to incorporate a biased lower-order stencil near discontinuities. However, the utilization of lower-order polynomials degrades the high resolution of the original scheme.

In [16], the authors proposed a novel third-order weighted CLS scheme as a first trial to implement the idea of “only reconstructing along smooth lines”. High resolution of the WCLS scheme was demonstrated in the work [16] where a third-order scheme resolves the shock at a resolution comparable to that of the fifth-order WENO-Z scheme [17]. However, several drawbacks still exist in the work of [16]. First, an additional third-party smoothness indicator is essential to identify the discontinuous interface in the flow. The smoothness indicator requires information from control volumes beyond von Neumann neighbors, introducing additional computational complexity and destroying the compact-stencil property of CLS schemes. Second, solving block-tridiagonal matrices is needed to obtain the polynomials for each control volume, and the computational overhead is significant, especially when dealing with multi-dimensional Navier-Stokes (NS) equations utilizing characteristic decomposition. Additionally, when extending to parallel simulations, connectivity between different processes should be taken care of. Parallel implementation is also challenging for other compact schemes [18,19].

This work extends the third-order WCLS scheme to higher-order accurate formulations, integrating with implicit time integration on curvilinear grids and resolving the three drawbacks in the previous work of [16]. The proposed implicit WCLS schemes are not a straightforward extension from the explicit third-order scheme. First, different from our previous work [16], smoothness indicators are now calculated only based on the information from von Neumann neighbors iteratively with existing derivative information possessed by the CLS schemes. However, with the traditional smoothness indicator by Jiang and Shu [20], the WCLS schemes can hardly converge by iteration near discontinuities. To cure this problem, an additional term based on an inexpensive linear reconstruction is added to the smoothness indicator. Second, additional boundary conditions behaving as the penalty matrices are proposed to apply on the ends of smooth reconstruction lines. The penalty matrices are self-adjusted through boundary value differences (BVDs). If the interface is located inside a discontinuity, the BVD becomes large, and the penalty matrix stabilizes the polynomials at this interface. If the interface is located inside a smooth reconstruction line, the BVD goes to zero and the formal order of accuracy is preserved even near smooth extrema. This novel mechanism keeps the proposed scheme not only essentially non-oscillatory but also high-resolution in smooth regions as shown in Fig. 4. Regarding the high-computational cost faced by the explicit third-order WCLS method [16], the coupling solution of implicit time integration for NS equations and the implicit iteration of CLS reconstruction system has been shown to be an effective approach to tackle the computational overhead faced by CLS reconstruction [11–13,21]. In this work, within the framework of implicit time integration, the reconstruction system of WCLS schemes is solved by generalized minimal residual method (GMRES). One step of implicit time iteration for NS equations and one step of GMRES method for WCLS reconstruction are conducted subsequently in one iteration loop and the solution of block tri-diagonal matrices is avoided. The boundary treatment between different processes faced by traditional compact schemes is unnecessary since the derivatives on virtual control volumes near connectivity can be synchronized during iterations. Drawing on all the ingredients, this paper proposes a series of high-resolution compact WCLS schemes that simultaneously achieve two seemingly contradictory goals: preserving high-order accuracy in smooth regions without polluting the solution and capturing discontinuities with sharp, non-oscillatory profiles. The proposed method appears promising for simulations of compressible turbulent flows.

The remainder of the paper is organized as follows. Section 2 introduces the basic idea of the WCLS schemes. Section 3 extends the WCLS schemes to the Navier-Stokes equations and presents the coupled solution procedure for implicit time integration and WCLS reconstruction. Section 4 presents numerical examples for inviscid Euler equations and viscous Navier-Stokes equations on curvilinear uniform and non-uniform grids to assess the resolution, robustness, and efficiency of the proposed WCLS schemes. Conclusions are given in Section 5.

## 2. Weighted Compact Least-Squares Scheme

To illustrate the basic idea of the weighted compact least-squares schemes, we take the one-dimensional (1D) advection equation as an example:

$$\frac{\partial u}{\partial t}+\frac{\partial f}{\partial x}=0. \tag{1}$$

$x$ and $t$ denote the spatial and temporal coordinates, respectively. $u$ is a scalar and $f = \lambda u$ is the flux function. $\lambda$ is constant. The computational domain denoted as $\Omega=[x_l,x_r]$ is subdivided into $N$ non-overlapping control volumes. The $i$-th control volume is denoted as $\Omega_i=[x_{i-1/2},x_{i+1/2}]$. $x_{-1/2}=x_l$ and $x_{N+1/2}=x_r$. The integration of Eq. (1) over control volume $\Omega_i$ results in the following semi-discretized form in finite volume method,

$$h_i\frac{\partial \overline{u}_i}{\partial t}+f_{i+1/2}-f_{i-1/2}=0, \tag{2}$$

where $h_i=x_{i+1/2}-x_{i-1/2}$. For control volume $\Omega_i$, a $k_{\text{th}}$-order polynomial is assumed as

$$p_i(x)=\overline{u}_i+\sum_{j=1}^{k}a_{i,j}\phi_{i,j}\,, \tag{3}$$

where $\phi_{i,j}$ is the $j$-th zero-mean basis function of control volume $\Omega_i$. $\phi_{i,j}$ is defined as

$$\phi_{i,j}=\left(\frac{x-x_{c,i}}{h_i}\right)^j-\frac{1}{h_i}\int_{x_{i-1/2}}^{x_{i+1/2}}\left(\frac{x-x_{c,i}}{h_i}\right)^j\mathrm{d}x\,, \tag{4}$$

where $x_{c,i}=\frac{1}{2}\left(x_{i-1/2}+x_{i+1/2}\right)$ is the center of $\Omega_i$ and $a_{i,j}, j=1,2,\cdots,k$ are unknowns to be reconstructed. Once $a_{i,j}$ are determined, the flux function $f_{i+1/2}$ is replaced by the numerical flux as

$$f_{i+1/2}=\frac{\lambda}{2}\left(u_{i+1/2}^L+u_{i+1/2}^R\right)-\frac{|\lambda|}{2}\left(u_{i+1/2}^R-u_{i+1/2}^L\right), \tag{5}$$

where $u_{i+1/2}^L=p_i(x_{i+1/2})$ and $u_{i+1/2}^R=p_{i+1}(x_{i+1/2})$. Equation (1) can be integrated by the multi-stage Runge-Kutta (RK) method [22–24].

The key innovation of the present work lies in the WCLS reconstruction of polynomials $p_i(x)$, $i=1,2,\cdots,N$ which are not only highly accurate but also non-oscillatory across discontinuities. To achieve this, a loss function is first defined for the control volume $\Omega_i$ as

$$\mathbb{I}_i=\mathbb{I}_{i,\mathrm{acc}}+\mathbb{I}_{i,\mathrm{pen}}+\mathbb{I}_{i,\mathrm{bc}}, \tag{6}$$

$$\mathbb{I}_{i,\mathrm{acc}}=\sum_{n\in\mathrm{ST}_i}\left(\beta_n\sum_{j=0}^{k}W_j\,l_\alpha^{2j}\left(\left.\frac{\mathrm{d}^jp_i}{\mathrm{d}x^j}\right|_{x_\alpha}-\left.\frac{\mathrm{d}^jp_n}{\mathrm{d}x^j}\right|_{x_\alpha}\right)^2\right), \tag{7}$$

$$\mathbb{I}_{i,\mathrm{pen}}=\sum_{n\in\mathrm{ST}_i}\left(\theta_n\sum_{j=1}^{k}\omega_j l_\alpha^{2j}\left(\left.\frac{\mathrm{d}^jp_i}{\mathrm{d}x^j}\right|_{x_\alpha}\right)^2\right), \tag{8}$$

$$\mathbb{I}_{i,\mathrm{bc}}=\sum_{j=0}^{k}W_{\mathrm{bc},j}l_i^{2j}\left(\left.\frac{\mathrm{d}^jp_i}{\mathrm{d}x^j}\right|_{x_{\mathrm{bc}}}-\left.\frac{\mathrm{d}^jp}{\mathrm{d}x^j}\right|_{x_{\mathrm{bc}}}\right)^2, \tag{9}$$

where $\mathrm{ST}_i$ is the set of indices for von Neumann neighbors of $\Omega_i$. For the 1D case in this paper, $\mathrm{ST}_i=\{i-1,i+1\}$ and $x_\alpha=\Omega_i\cap\Omega_n$.

$\mathbb{I}_{i,\mathrm{acc}}$ consists of the continuous equations for all derivatives between polynomial $p_i(x)$ in $\Omega_i$ and $p_n(x)$ in its neighboring $\Omega_n$. $l_\alpha$ is the length scale utilized to keep terms of different orders dimensionally consistent and is chosen as $l_\alpha = h_i$. $W_j$, $j = 1,2,\cdots,k$ are linear coefficients which can be optimized to reach a broad range of resolved wavenumbers [16,25]. $\beta_n$ is the weight of interface $\Omega_i \cap \Omega_n$. $\mathbb{I}_{i,\mathrm{bc}}$ is the boundary term and is included if $\Omega_i$ includes boundaries. $\mathrm{d}^j p/\mathrm{d}x^j$ is the desired $j$-th derivative at boundary point $x_{\mathrm{bc}}$ and $W_{\mathrm{bc},j}$ is the respective weight for the boundary condition.

For CLS schemes, the loss function constitutes only $\mathbb{I}_{i,\mathrm{acc}}$ (with $\beta_n$ all equals to 1) and $\mathbb{I}_{i,\mathrm{bc}}$ (if the interface is a boundary), thus the reconstruction line extends along the axes to the whole computational domain even across discontinuities. However, this is mathematically invalid since the derivatives are discontinuous on the two sides of discontinuities. Enforcing the continuous equations across discontinuity is the source of Gibbs phenomenon in a compact scheme.

To cure the Gibbs phenomenon, the weight $\beta_n$ should be designed to approach 0 if the interface $\Omega_i \cap \Omega_n$ is located inside a discontinuity. In this way, the reconstruction lines of compact schemes are separated into a series of smooth segments. The separation of reconstruction lines into smooth segments is only the first step. As observed numerically, undesired gradients still emerge at the two ends if no boundary conditions are enforced for the smooth reconstruction lines. A natural idea is to apply a penalty term $\mathbb{I}_{i,\mathrm{pen}}$ when $\beta_n \to 0$. In $\mathbb{I}_{i,\mathrm{pen}}$, $\omega_j = (1/j!)^2$ and $\theta_n$ is the weight of the penalty term. When $\mathbb{I}_{i,\mathrm{pen}}$ takes effect, the derivatives approach zero near $x_\alpha$, which further enhances non-oscillatory properties of the proposed schemes. Importantly, $\theta_n$ should be carefully designed to avoid order degradation in smooth regions, especially near extrema.

Taking derivatives of $\mathbb{I}_i$ over unknown coefficients $a_{i,j}$, $j = 1,2,\cdots,k$, the reconstruction matrix for $\Omega_i$ can be derived as

$$\boldsymbol{M}_i^{(-1)}\boldsymbol{a}_{i-1} + \boldsymbol{M}_i^{(0)}\boldsymbol{a}_i + \boldsymbol{M}_i^{(1)}\boldsymbol{a}_{i+1} = \boldsymbol{b}_i, \tag{10}$$

where $\boldsymbol{a}_i = \left[a_{i,1}, a_{i,2}, \cdots, a_{i,k}\right]^{\mathrm{T}}$, $\boldsymbol{M}_i^{(0)} = \boldsymbol{M}_{i,\mathrm{acc}}^{(0)} + \boldsymbol{M}_{i,\mathrm{pen}}^{(0)} + \boldsymbol{M}_{i,\mathrm{bc}}^{(0)}$, $\boldsymbol{b}_i = \boldsymbol{b}_{i,\mathrm{acc}} + \boldsymbol{b}_{i,\mathrm{bc}}$, $\boldsymbol{M}_i^{(m)} = \boldsymbol{M}_{i,\mathrm{acc}}^{(m)}, m = -1,1$. $\boldsymbol{M}_{i,\mathrm{acc}}^{(m)}, m = -1,0,1$ and $\boldsymbol{b}_{i,\mathrm{acc}}$ result from $\partial\mathbb{I}_{i,\mathrm{acc}}/\partial\boldsymbol{a}_i$; $\boldsymbol{M}_{i,\mathrm{bc}}^{(0)}$ and $\boldsymbol{b}_{i,\mathrm{bc}}$ result from $\partial\mathbb{I}_{i,\mathrm{bc}}/\partial\boldsymbol{a}$; $\boldsymbol{M}_{i,\mathrm{pen}}^{(0)}$ results from $\partial\mathbb{I}_{i,\mathrm{pen}}/\partial\boldsymbol{a}$. $[\cdot]^{\mathrm{T}}$ denotes the transpose of a vector. It should be noted that when taking derivatives of the loss function $\mathbb{I}$ over the unknowns $\boldsymbol{a}_i$, the smoothness indicator $\beta$ and the coefficients of penalty terms $\theta$ are treated as independent of the unknowns $\boldsymbol{a}_i$.

Coupling the reconstruction matrices for all the control volumes, a global sparse block tri-diagonal system is obtained and denoted as

$$\boldsymbol{M}\boldsymbol{a} = \boldsymbol{b}, \tag{11}$$

where $\boldsymbol{a} = [\boldsymbol{a}_1^{\mathrm{T}}, \boldsymbol{a}_2^{\mathrm{T}}, \cdots, \boldsymbol{a}_N^{\mathrm{T}}]^{\mathrm{T}}$, $\boldsymbol{b} = [\boldsymbol{b}_1^{\mathrm{T}}, \boldsymbol{b}_2^{\mathrm{T}}, \cdots, \boldsymbol{b}_N^{\mathrm{T}}]^{\mathrm{T}}$. $\boldsymbol{M}_i^{(-1)}$, $\boldsymbol{M}_i^{(0)}$ and $\boldsymbol{M}_i^{(1)}$ are the sub diagonal, diagonal, and super diagonal elements of $\boldsymbol{M}$, respectively.

The direct solution of reconstruction linear system of Eq. (11) results in computational overhead compared with the traditional WENO schemes, especially when later extended to NS equations utilizing characteristic decomposition where different conservation variables are coupled with each other [16]. In this work, Eq. (11) is solved by the GMRES method and is embedded into the implicit time integration. In the coupling solution strategy [11,13,21,26], the implicit iteration of NS equations and the searching for solutions in Krylov subspace of reconstruction system are conducted subsequently only once in a single iteration step. The details of this coupling solution strategy are presented in the next section. In this way, the WCLS schemes avoid the computational overhead led by the direct solution of block tri-diagonal linear system and achieve comparable computational complexity in a single iteration step with the traditional WENO schemes.

Now, let us introduce the design of key parameters $\beta_n$ and $\theta_n$ in Eqs. (7) and (8).

$\beta_n$ is set as

$$\mathrm{SI}_n = \sum_{j=3}^{k} \frac{h_i^{2j}}{(j!)^2} \int_{x_{i-1/2}}^{x_{i+1/2}} \left(\frac{\mathrm{d}^j p_n}{\mathrm{d}x^j}\right)^2 \mathrm{d}x + \frac{h_i^2}{2}\left(\int_{x_{i-1/2}}^{x_{i+1/2}} \left(\frac{\mathrm{d}p_n}{\mathrm{d}x}\right)^2 \mathrm{d}x + \int_{x_{i-1/2}}^{x_{i+1/2}} \left(\frac{\mathrm{d}q_n}{\mathrm{d}x}\right)^2 \mathrm{d}x\right), \tag{12}$$

$$\phi_n = \frac{1}{\mathrm{SI}_n + \varepsilon}, \tag{13}$$

$$\beta_n = \phi_n / \sum_{m\in \mathrm{ST}_i} \phi_m, n \in \mathrm{ST}_i, \tag{14}$$

where $\varepsilon$ is a small number avoiding division by zero and is chosen as $\max(1.0 \times 10^{-16}(|u_i| + |u_{i-1}| + |u_{i+1}|)^2/9, 1 \times 10^{-50})$; $q_i(x)$ is a first-order polynomial reconstructed in $\Omega_i$ as

$$q_i(x) = \bar{u}_i + \frac{\Delta x_{i-1/2}(\bar{u}_i - \bar{u}_{i-1}) + \Delta x_{i+1/2}(\bar{u}_{i+1} - \bar{u}_i)}{\Delta x_{i-1/2}^2 + \Delta x_{i+1/2}^2}(x - x_{c,i}), \tag{15}$$

where $\Delta x_{i\pm 1/2} = (h_{i\pm 1} + h_i)/2$.

The smoothness indicator SI in Eq. (12) is different from the one in WENO schemes proposed by Jiang and

Shu [20,27]. The smoothness indicator in the form of Jiang and Shu is

$$\mathrm{SI}_n = \sum_{j=1}^{k} h_i^{2j} \int_{x_{i-1/2}}^{x_{i+1/2}} \left(\frac{\mathrm{d}^j p_n}{\mathrm{d}x^j}\right)^2 \mathrm{d}x. \tag{16}$$

There are two key improvements in the proposed form of Eq. (12). First, Eq. (12) includes the first-order derivatives from $q_i(x)$ and the differences of neighboring cell-averaged values. $q_i(x)$ and the differences of neighboring cell-averaged values are essential to start the iteration smoothly. Consider a control volume on the left side of an ideal discontinuity of a step function. At the beginning of iteration, the derivatives on all control volumes are assigned as 0, then the two smoothness indicators of this control volume are equal, and oscillations are triggered at the very first iteration step. Once oscillations are triggered, the smoothness indicators are biased towards the smoother side and the derivatives in the control volume are smoothed out again in the next iteration. Eventually, the system can hardly converge to a steady state. This is illustrated in Fig. 1. Second, the proposed smoothness indicator includes a reduced factor $1/(j!)^2$ for higher-order derivatives. The idea of reducing the weights of higher-order derivatives was originally proposed by Li et al. in the work of WENO limiter for discontinuous Galerkin method [28]. Numerical experiments indicate that the factor $1/(j!)^2$ is beneficial to improve the resolution in smooth regions. Here, the evaluation of Eq. (12) is based on the known intermediate polynomial of the WCLS scheme, except for the inexpensive linear polynomial $q_i(x)$.

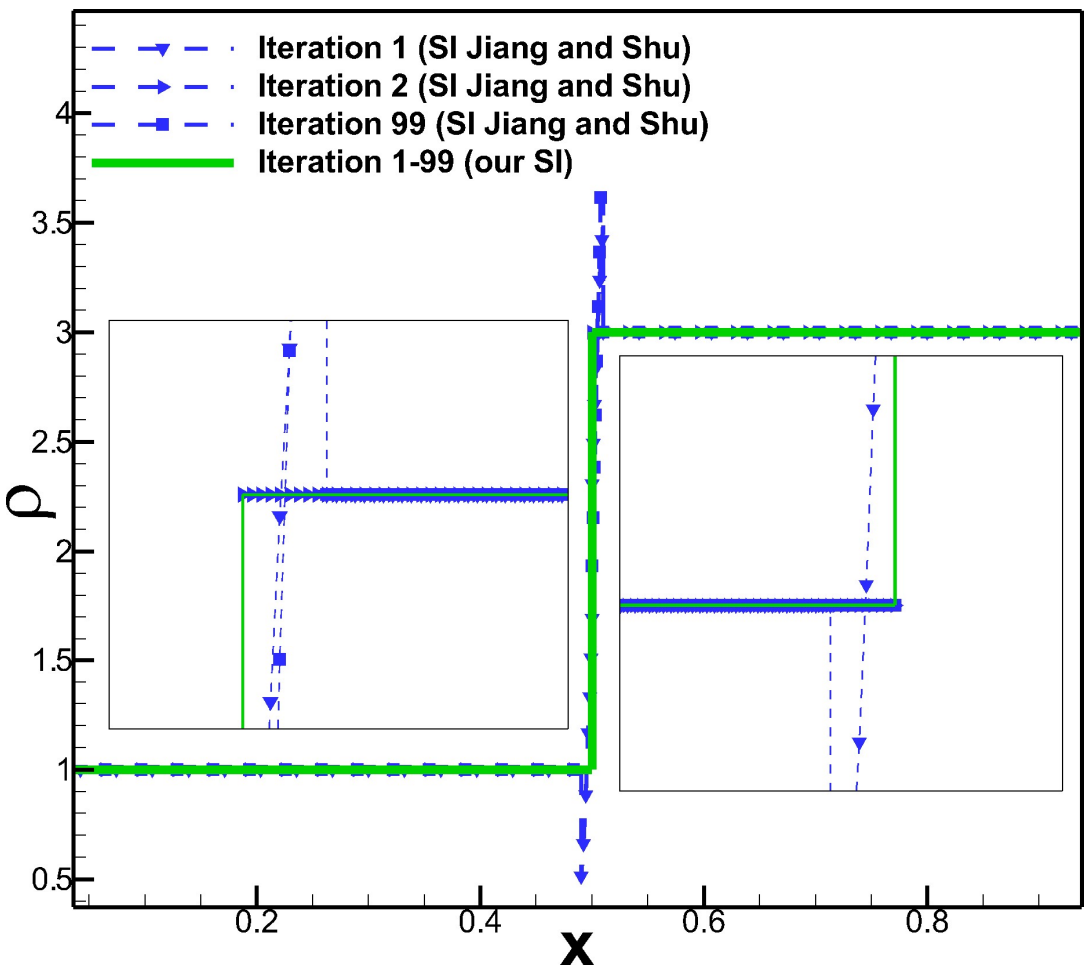

Figure 1: Reconstruction based on the proposed SI of Eq. (12) and SI by Jiang and Shu [20,27].

With the smoothness indicator, a shock detector can be constructed without additional computational cost as

$$\sigma_i = \frac{2\mathrm{SI}_{i,\min}\mathrm{SI}_{i,\max}}{\mathrm{SI}_{i,\min}^2 + \mathrm{SI}_{i,\max}^2}, \tag{17}$$

where

$$\mathrm{SI}_{i,\min} = \min_{n\in \mathrm{ST}_i} \mathrm{SI}_n\,, \mathrm{SI}_{i,\max} = \max_{n\in \mathrm{ST}_i} \mathrm{SI}_n\,. \tag{18}$$

If $\sigma_i \geq \sigma_c$, the flow is treated as continuous and $\beta_n, n \in \mathrm{ST}_i$ are set as equal for neighbors of $\Omega_i$. Note that $\sigma_c$ in this work is only utilized to improve computational efficiency to avoid characteristic decomposition of Euler equations and can be set as conservative as possible.

SI in Eq. (12) depends on $a_{n,j}$, $j = 1,2,\cdots,k$ and is compact, i.e., Eq. (12) only depends on variables within von Neumann neighbors, which is meaningful when further extending the proposed schemes to unstructured grids. However, this also leads to a nonlinear form of Eq. (11). To make Eq. (11) solvable, at the beginning of each iteration step of GMRES, the weights $\beta$ on all control volumes are first calculated and treated as constant within this iteration step. Thus, Eq. (11) is linearized and the Krylov subspace searching algorithm still applies.

$\theta_n$ is the second key parameter and acts in penalty matrix. In smooth regions, $\theta_n = O(\Delta x^k)$ is a sufficient condition for Eq. (3) to achieve a prior $k$-th order of accuracy.

**Theorem 1**. *Without considering the boundary conditions, if $\theta_n = O(\Delta x^k)$, the polynomial $p_i(x)$ reconstructed by the WCLS schemes achieves a prior $k$-th order of accuracy.*

Proof. *If $\theta_n = O(\Delta x^k)$, then all the elements in the penalty matrix $\boldsymbol{M}_{\mathrm{pen}}$ are of order $O(\Delta x^k)$. The solution of Eq. (11) satisfies*

$$
\begin{aligned}
\boldsymbol{a} &= \left(\boldsymbol{M}_{\text{acc}} + \boldsymbol{M}_{\text{pen}}\right)^{-1}\boldsymbol{b}_{\text{acc}} = \left(\boldsymbol{I} + \boldsymbol{M}_{\text{acc}}^{-1}\boldsymbol{M}_{\text{pen}}\right)^{-1}\boldsymbol{M}_{\text{acc}}^{-1}\boldsymbol{b}_{\text{acc}} \\
&= \left(\boldsymbol{I} + \sum_{j=1}^{+\infty}(-1)^j\left(\boldsymbol{M}_{\text{acc}}^{-1}\boldsymbol{M}_{\text{pen}}\right)^j\right)\boldsymbol{M}_{\text{acc}}^{-1}\boldsymbol{b}_{\text{acc}} = \boldsymbol{M}_{\text{acc}}^{-1}\boldsymbol{b}_{\text{acc}} + O(\Delta x^k)\boldsymbol{b}_{\text{acc}},
\end{aligned} \tag{19}
$$

*where* $\boldsymbol{I}$ *is an identity matrix. Through the definition of loss function* $\mathbb{I}_{\text{acc}}$, $\boldsymbol{b}_{\text{acc}}$ *is of order* $O(\Delta x)$*, thus*

$$
\boldsymbol{a} = \boldsymbol{M}_{\text{acc}}^{-1}\boldsymbol{b}_{\text{acc}} + O(\Delta x^{k+1}). \tag{20}
$$

*According to the* $k$*-exact property of* $\mathbb{I}_{\text{acc}}$ *[13], the reconstructed polynomials* $p_i(x)$ *by* $\boldsymbol{M}_{\text{acc}}^{-1}\boldsymbol{b}_{\text{acc}}$ *for each control volume achieve a prior* $k$*-th order of accuracy. The disturbances of small quantities* $O(\Delta x^{k+1})$ *in* $\boldsymbol{a}$ *do not affect the numerical order for the WCLS schemes.*

In this work, the penalty matrix parameter $\theta_n$ is set as

$$
\begin{cases}
\theta_n = (c_b + c_s S_n)\left|\dfrac{1}{2} - \beta_n\right| \dfrac{|p_i(x_\alpha) - p_n(x_\alpha)|}{\sum_{m\in \text{ST}_i}|\overline{u}_i - \overline{u}_m|}, x_\alpha = \Omega_i \cap \Omega_n, \\
S_{i-1} = \dfrac{1}{2}|p_i - p_{i-1}|\left(\dfrac{1}{p_{i-1}} + \dfrac{1}{p_i}\right) + \dfrac{1}{2}|\hat{u}_i - \hat{u}_{i-1}|\left(\dfrac{1}{c_{i-1}} + \dfrac{1}{c_i}\right), \\
S_{i+1} = \dfrac{1}{2}|p_i - p_{i+1}|\left(\dfrac{1}{p_{i+1}} + \dfrac{1}{p_i}\right) + \dfrac{1}{2}|\hat{u}_i - \hat{u}_{i+1}|\left(\dfrac{1}{c_{i+1}} + \dfrac{1}{c_i}\right),
\end{cases} \tag{21}
$$

where $\hat{u}$ is the normal velocity along the reconstruction lines and $c$ is the sound speed given by $c = \sqrt{\gamma p/\rho}$. $\theta_n$ is the multiplication of three components: an adaptive dimensionless penalty parameter $(c_b + c_s S_n)$, the absolute difference between $\beta_n$ and the ideal weight $\frac{1}{2}$, and a BVD term. The penalty parameter $(c_b + c_s S_n)$ is decomposed into a baseline part and an adaptive part. The coefficient $c_b$ supplies a background penalty level to maintain basic numerical stability, whereas $c_s$ scales the sensor contribution $S_n$, thereby increasing the penalty strength around nonlinear discontinuities of shocks. $\left|\frac{1}{2} - \beta_n\right|$ vanishes in smooth regions. The last BVD term, i.e., $|p_i(x_\alpha) - p_n(x_\alpha)|$ is the key to automatically distinguishing the resolved and under-resolved regions. In smooth regions, the flow is well-resolved and $p_i(x_\alpha) - p_n(x_\alpha) = O(\Delta x^{k+1})$. Since $c_b$ and $c_s$ are fixed $O(1)$ coefficients and the sensor $S_n$ remains bounded in smooth regions, the prefactor associated with the penalty remains $O(1)$. Therefore, $\theta_n = O(\Delta x^k)$ considering $\overline{u}_i - \overline{u}_m = O(\Delta x)$, $m \in \text{ST}_i$. According to Theorem 1, the prior $k$-th order of accuracy can be achieved. If the flow is under-resolved, e.g., when $\Omega_i$ lies near discontinuities, $\theta_n \to O(1)$ and dissipation is introduced to stabilize the simulation, enhancing the non-oscillatory property of the WCLS schemes. Besides the work of CLS reconstruction, the concept of BVD has also been widely applied in many other works, especially the BVD schemes proposed by Sun et al. [29] and Cheng et al. [30], where the BVD is utilized to select candidate polynomials near discontinuities.

In this work, $(c_b, c_s) = (0.1, 0.8)$ is used for the third-order WCLS scheme, while $(c_b, c_s) = (0.05, 0.8)$ is adopted for the fourth- and fifth-order WCLS schemes. Figure 2 shows the non-oscillatory profile on stationary shocks with Mach numbers 1.001, 100, and 1000, demonstrating the robustness of the present penalty formulation over a wide Mach-number range. It should be noted again that $\theta_n$ in Eq. (21) depends on the unknowns $\boldsymbol{a}$, is calculated at the beginning of each GMRES step and fixed as constant within the corresponding iteration step.

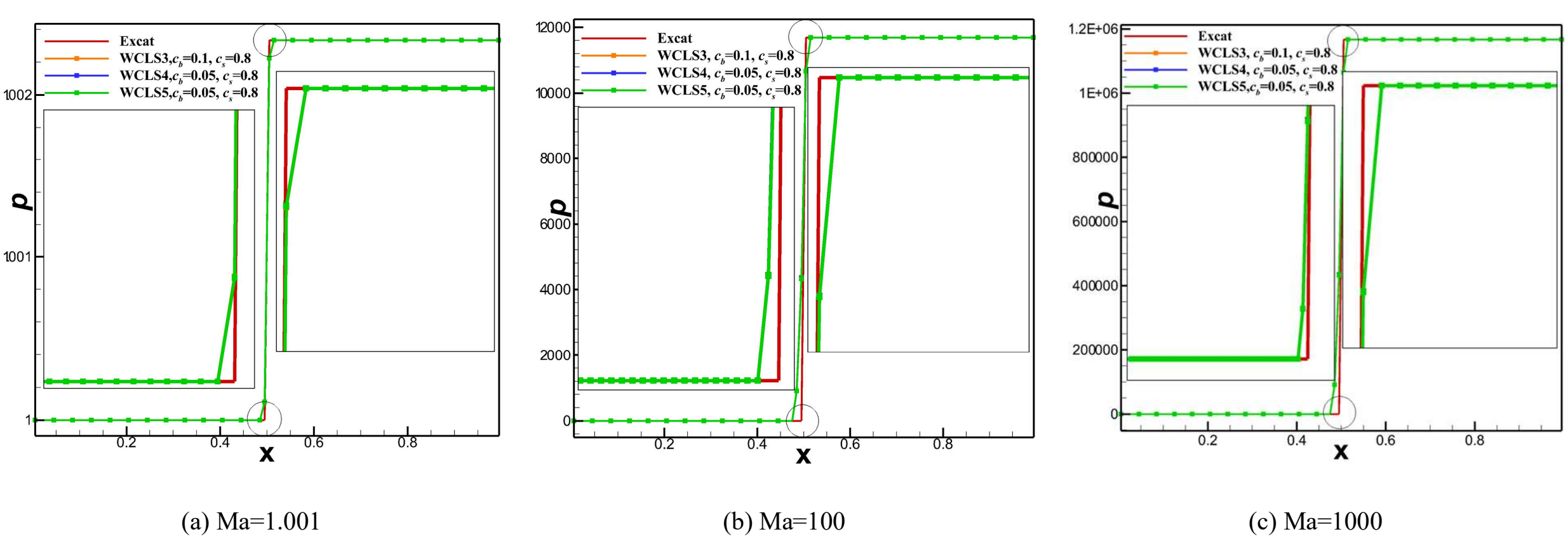

(a) Ma=1.001 (b) Ma=100 (c) Ma=1000

Figure 2: Performance of the adopted penalty coefficients for stationary shock problems at different Mach numbers.

The linear coefficients of $W_i$, $i = 1,2,\cdots,k$ do not affect the prior order of accuracy and are determined through the optimization of spectral properties. The modified non-dimensional wavenumber $\kappa'$ can be derived through the analysis in the work of Wang et al. [11] and Pan et al. [16,25] assuming a harmonic initial condition for periodic linear advection equation and is not replicated here. Once $\kappa'$ is obtained, the maximum resolved wavenumber $\kappa_{\mathrm{c,Re}}$ and $\kappa_{\mathrm{c,Im}}$ according to the real and imaginary part of $\kappa'$ are defined as in Eqs. (22) and (23), respectively,

$$\forall \kappa \leq \kappa_{\mathrm{c,Re}}, \left|\frac{\mathrm{Re}(\kappa')}{\kappa} - 1\right| \leq 0.5\%, \quad \text{and} \quad \exists \kappa > \kappa_{\mathrm{c,Re}}, \left|\frac{\mathrm{Re}(\kappa')}{\kappa} - 1\right| > 0.5\%. \tag{22}$$

$$\forall \kappa \leq \kappa_{\mathrm{c,Im}}, \left|\mathrm{Im}(\kappa')\right| \leq 0.5\%, \quad \text{and} \quad \exists \kappa > \kappa_{\mathrm{c,Im}}, \left|\mathrm{Im}(\kappa')\right| > 0.5\%. \tag{23}$$

It has been demonstrated that appropriate dissipation is essential for a robust (W)CLS scheme [16,31]. Given $\kappa_{\mathrm{c,Im}}$, a genetic algorithm is utilized to maximize $\kappa_{\mathrm{c,Re}}$. In this work, $\kappa_{\mathrm{c,Im}}$ is chosen as $1.0, 1.2, 1.2$ for the third-, fourth- and fifth-order WCLS schemes and the optimized coefficients are listed in Table 1. The dispersion and dissipation properties are presented as in Fig. 3 with $r$-index [32] defined by

$$r = \frac{\left|\frac{\partial \mathrm{Re}(\kappa')}{\partial \kappa} - 1\right| + 0.001}{|\mathrm{Im}(\kappa')| + 0.001}. \tag{24}$$

It is shown that all the optimized schemes satisfy the condition $r \leq 10$ and appropriate dissipation is applied for under-resolved flow structures.

Table 1: Optimized weights for WCLS schemes.

| Order | $\kappa_{c,\mathrm{Im}}$ | $W_1$ | $W_2$ | $W_3$ | $W_4$ |
|---|---|---|---|---|---|
| 3rd | 1.0 | 8.238e-02 | 1.072e-02 | | |
| 4th | 1.2 | 1.024e-01 | 1.172e-02 | 3.719e-03 | |
| 5th | 1.2 | 1.154e-01 | 2.249e-02 | 7.237e-03 | 2.681e-03 |

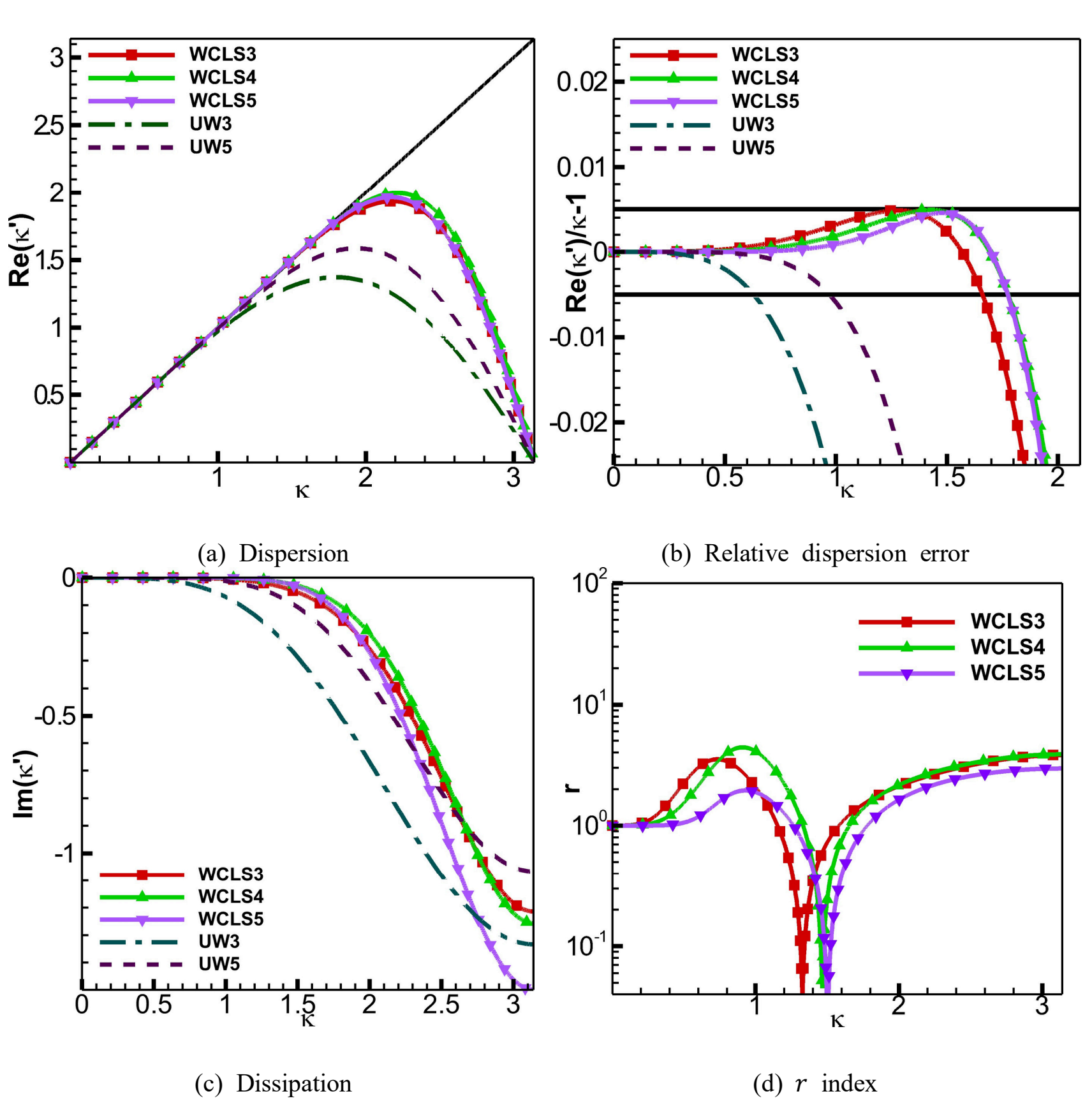

(a) Dispersion (b) Relative dispersion error

(c) Dissipation (d) $r$ index

Figure 3: Spectral properties for the optimized WCLS schemes. UW3: the third-order upwind scheme; UW5: the fifth-order upwind scheme.

Eventually, to demonstrate the effectiveness of the proposed WCLS schemes, Fig. 4 presents the reconstructed quartic polynomials in a smooth high-gradient region with profile defined by $2 + \tanh(80(x - 0.5))$ and in a discontinuous region with profile defined by a step function. Although the CLS scheme (blue lines) captures the smooth high-gradient profile as in Fig. 4 (a), it triggers high-frequency oscillations near discontinuities. To damp out the high-frequency oscillations, the CLS scheme must be coupled with a gradient limiter. Unfortunately, the utilized WBAP limiter still cannot eliminate the oscillations thoroughly and even pollutes the smooth solution as shown in Fig. 4 (a). On the contrary, the proposed WCLS scheme is oscillation-free and maintains high resolution in high-gradient smooth regions simultaneously. Figure 4 confirms that the WCLS scheme is not simply a replacement of the post-processing technique in the original CLS schemes but is a substantial improvement through the philosophy of only reconstructing along smooth reconstruction lines.

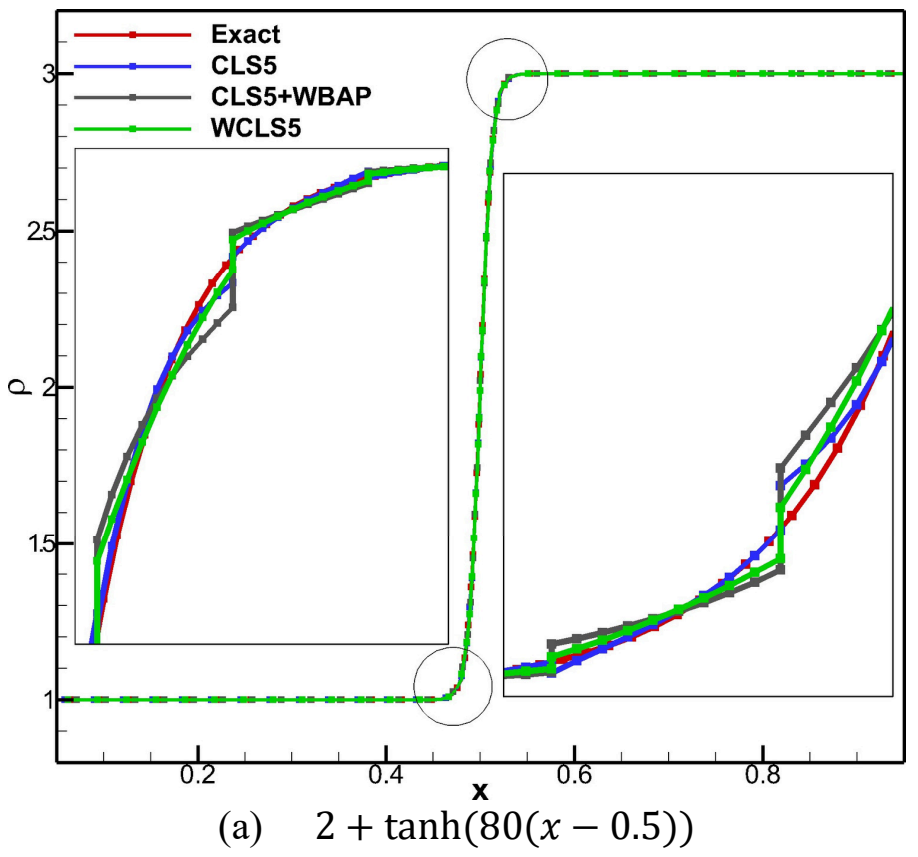

(a) $2 + \tanh(80(x - 0.5))$

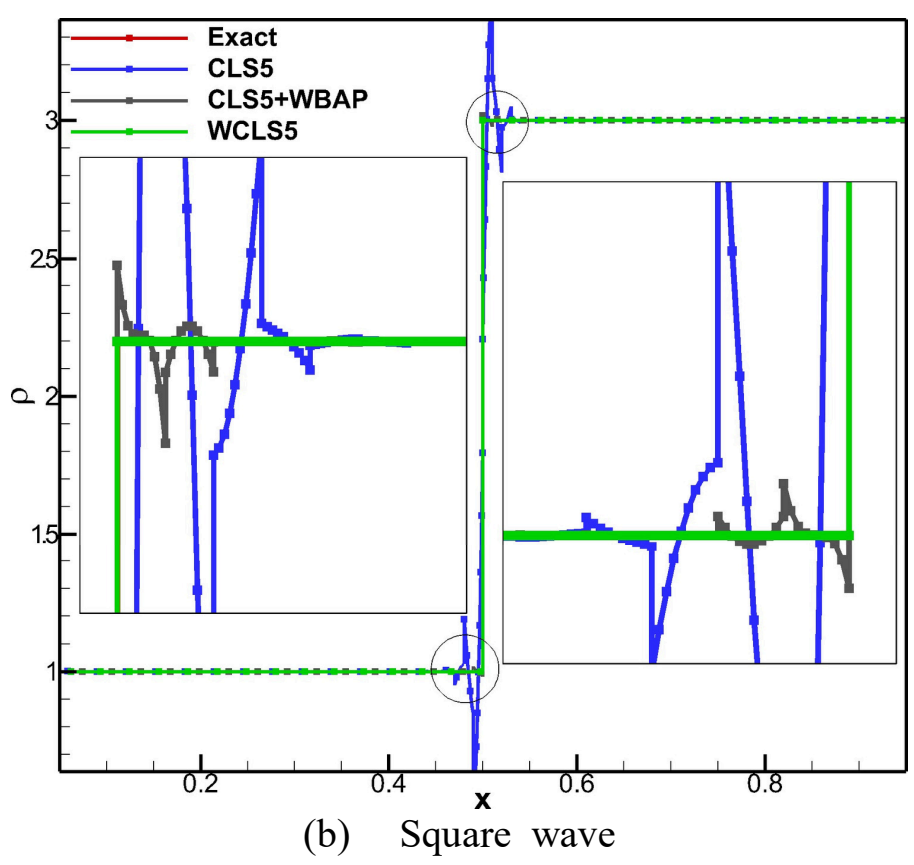

(b) Square wave

Figure 4: Reconstructed polynomials near smooth high-gradient regions and discontinuities.

## 3. Extension to Navier-Stokes Equations

As illustrated in the last section, the WCLS schemes proposed in this paper are nonlinear since the smoothness indicator $\beta$ and penalty coefficients $\theta$ depend on the solution $\boldsymbol{a}$ itself. However, when embedding the WCLS schemes into the implicit time integration, coupling the iteration of NS equations and the iteration of WCLS reconstruction [12,13,26,31], the additional solution complexity introduced by the nonlinear Eq. (11) is avoided, leading to an efficient, high-resolution WCLS framework. Before introducing the coupling strategy, the NS equations are first briefly reviewed.

The two-dimensional (2D) NS equations are

$$\frac{\partial \boldsymbol{U}}{\partial t} + \frac{\partial(\boldsymbol{F} - \boldsymbol{F}_v)}{\partial x} + \frac{\partial(\boldsymbol{G} - \boldsymbol{G}_v)}{\partial y} = 0. \tag{25}$$

In Eq. (25), $\boldsymbol{U}$ is the vector of conservation variables as

$$\boldsymbol{U} = [\rho, \rho u_1, \rho u_2, \rho E]^{\mathrm{T}}. \tag{26}$$

$\boldsymbol{F}$ and $\boldsymbol{G}$ are the inviscid components of fluxes as

$$\boldsymbol{F} = [\rho u_1, \rho u_1^2 + P, \rho u_1 u_2, \rho H u_1]^{\mathrm{T}}, \tag{27}$$

$$\boldsymbol{G} = [\rho u_2, \rho u_1 u_2, \rho u_2^2 + P, \rho H u_2]^{\mathrm{T}}. \tag{28}$$

$\boldsymbol{F}_v$ and $\boldsymbol{G}_v$ are the viscous components as

$$\boldsymbol{F}_v = \left[0, \tau_{11}, \tau_{21}, u_1\tau_{11} + u_2\tau_{12} + k_h \frac{\partial \mathcal{T}}{\partial x}\right]^{\mathrm{T}}, \tag{29}$$

$$\boldsymbol{G}_v = \left[0, \tau_{12}, \tau_{22}, u_1\tau_{21} + u_2\tau_{22} + k_h \frac{\partial \mathcal{T}}{\partial y}\right]^{\mathrm{T}}. \tag{30}$$

$\rho$, $E$, $H$ and $P$ denote the density, specific total energy, specific total enthalpy and pressure for the fluid, respectively. $u_1$ and $u_2$ denote the $x$ and $y$ components of velocity, respectively. $\mathcal{T}$ is the temperature. The specific total energy and specific total enthalpy are described as

$$E = e + \frac{1}{2}(u_1^2 + u_2^2), \tag{31}$$

$$H = e + \frac{1}{2}(u_1^2 + u_2^2) + \frac{P}{\rho}, \tag{32}$$

where $e$ is the specific internal energy. This work assumes perfect gas and $e = \frac{P}{(\gamma-1)\rho}$, where $\gamma$ is the ratio of specific heat capacity or the adiabatic index. $\gamma = 1.4$ in this work. The viscous stress tensor $\tau_{ij}$ is defined as

$$\tau_{ij} = 2\mu S_{ij} - \frac{2}{3} S_{kk}\delta_{ij}, \tag{33}$$

where $\delta_{ij}$ is the Kronecker delta function, and $S_{ij}$ is the symmetric part of velocity gradient as

$$S_{ij} = \frac{1}{2}\left(\frac{\partial u_j}{\partial x_i} + \frac{\partial u_i}{\partial x_j}\right). \tag{34}$$

$\mu$ is the viscosity coefficient. The conductivity coefficient $k_h$ is related to the viscosity $\mu$ through the Prandtl number as

$$\mathrm{Pr} = \frac{c_p \mu}{k_h}, \tag{35}$$

where $c_p$ is the specific heat capacity at constant pressure. In this work, $\mathrm{Pr} = 0.72$.

When simulating fluids with shocks, contact discontinuities or rarefaction waves, the conservation variables are first mapped into the characteristic space. Taking the control volume $\Omega_i$ as an example,

$$\overline{\boldsymbol{V}}_m = \boldsymbol{L}_i \overline{\boldsymbol{U}}_m, m \in \mathrm{ST}_i \cup \{i\}, \tag{36}$$

$$\boldsymbol{a}_{m,j}^{V} = \boldsymbol{L}_i \boldsymbol{a}_{m,j}^{U}, m \in \mathrm{ST}_i \cup \{i\}, j = 1,2,\cdots,k, \tag{37}$$

where $\boldsymbol{V}$ denotes the vector of characteristic variables and the overline $\overline{\cdot}$ denotes that the variable is cell-averaged; $\boldsymbol{L}_i$ is the left eigen matrix of the Jacobian $\partial(\boldsymbol{F} \cdot n_x + \boldsymbol{G} \cdot n_y)/\partial \boldsymbol{U}|_{\boldsymbol{U}=\overline{\boldsymbol{U}}_i}$ with $n_x$ and $n_y$ denoting the normal vector defined at the center of $\Omega_i$ along the reconstruction axes. $\boldsymbol{a}_{m,j}^{U}$ and $\boldsymbol{a}_{m,j}^{V}$ are vectors of the $j$-th unknowns for conservation and characteristic variables, respectively, as

$$\boldsymbol{a}_{m,j}^{U} = \left[a_{m,j}^{\rho}, a_{m,j}^{\rho u_1}, a_{m,j}^{\rho u_2}, a_{m,j}^{\rho E}\right]^{\mathrm{T}}, \tag{38}$$

$$\boldsymbol{a}_{m,j}^{V} = \left[a_{m,j}^{V_1}, a_{m,j}^{V_2}, a_{m,j}^{V_3}, a_{m,j}^{V_4}\right]^{\mathrm{T}}. \tag{39}$$

Then for each characteristic variable, the $\beta$ and $\theta$ parameters are calculated and the linearized reconstruction matrix for Euler/NS equations is obtained. The reconstruction matrix is written as

$$\widetilde{\boldsymbol{\mathcal{M}}}\widetilde{\boldsymbol{a}} = \widetilde{\boldsymbol{b}}, \tag{40}$$

where

$$\widetilde{\boldsymbol{\mathcal{M}}} = \begin{pmatrix} \boldsymbol{I} & \widetilde{\boldsymbol{R}}_1\widetilde{\boldsymbol{M}}_1^{(1)}\widetilde{\boldsymbol{L}}_1 & \boldsymbol{0} & \cdots & \cdots & \boldsymbol{0} & \widetilde{\boldsymbol{R}}_1\widetilde{\boldsymbol{M}}_1^{(-1)}\widetilde{\boldsymbol{L}}_1 \\ \widetilde{\boldsymbol{R}}_2\widetilde{\boldsymbol{M}}_2^{(-1)}\widetilde{\boldsymbol{L}}_2 & \boldsymbol{I} & \widetilde{\boldsymbol{R}}_2\widetilde{\boldsymbol{M}}_2^{(1)}\widetilde{\boldsymbol{L}}_2 & \boldsymbol{0} & \cdots & \cdots & \boldsymbol{0} \\ \boldsymbol{0} & \widetilde{\boldsymbol{R}}_3\widetilde{\boldsymbol{M}}_3^{(-1)}\widetilde{\boldsymbol{L}}_3 & \boldsymbol{I} & \widetilde{\boldsymbol{R}}_3\widetilde{\boldsymbol{M}}_3^{(1)}\widetilde{\boldsymbol{L}}_3 & \boldsymbol{0} & \cdots & \boldsymbol{0} \\ \vdots & \vdots & \vdots & \vdots & \vdots & \ddots & \vdots \\ \widetilde{\boldsymbol{R}}_N\widetilde{\boldsymbol{M}}_N^{(1)}\widetilde{\boldsymbol{L}}_N & \boldsymbol{0} & \cdots & \cdots & \boldsymbol{0} & \widetilde{\boldsymbol{R}}_N\widetilde{\boldsymbol{M}}_N^{(-1)}\widetilde{\boldsymbol{L}}_N & \boldsymbol{I} \end{pmatrix}, \tag{41}$$

$$\widetilde{\boldsymbol{a}} = [(\widetilde{\boldsymbol{a}}_1)^{\mathrm{T}}, (\widetilde{\boldsymbol{a}}_2)^{\mathrm{T}}, (\widetilde{\boldsymbol{a}}_3)^{\mathrm{T}}, \cdots, (\widetilde{\boldsymbol{a}}_N)^{\mathrm{T}}]^{\mathrm{T}}, \tag{42}$$

$$\widetilde{\boldsymbol{b}} = \left[\left(\widetilde{\boldsymbol{R}}_1\widetilde{\boldsymbol{M}}_1^{-1}\widetilde{\boldsymbol{b}}_1\right)^{\mathrm{T}}, \left(\widetilde{\boldsymbol{R}}_2\widetilde{\boldsymbol{M}}_2^{-1}\widetilde{\boldsymbol{b}}_2\right)^{\mathrm{T}}, \left(\widetilde{\boldsymbol{R}}_3\widetilde{\boldsymbol{M}}_3^{-1}\widetilde{\boldsymbol{b}}_3\right)^{\mathrm{T}}, \cdots, \left(\widetilde{\boldsymbol{R}}_N\widetilde{\boldsymbol{M}}_N^{-1}\widetilde{\boldsymbol{b}}_N\right)^{\mathrm{T}}\right]^{\mathrm{T}}, \tag{43}$$

and

$$\widetilde{\boldsymbol{a}}_i = \left[a_{i,1}^{\rho}, a_{i,2}^{\rho}, \cdots, a_{i,k}^{\rho}, a_{i,1}^{\rho u_1}, a_{i,2}^{\rho u_1}, \cdots, a_{i,k}^{\rho u_1}, a_{i,1}^{\rho u_2}, a_{i,2}^{\rho u_2}, \cdots, a_{i,k}^{\rho u_2}, a_{i,1}^{\rho E}, a_{i,2}^{\rho E}, \cdots, a_{i,k}^{\rho E}\right]^{\mathrm{T}}, \tag{44}$$

$$\widetilde{\boldsymbol{b}}_i = \left[b_{i,1}^{V_1}, b_{i,2}^{V_1}, \cdots, b_{i,k}^{V_1}, b_{i,1}^{V_2}, b_{i,2}^{V_2}, \cdots, b_{i,k}^{V_2}, b_{i,1}^{V_3}, b_{i,2}^{V_3}, \cdots, b_{i,k}^{V_3}, b_{i,1}^{V_4}, b_{i,2}^{V_4}, \cdots, b_{i,k}^{V_4}\right]^{\mathrm{T}}. \tag{45}$$

The matrix $\widetilde{\boldsymbol{M}}_i^{-1}$ is block diagonal defined as

$$\widetilde{\boldsymbol{M}}_i^{-1} = \begin{pmatrix} \left[\boldsymbol{M}_{V_1,i}^{(0)}\right]^{-1} & \boldsymbol{0} & \boldsymbol{0} & \boldsymbol{0} \\ \boldsymbol{0} & \left[\boldsymbol{M}_{V_2,i}^{(0)}\right]^{-1} & \boldsymbol{0} & \boldsymbol{0} \\ \boldsymbol{0} & \boldsymbol{0} & \left[\boldsymbol{M}_{V_3,i}^{(0)}\right]^{-1} & \boldsymbol{0} \\ \boldsymbol{0} & \boldsymbol{0} & \boldsymbol{0} & \left[\boldsymbol{M}_{V_4,i}^{(0)}\right]^{-1} \end{pmatrix}, \tag{46}$$

and the matrices $\widetilde{\boldsymbol{M}}_i^{(m)}$, $m = -1, 1$, are

$$\widetilde{\boldsymbol{M}}_i^{(m)} = \begin{pmatrix} \left[\boldsymbol{M}_{V_1,i}^{(0)}\right]^{-1}\boldsymbol{M}_{V_1,i}^{(m)} & \boldsymbol{0} & \boldsymbol{0} & \boldsymbol{0} \\ \boldsymbol{0} & \left[\boldsymbol{M}_{V_2,i}^{(0)}\right]^{-1}\boldsymbol{M}_{V_2,i}^{(m)} & \boldsymbol{0} & \boldsymbol{0} \\ \boldsymbol{0} & \boldsymbol{0} & \left[\boldsymbol{M}_{V_3,i}^{(0)}\right]^{-1}\boldsymbol{M}_{V_3,i}^{(m)} & \boldsymbol{0} \\ \boldsymbol{0} & \boldsymbol{0} & \boldsymbol{0} & \left[\boldsymbol{M}_{V_4,i}^{(0)}\right]^{-1}\boldsymbol{M}_{V_4,i}^{(m)} \end{pmatrix}, m = -1,1. \tag{47}$$

$\boldsymbol{M}_{V_j,i}^{(m)}$ $m = -1,0,1$ results from $\partial \mathbb{I}_i^{V_j}/\partial \boldsymbol{a}_i^{V_j}$, $j = 1,2,3,4$ assuming $\beta_n$ in Eq. (14) and $\theta_n$ in Eq.(21) independent of $\boldsymbol{a}_i^{V_j}$ in the form as

$$\boldsymbol{M}_{V_j,i}^{(-1)}\boldsymbol{a}_{i-1}^{V_j} + \boldsymbol{M}_{V_j,i}^{(0)}\boldsymbol{a}_i^{V_j} + \boldsymbol{M}_{V_j,i}^{(1)}\boldsymbol{a}_{i+1}^{V_j} = \boldsymbol{b}_i^{V_j}, \tag{48}$$

where $\boldsymbol{a}_i^{V_j} = \left[a_{i,1}^{V_j}, a_{i,2}^{V_j}, \cdots, a_{i,k}^{V_j}\right]^{\mathrm{T}}$ and $\boldsymbol{b}_i^{V_j} = \left[b_{i,1}^{V_j}, b_{i,2}^{V_j}, \cdots, b_{i,k}^{V_j}\right]^{\mathrm{T}}$.

$\widetilde{\boldsymbol{L}}_i$ and $\widetilde{\boldsymbol{R}}_i$ are extended eigen matrices arranged in a block matrix form as

$$\widetilde{\boldsymbol{L}}_i = \left[L_{i,mn}\boldsymbol{I}_{k\times k}\right]_{4\times 4}, \tag{49}$$

$$\widetilde{\boldsymbol{R}}_i = \left[R_{i,mn}\boldsymbol{I}_{k\times k}\right]_{4\times 4}. \tag{50}$$

The size of the block matrix is $4\times 4$ in 2D cases, with each block being of size $k\times k$. $k$ is the number of unknowns for polynomial of each conservation variable. $L_{i,mn}$ and $R_{i,mn}$ are the $mn$-th element of left and right eigen matrices $\boldsymbol{L}_i$ and $\boldsymbol{R}_i$, respectively. $\boldsymbol{I}_{k\times k}$ is an identity matrix of size $k\times k$.

In this work, Eq. (40) is solved by GMRES method with $N_{\mathrm{GMRES}}+1$ basis vectors. Defining $\tilde{\boldsymbol{r}} = \tilde{\boldsymbol{b}} - \widetilde{\boldsymbol{\mathcal{M}}}\tilde{\boldsymbol{a}}$, then $\tilde{\boldsymbol{a}}$ is updated to $\tilde{\boldsymbol{a}}'$ by

$$\tilde{\boldsymbol{a}}' = \tilde{\boldsymbol{a}} + \left[\tilde{\boldsymbol{r}}, \widetilde{\boldsymbol{\mathcal{M}}}\tilde{\boldsymbol{r}}, \widetilde{\boldsymbol{\mathcal{M}}}^2\tilde{\boldsymbol{r}}, \cdots, \widetilde{\boldsymbol{\mathcal{M}}}^{N_{\mathrm{GMRES}}}\tilde{\boldsymbol{r}}\right]\cdot \boldsymbol{y} \tag{51}$$

where $\boldsymbol{y}$ is chosen such that $\tilde{\boldsymbol{a}}'$ minimizes $||\tilde{\boldsymbol{b}} - \widetilde{\boldsymbol{\mathcal{M}}}\tilde{\boldsymbol{a}}'||_2$. Unless otherwise specified, $N_{\mathrm{GMRES}}$ is set to 1, which is sufficient for convergence of the NS solver.

For multi-dimensional curvilinear non-uniform grids, the WCLS reconstruction is applied in a dimension-by-dimension manner in this work. The inviscid flux is approximated by Roe Riemann solver with $h$-type entropy correction [33] to eliminate the Carbuncle phenomena. Other approximate Riemann solver such as the Osher flux with entropy fix can also be utilized [34]. The viscous flux is discretized by second-order Gauss formula as in the work of Pan [25], the derivatives along the reconstruction axes are corrected by the derivatives reconstructed by the WCLS schemes as in the work by Wang and Ren [35]. The use of derivatives reconstructed by the WCLS schemes can significantly improve the accuracy of the viscous terms, as demonstrated by Sec. 4.11.

Once the inviscid and viscous fluxes are obtained, Eq. (25) can be written in the following semi-discretized form,

$$\frac{\partial \overline{\boldsymbol{U}}_i}{\partial t} = -\frac{1}{|\Omega_i|}\sum_{f\in\partial\Omega_i} S_f\left((\widehat{\boldsymbol{F}} - \widehat{\boldsymbol{F}}_v)n_x + (\widehat{\boldsymbol{G}} - \widehat{\boldsymbol{G}}_v)n_y\right) = \boldsymbol{R}_i\left(\vec{\overline{\boldsymbol{U}}}\right), \tag{52}$$

where $|\Omega_i|$ denotes the volume of $\Omega_i$; $S_f$ is the area of face belonging to $\partial\Omega_{\mathrm{i}}$; $\widehat{\cdot}$ denotes the numerical flux evaluated at the center of face $f$.

For unsteady problems, multi-stage singly diagonal implicit Runge-Kutta (SDIRK) method is utilized to advance Eq. (52) in temporal direction. For a $s$-stage $p$-th order SDIRK method denoted as S$s$P$p$, Eq. (52) is advanced by

$$\frac{\vec{\overline{\boldsymbol{U}}}^{n+1} - \vec{\overline{\boldsymbol{U}}}^{n}}{\Delta t} = \sum_{i=1}^{s}\mathcal{B}_i\vec{\boldsymbol{R}}\left(\vec{\overline{\boldsymbol{U}}}^{(i)}\right), \tag{53}$$

where

$$\frac{\vec{\overline{\boldsymbol{U}}}^{(i)} - \vec{\overline{\boldsymbol{U}}}^{n}}{\Delta t} = \sum_{j=1}^{i}\mathcal{A}_{ij}\vec{\boldsymbol{R}}\left(\vec{\overline{\boldsymbol{U}}}^{(j)}\right), \tag{54}$$

and $\vec{\overline{\boldsymbol{U}}} = \left[\overline{\boldsymbol{U}}_1^{\mathrm{T}}, \overline{\boldsymbol{U}}_2^{\mathrm{T}}, \cdots, \overline{\boldsymbol{U}}_N^{\mathrm{T}}\right]^{\mathrm{T}}$ and $\vec{\boldsymbol{R}}\left(\vec{\overline{\boldsymbol{U}}}\right) = \left[\boldsymbol{R}_1^{\mathrm{T}}, \boldsymbol{R}_2^{\mathrm{T}}, \cdots, \boldsymbol{R}_N^{\mathrm{T}}\right]^{\mathrm{T}}$; $\mathcal{A}_{ij}$ and $\mathcal{B}_i$ are coefficients of the SDIRK method. To solve the implicit Eq. (54), a dual time stepping technique is introduced as

$$\left(\frac{1}{\Delta\tau} + \frac{1}{\Delta t} - \mathcal{A}_{ii}\frac{\partial\vec{\boldsymbol{R}}\left(\vec{\overline{\boldsymbol{U}}}^{[m]}\right)}{\partial\vec{\overline{\boldsymbol{U}}}^{[m]}}\right)\Delta\vec{\overline{\boldsymbol{U}}} = \vec{\boldsymbol{R}}_t^{[m]}, \tag{55}$$

with

$$\Delta\vec{\overline{\boldsymbol{U}}} = \vec{\overline{\boldsymbol{U}}}^{[m+1]} - \vec{\overline{\boldsymbol{U}}}^{[m]},\quad \vec{\boldsymbol{R}}_t^{[m]} = \frac{1}{\Delta t}\left(\vec{\overline{\boldsymbol{U}}}^{n} - \vec{\overline{\boldsymbol{U}}}^{[m]}\right) + \mathcal{A}_{ii}\vec{\boldsymbol{R}}\left(\vec{\overline{\boldsymbol{U}}}^{[m]}\right) + \textstyle\sum_{j=1}^{i-1}\mathcal{A}_{ij}\vec{\boldsymbol{R}}\left(\vec{\overline{\boldsymbol{U}}}^{[j]}\right). \tag{56}$$

where $[m]$ denotes the iteration index.

Equation (54) can be solved by lower-upper symmetric-Gauss-Seidel (LUSGS) algorithm with an approximate calculation of the Jacobian matrix. When evaluating $\partial\vec{\boldsymbol{R}}\left(\vec{\overline{\boldsymbol{U}}}^{[m]}\right)/\partial\vec{\overline{\boldsymbol{U}}}^{[m]}$, zero-th order reconstruction is assumed for the inviscid flux with Lax flux function. The influence of viscous flux is introduced by modifying the maximum eigen value in Lax flux function of face $f$ to

$$\lambda_{max} = \left|u_1 n_x + u_2 n_y\right| + \sqrt{\gamma\frac{p}{\rho}} + 2\max\left(\frac{4}{3}, \frac{\gamma}{\mathrm{Pr}}\right)\frac{\mu}{\rho h_f}, \tag{57}$$

where $h_f$ is the characteristic length scale on face $f$ and can be set as

$$h_f = \min(h_l, h_r). \tag{58}$$

$h_l$ and $h_r$ are the axis length for control volumes locating on the left and right sides of face $f$.

The coupled iteration algorithm for the NS equations and the WCLS reconstruction is summarized in **Algorithm I**, where an unsteady flow solver is taken as an example. In **Algorithm I**, $N_{\mathrm{RK}}$ is the number of stages for the RK method, $N_{\mathrm{inner}}$ is the maximum steps of inner iterations and $N_{\mathrm{GMRES}}+1$ is the number of vectors in Krylov sub-space. Unless otherwise specified, $N_{\mathrm{GMRES}}$ is set to 1 in the simulations. Additionally, $\epsilon_{\mathrm{NS/Euler}}$ and $\epsilon_{\mathrm{Recon}}$ are the convergence criteria for the NS/Euler equations and the WCLS reconstruction, respectively. The steady solver can be implemented in the same way as in the work of Ni et al. [21].

**Algorithm I:** Unsteady flow solver with coupled implicit NS/Euler iteration and WCLS reconstruction.

**Input:** Initial condition $\vec{\overline{\boldsymbol{U}}}^{n}$.

**Output:** Solution at next time step $\vec{\overline{\boldsymbol{U}}}^{n+1}$.

1 $\vec{\overline{\boldsymbol{U}}}^{(0)} = \vec{\overline{\boldsymbol{U}}}^{n}$ //initialization
2 **Do** $i = 1, 2, \cdots N_{\mathrm{RK}}$
3 $\quad \vec{\overline{\boldsymbol{U}}}^{[1]} = \vec{\overline{\boldsymbol{U}}}^{(i-1)}$
4 $\quad \widetilde{\boldsymbol{a}}^{[1]} = \widetilde{\boldsymbol{a}}^{(i-1)}$
5 $\quad$ **Do** $s = 1$ **to** $N_{\mathrm{inner}}$
6 $\quad\quad$ Calculate $\theta^{[s]}$ and $\beta^{[s]}$; //Smoothness indicators and coefficients of penalty matrices
7 $\quad\quad \boldsymbol{r}_0 = \widetilde{\boldsymbol{b}}^{[s]} - \widetilde{\boldsymbol{\mathcal{M}}}^{[s]}\widetilde{\boldsymbol{a}}^{[s]}$
8 $\quad\quad \boldsymbol{Z}_1 = \boldsymbol{r}_0/\|\boldsymbol{r}_0\|$
9 $\quad\quad$ **Do** $m = 1, N_{\mathrm{GMRES}}$
10 $\quad\quad\quad \boldsymbol{Z}_{m+1} = \widetilde{\boldsymbol{\mathcal{M}}}^{[s]}\boldsymbol{Z}_m$
11 $\quad\quad\quad$ **Do** $j = 1$ **to** $m$
12 $\quad\quad\quad\quad \boldsymbol{H}(j,m) = (\boldsymbol{Z}_{m+1}, \boldsymbol{Z}_j)$
13 $\quad\quad\quad\quad \boldsymbol{Z}_{m+1} = \boldsymbol{Z}_{m+1} - \boldsymbol{H}(j,m)\boldsymbol{Z}_j$
14 $\quad\quad\quad$ **End**
15 $\quad\quad\quad \boldsymbol{Z}_{m+1} \leftarrow \boldsymbol{Z}_{m+1}/\|\boldsymbol{Z}_{m+1}\|$
16 $\quad\quad$ **End**
17 $\quad\quad \widetilde{\boldsymbol{a}}^{[s+1]} = \widetilde{\boldsymbol{a}}^{[s]} + [\boldsymbol{Z}_1, \boldsymbol{Z}_2, \boldsymbol{Z}_3, \cdots, \boldsymbol{Z}_{N_{\mathrm{GMRES}}+1}] \cdot \boldsymbol{y}$
18 $\quad\quad$ where $\boldsymbol{y} = \mathrm{argmin}_{\boldsymbol{y}} \|\boldsymbol{r}_0 - \widetilde{\boldsymbol{\mathcal{M}}}^{[s]} \cdot ([\boldsymbol{Z}_1, \boldsymbol{Z}_2, \boldsymbol{Z}_3, \cdots, \boldsymbol{Z}_{N_{\mathrm{GMRES}}+1}] \cdot \boldsymbol{y})\|_2$
19 $\quad\quad$ Solve $\left(\frac{1}{\Delta t} + \frac{1}{\Delta t} - \mathcal{A}_{ij}\,\partial \vec{\boldsymbol{\mathcal{R}}}\left(\vec{\overline{\boldsymbol{U}}}^{[s]}\right)/\partial \vec{\overline{\boldsymbol{U}}}^{[s]}\right)\Delta\vec{\overline{\boldsymbol{U}}} = \vec{\boldsymbol{R}}_t^{[s]}$
20 $\quad\quad \vec{\overline{\boldsymbol{U}}}^{[s+1]} = \vec{\overline{\boldsymbol{U}}}^{[s]} + \Delta\vec{\overline{\boldsymbol{U}}}$; //update the cell-averages
21 $\quad\quad$ **if** $\left\|\vec{\boldsymbol{R}}_t^{[s]}\right\|_2 / \left\|\sum_{j=1}^{i} \mathcal{A}_{ij}\vec{\boldsymbol{R}}\left(\vec{\overline{\boldsymbol{U}}}^{(j)}\right)\right\|_2 < \epsilon_{\mathrm{NS/Euler}}$ and $\|\boldsymbol{r}_0\|_2/\|\widetilde{\boldsymbol{b}}^{[s]}\|_2 < \epsilon_{\mathrm{Recon}}$ **then**
22 $\quad\quad\quad$ break
23 $\quad\quad$ **end**
24 $\quad$ **End**
25 $\quad \vec{\overline{\boldsymbol{U}}}^{(i)} = \vec{\overline{\boldsymbol{U}}}^{[s+1]}$ //the converged solution of the inner iterations
26 **End**
27 $\vec{\overline{\boldsymbol{U}}}^{n+1} = \vec{\overline{\boldsymbol{U}}}^{n} + \Delta t\left(\sum_{i=1}^{N_{\mathrm{RK}}} \mathfrak{B}_i \vec{\boldsymbol{\mathcal{R}}}\left(\vec{\overline{\boldsymbol{U}}}^{(i)}\right)\right)$.

## 4. Numerical Experiments

A series of inviscid Euler and viscous NS problems are utilized in this section to validate the high resolution and shock-capturing capability of the proposed WCLS schemes. The discretization method in temporal direction is the four-stage third-order explicit singly diagonal implicit Runge-Kutta (ESDIRKS4P3) scheme [24] and $N_{\mathrm{GMRES}} = 1$ unless otherwise specified. For viscous boundaries, the weight of boundary values are set a $W_{\mathrm{bc},0} = 100$ and $W_{\mathrm{bc},j} = 0, j = 1,2,\cdots,k$ in Eq. (9) to enforce a non-slip boundary condition; For inviscid boundaries, the weight of boundary values are set as $W_{\mathrm{bc},0} = 0.0$ and $W_{\mathrm{bc},j} = 0.5$, $j = 1,2,\cdots,k$ with constraints of derivatives $\mathrm{d}^j p/\mathrm{d}x^j = 0$, $j = 1,2,\cdots,k$. The critical value of shock detector is chosen as $\sigma_c = 0.6$. The WENO3-JS scheme is a version for uniform grids, and the WENO5-Z scheme is a version modified for non-uniform grids [36]. The CLS schemes are all coupled with WBAP limiter [11–15] (denoted as WBAP-CLS) with a third-party smoothness indicator proposed by Pan et al.

[16].

4.1 Numerical Accuracy

In this section, 1D Euler equations with periodic boundary conditions are utilized to check the accuracy of WCLS schemes. The computational domain is $[0,2]$ discretized by uniform grids and the initial condition is

$$[\rho, u, P] = [1.1 + \sin^m(\pi x), 1.0, 1.0], \tag{59}$$

where $m$ is a constant chosen as 1.0, 2.0 and 3.0.

There exist first-order extrema when $m = 1.0$ or 2.0 and second-order extrema when $m = 3.0$. The simulation is advanced until $t = 2.0$. The convergence criteria for the Euler equations and WCLS reconstruction are set as $\epsilon_{\text{NS}} = 1 \times 10^{-11}$ and $\epsilon_{\text{Recon}} = 1 \times 10^{-1}$, respectively. $N_{\text{GMRES}} = 3$ is used in this case to ensure fully convergent results. The $L_2$ and $L_\infty$ norm of errors for the density are presented in Fig. 5. The results of WENO3-JS scheme [20] and WENO5-Z scheme [17] are also presented in comparison. To further validate the accuracy-preserving property of the WCLS schemes, the weighting process is turned on across the computational domain without shock detector in this example. For WCLS3 scheme, ESDIRKS4P3 [24] is adopted in temporal discretization; for WCLS4 and WCLS5 schemes, ESDIRKS6P5 [24] schemes are utilized to ensure sufficient temporal accuracy.

The prior third, fourth and fifth orders are all preserved for the WCLS3, WCLS4 and WCLS5 schemes, respectively, from $m = 1$ to 3. Besides, the errors of the WCLS schemes are orders of magnitude smaller than those of the traditional WENO schemes with the same order of accuracy. On the other hand, due to the inability to distinguish smooth extrema and discontinuities, there exists order degradation for the WENO3-JS scheme whenever $m = 1$, 2 or 3. For the improved WENO5-Z scheme, when second-order extreme exists, i.e., when $m = 3$, order degradation is observed.

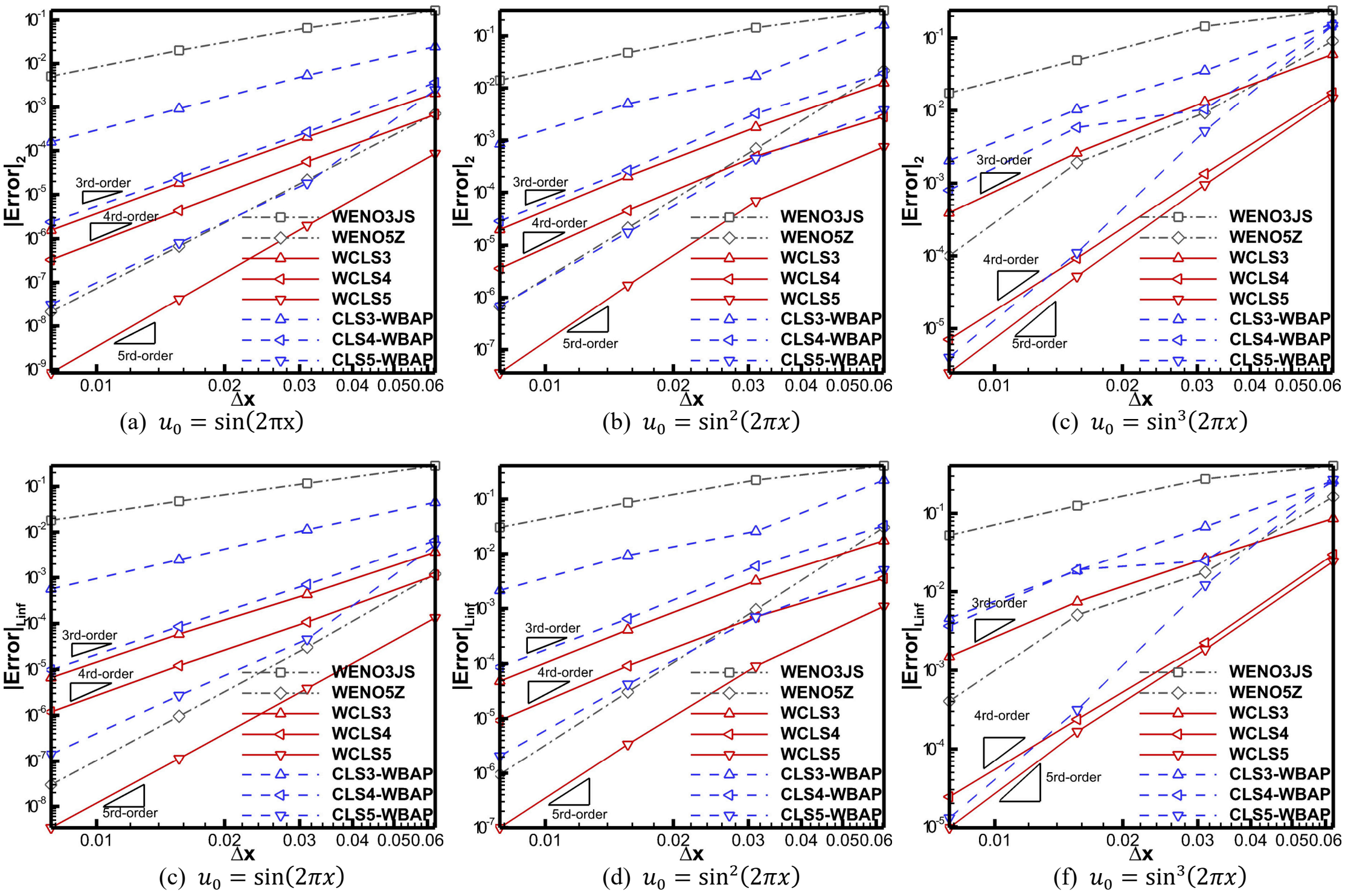

(a) $u_0 = \sin(2\pi x)$ (b) $u_0 = \sin^2(2\pi x)$ (c) $u_0 = \sin^3(2\pi x)$

(c) $u_0 = \sin(2\pi x)$ (d) $u_0 = \sin^2(2\pi x)$ (f) $u_0 = \sin^3(2\pi x)$

Figure 5: Accuracy test for the WCLS schemes.

In Fig. 5, the errors of the WBAP-CLS schemes are also presented, demonstrating a significant improvement of the proposed WCLS schemes over the WBAP-CLS schemes.

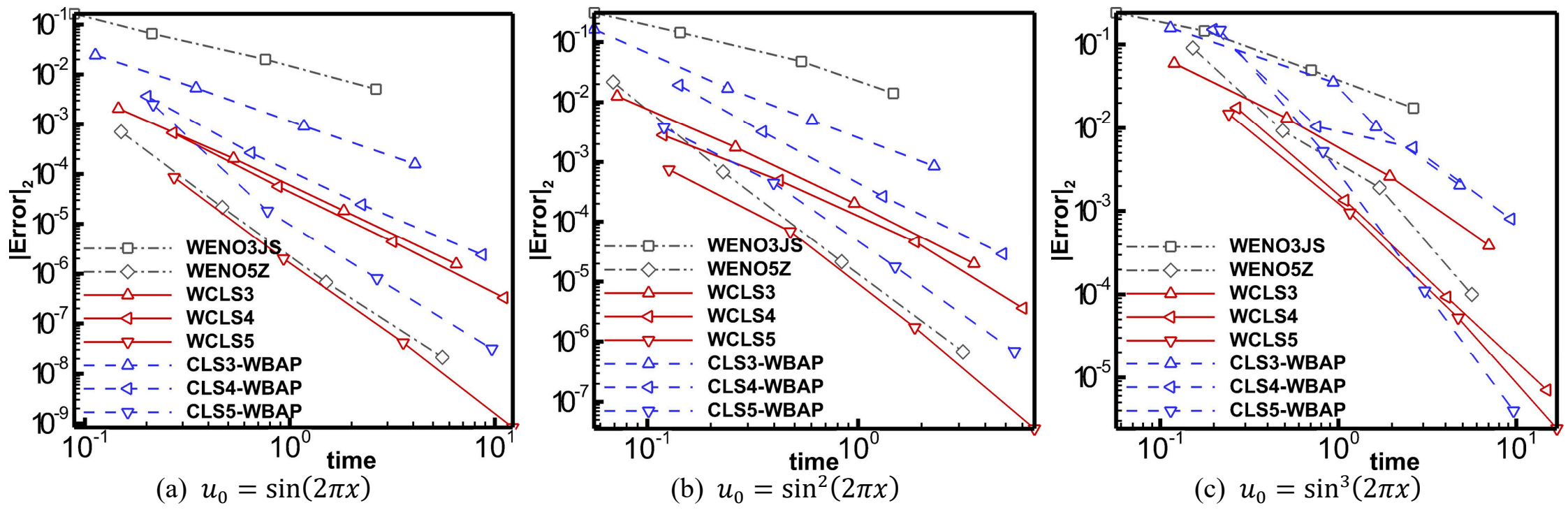

(a) $u_0 = \sin(2\pi x)$ (b) $u_0 = \sin^2(2\pi x)$ (c) $u_0 = \sin^3(2\pi x)$

Figure 6: Computational efficiency for the WCLS schemes.

Figure 6 also presents the simulation cost versus numerical errors. The proposed WCLS schemes consume less simulation time when achieving the same level of errors, confirming again their high efficiency when embedded with implicit time discretization framework.

## 4.2 Gaussian-Square-Triangle-Ellipse Waves

After confirming the formal accuracy of the WCLS schemes, we next examine their discontinuity-capturing performance in a problem containing coexisting Gaussian, square, triangle, and ellipse waves. The computational domain is [0, 2], discretized by 200 uniform control volumes with periodic boundary conditions. The flow velocity and pressure are both set to 1.0, and the density is initialized as

$$\rho(x) = 1 + \begin{cases} \frac{1}{6}\big(G(x,\psi,z-\delta) + G(x,\psi,z+\delta) + 4G(x,\psi,z)\big), & 0.2 \le x \le 0.4, \\ 1, & 0.6 \le x \le 0.8, \\ 1 - |10(x-1.1)|, & 1.0 \le x \le 1.2, \\ \frac{1}{6}\big(F(x,\alpha,a-\delta) + F(x,\alpha,a+\delta) + 4F(x,\alpha,a)\big), & 1.4 \le x \le 1.6, \\ 0.001, & \text{elsewise}, \end{cases} \tag{60}$$

where

$$G(x,\psi,z) = e^{\psi(x-1-z)^2}, F(x,\alpha,a) = \sqrt{\max(1-\alpha^2(x-1-a^2),0)}, \tag{61}$$

and $z = -0.7$, $\delta = 0.005$, $\psi = \log_{10} 2/(36\delta^2)$, $\alpha = 10$ and $a = 0.5$. The convergence criteria of $\epsilon_{\text{Euler}}$ and $\epsilon_{\text{Recon}}$ are both set as $1.0 \times 10^{-4}$ and the maximum number of inner iterations is set as 6. The Courant number (CFL) is set as 1.0 and the simulation is advanced until $t = 10.0$.

The results are presented in Fig. 7 and the close views of the four waves are illustrated in Fig. 8. As shown, all four waves are dissipated by the WENO3-JS scheme. The WENO5-Z scheme can describe the Gaussian, triangle and ellipse waves in high resolution. Nevertheless, the square wave is almost dissipated by the WENO5-Z scheme and there is no constant region between $x \in [0.65,0.75]$. The results of the WCLS schemes perform better than the WENO5-Z scheme. It should be mentioned that the square wave can still be well-described by the WCLS schemes, even by the third-order WCLS3 scheme. This example confirms not only the high accuracy in smooth regions but also the high resolution near contact discontinuities of the proposed WCLS schemes. The results of the WBAP-CLS schemes are also presented. The proposed WCLS schemes show much better resolution than the WBAP-CLS schemes near square and ellipse waves.

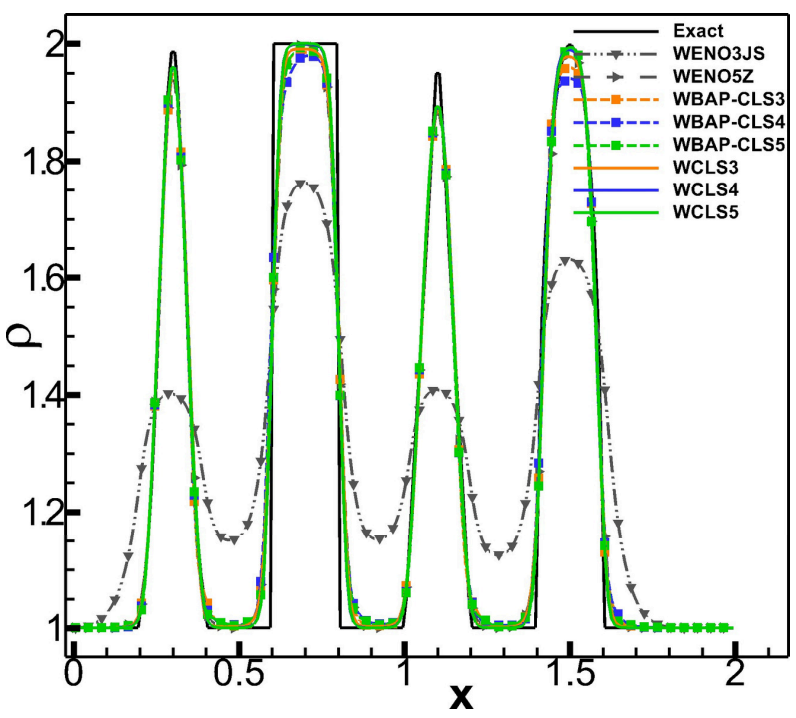

Figure 7: Results of the WCLS schemes for Gaussian-square-triangle-ellipse waves.

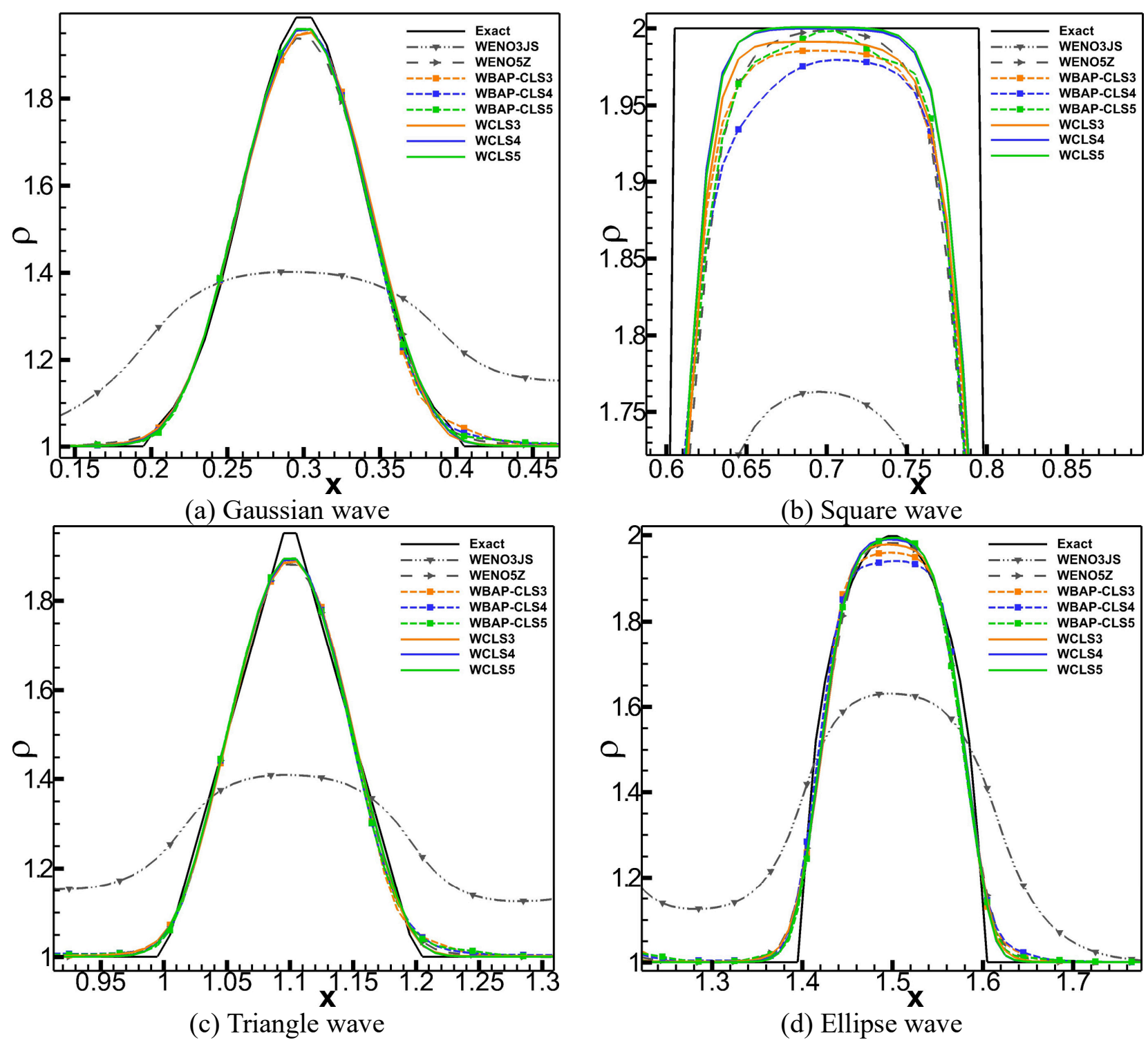

Figure 8: Close views of the Gaussian-square-triangle-ellipse waves.

### 4.3 Propagation of Broadband Sound Waves

In this section, broadband sound-wave propagation is used to assess the optimized spectral properties of the proposed WCLS schemes. The initial condition is

$$\begin{cases} P(x,t=0) = P_0\left(1+\varepsilon\sum_{k=1}^{N/2}\left[E_p(k)\right]^{0.5}\sin\left(2\pi k(x+\psi_k)\right)\right), \\ \rho(x,t=0) = \rho_0\left(\dfrac{P(x,t=0)}{P_0}\right)^{1/\gamma}, \\ u(x,t=0) = u_0+\dfrac{2}{\gamma-1}\left(c(x,t=0)-c_0\right). \end{cases} \tag{62}$$

The energy spectrum is $E_p(k) = (k/k_0)^4 e^{-2(k/k_0)^2}$, where the characteristic wave number $k_0$ is chosen from $\{4, 8, 12\}$. The mean flow is $P_0 = 1$ , $u_0 = 1$ and $\rho_0 = 1$. $\varepsilon$ denotes the intensity of acoustic energy and is set as 0.001. $\psi_k \in [0, 1]$, $k = 1, 2, 3, \cdots, N/2$ are random numbers. $c$ is the speed of sound defined by $c = \sqrt{\gamma P/\rho}$. The periodic simulation domain locates in [0, 1] and discretized by $N$ = 128 uniform control volumes. The simulation parameters are $\epsilon_{\text{Euler}} = \epsilon_{\text{Recon}} = 1 \times 10^{-4}$, $N_{\text{inner}} = 6$ and CFL = 1.0.

The pressure profile at $t = 1.0/(1+\sqrt{\gamma})$ is shown in Fig. 9. When $k_0 = 4$, all the schemes are successful in resolving the flow structures. When $k_0 = 8$ , numerical dissipation takes effect. WENO3-JS scheme is the most dissipative one. The WCLS schemes perform better than the WENO3-JS, WENO5-Z and WBAP-CLS schemes. When $k_0$ is increased to 12. The flow structures are further dissipated by the WENO3-JS scheme. However, as shown in Fig. 9 (d), the proposed three orders of WCLS schemes resolve finer structures than the WENO5-Z scheme and the WBAP-CLS schemes.

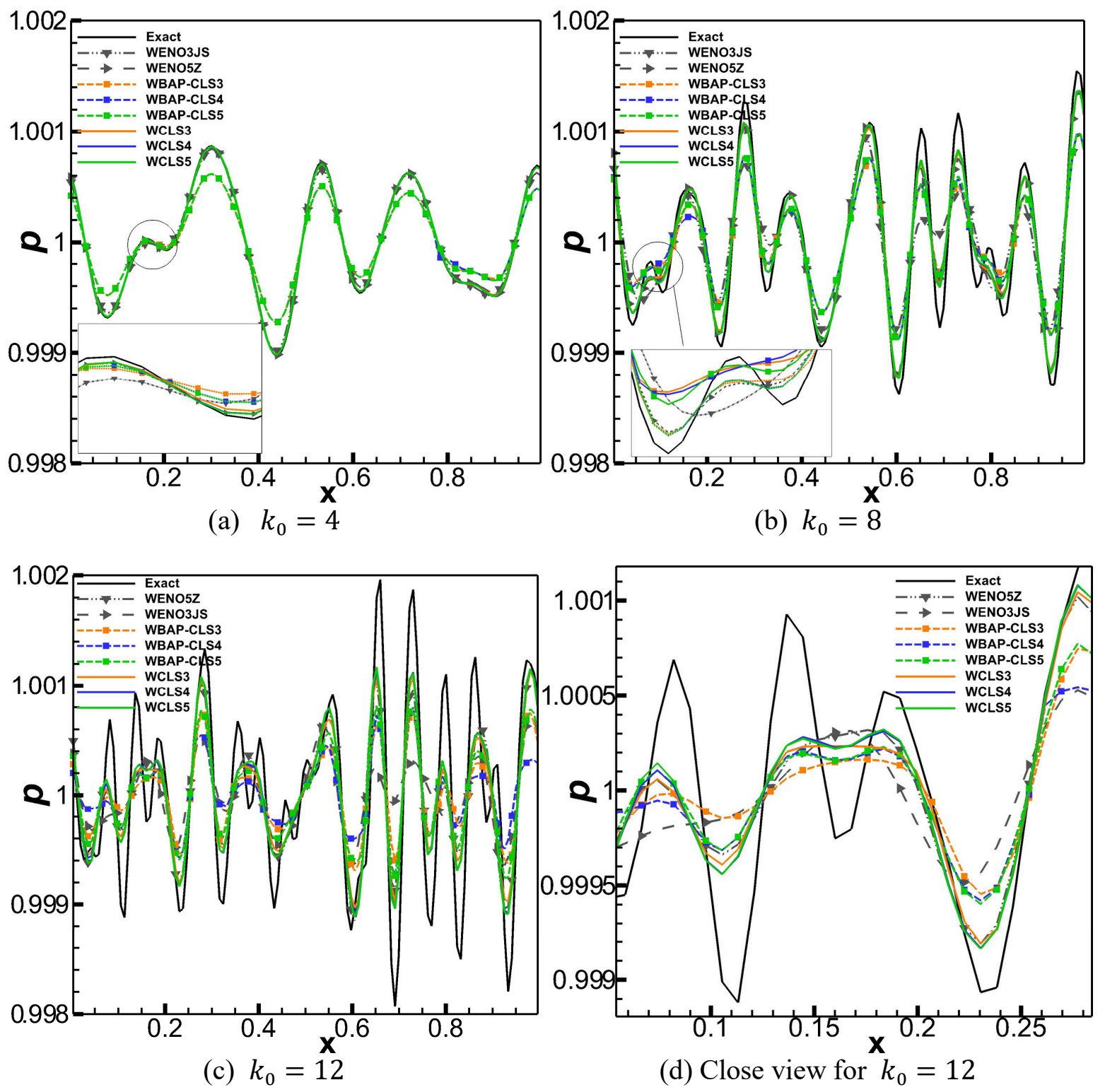

(a) $k_0 = 4$ (b) $k_0 = 8$

(c) $k_0 = 12$ (d) Close view for $k_0 = 12$

Figure 9: Acoustic waves for different $k_0$.

## 4.4 Lax Problem

In this section, a shock tube problem, denoted as Lax problem [37], is utilized to further validate the robustness of the proposed WCLS schemes near shocks. The simulation domain is [0,1] discretized by 200 uniform control volumes. $\epsilon_{\text{Euler}}$ and $\epsilon_{\text{Recon}}$ are both set as $1.0 \times 10^{-4}$. The maximum number of inner iterations is 6. The Courant number is set as 1.0 and the simulation is advanced until $t = 0.1$.

The results are shown in Fig. 10. No oscillations are observed near the contact and shock discontinuities. The WCLS schemes outperform the WENO3-JS and WENO5-Z schemes in resolving these discontinuities. Compared with the CLS schemes, the WCLS schemes provide both higher resolution and fewer oscillations.

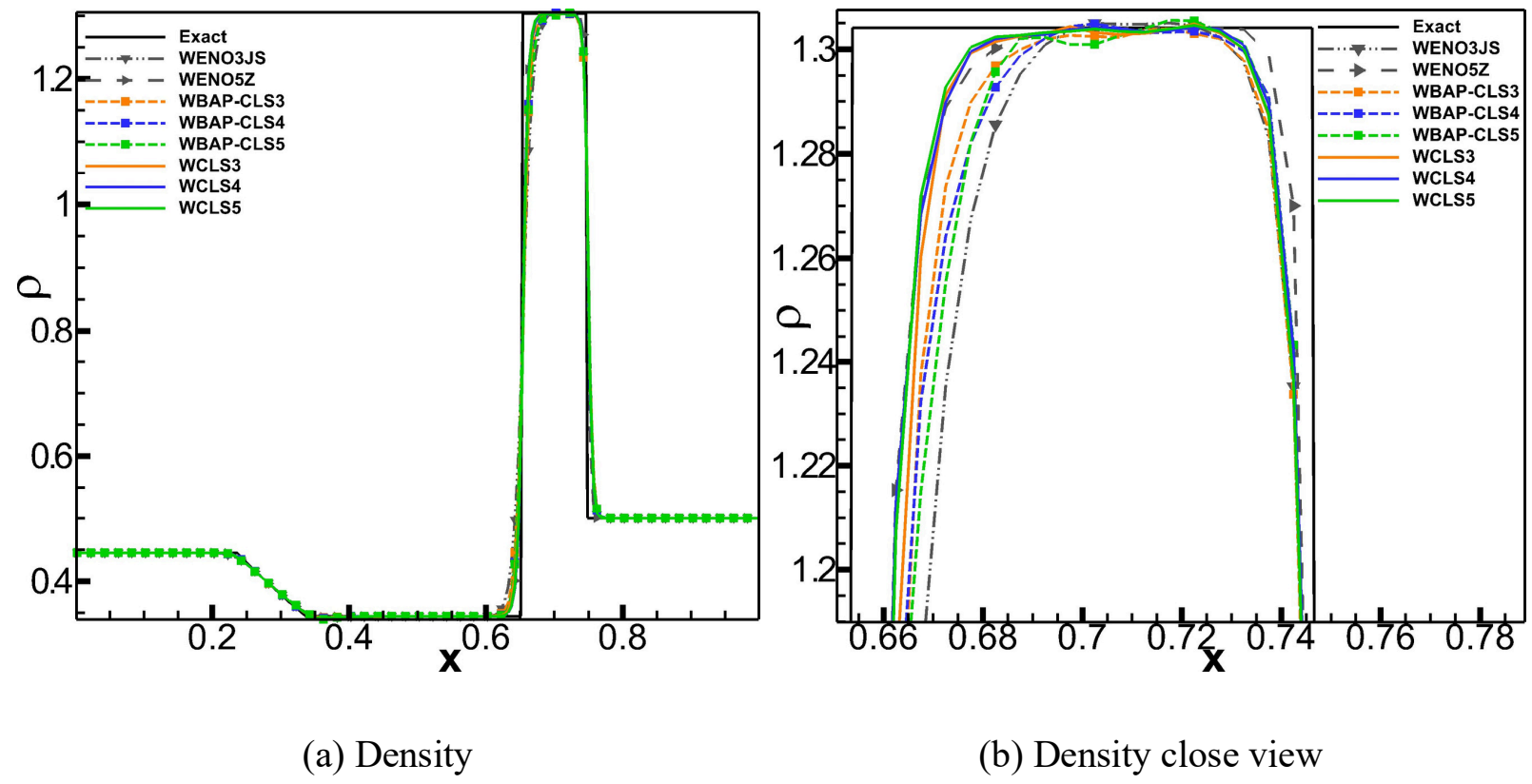

(a) Density (b) Density close view

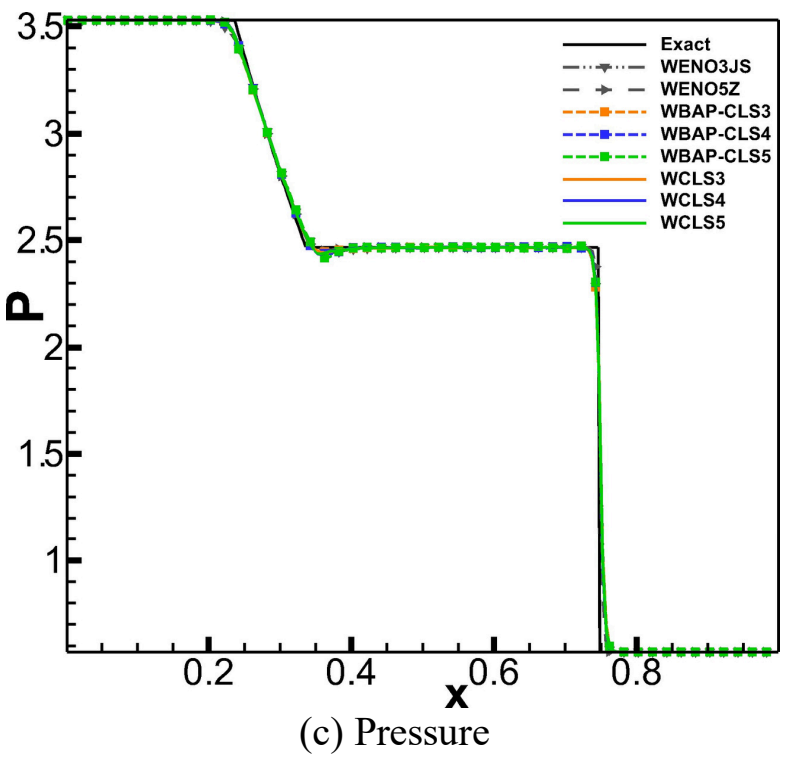

(c) Pressure

Figure 10: Lax shock tube problem. $N = 200$, $t = 0.1$ and CFL = 1.0.

## 4.5 Shu-Osher Problem [38]

The Shu-Osher problem is a standard benchmark for assessing both the resolution and shock-capturing capability of reconstruction methods for compressible flows. In this problem, a right-moving shock interacts with density perturbations and generates multiscale flow structures. The initial condition for the Shu-Osher problem is

$$[\rho, u, P] = \begin{cases} 3.857143, 2.629369, 10.333333, & \text{if } 0.0 \le x < 1.0, \\ 1 + 0.2\sin(5x), 0, 1, & \text{if } 1.0 \le x \le 10. \end{cases} \tag{63}$$

The simulation domain is [0,10] discretized by 200 uniform control volumes. The convergence criteria for both Euler equations and the WCLS reconstruction are $1 \times 10^{-4}$ and the maximum inner iteration step is 6. The Courant number is 1.0. The density profile at $t$ =1.8 is shown in Fig. 11. Near shocks, the proposed WCLS schemes capture the discontinuities sharply without oscillations. In smooth regions, the WCLS schemes simulate the flow structures accurately in higher resolution than the WENO5-Z and WBAP-CLS schemes.

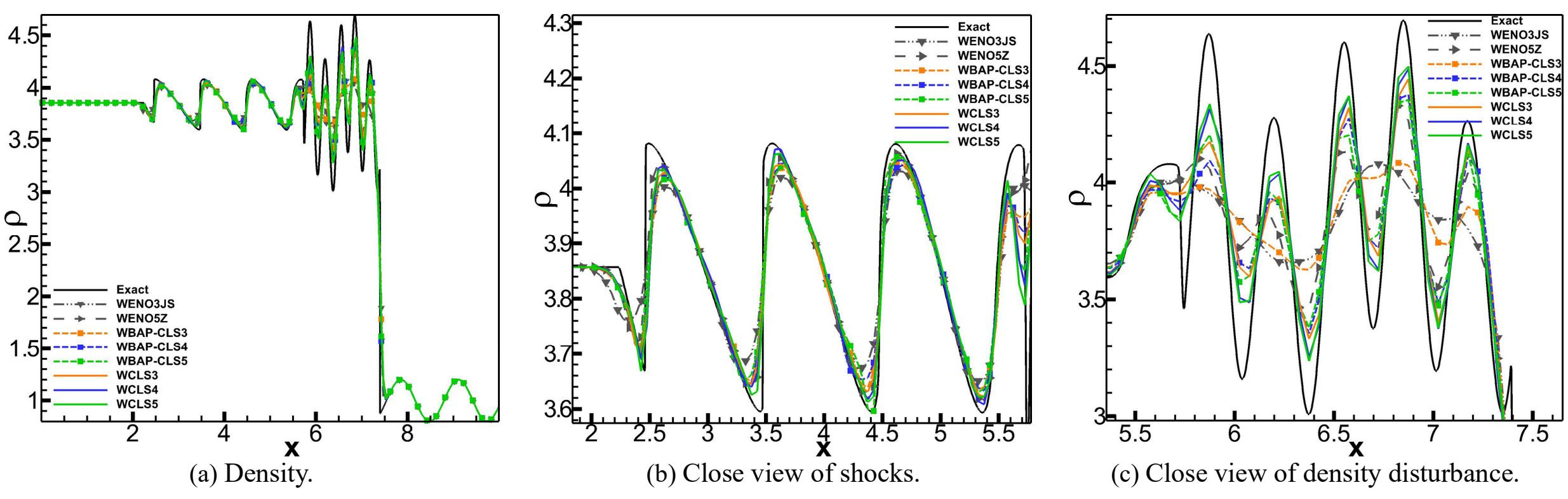

(a) Density. (b) Close view of shocks. (c) Close view of density disturbance.

Figure 11: Shu-Osher problem. $N = 200$, $t = 1.8$ and CFL = 1.0.

## 4.6 Isentropic Vortex Problem on Non-uniform Grids

A smooth isentropic vortex is first used to assess the performance of the WCLS schemes on two-dimensional non-uniform curvilinear structured grids. An isentropic disturbance is superimposed on the mean flow, and the initial condition is given by

$$\begin{cases} u_1 = 1 - \dfrac{5}{2\pi} e^{0.5(1-r^2)} y, \\ u_2 = 1 + \dfrac{5}{2\pi} e^{0.5(1-r^2)} x, \\ T = 1 - \dfrac{25(\gamma - 1)}{8\gamma\pi^2} e^{1-r^2}, \\ p = \rho^{\gamma}, \end{cases} \tag{64}$$

where $r = x^2 + y^2$. The computational domain is $[-5,5] \times [-5,5]$ with sin-wave disturbance. The non-uniform grids are constructed following Sec. 4.3.1 of Wang and Ren [35], and are shown in Fig. 12. The grid number is $128 \times 128$. The CFL number is 2.0. The convergence criteria for both Euler equations and the WCLS reconstruction are $1 \times 10^{-6}$ and the maximum inner iteration step is 6.

Density contours at time $t = 10$ are shown as in Fig. 13. The WENO3-JS scheme is implemented assuming uniform grids, failing to preserve the structure of contours which should be circles. On the other hand, the proposed WCLS schemes and the WENO5-Z schemes for non-uniform grids are all successful in capturing the flow structures with high resolution.

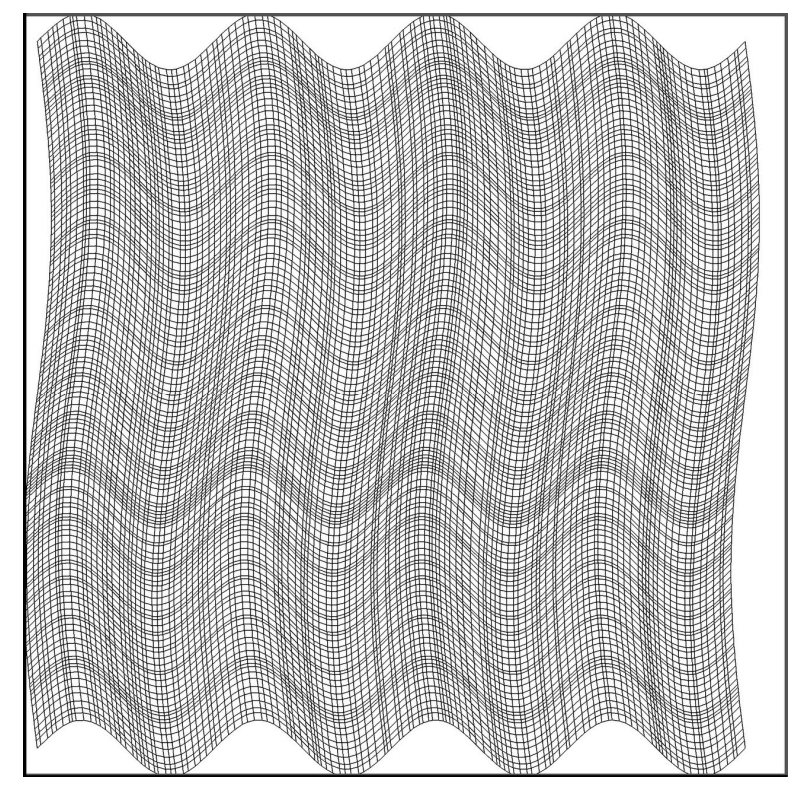

Figure 12: Distorted-rand grids for isentropic vortex problem.

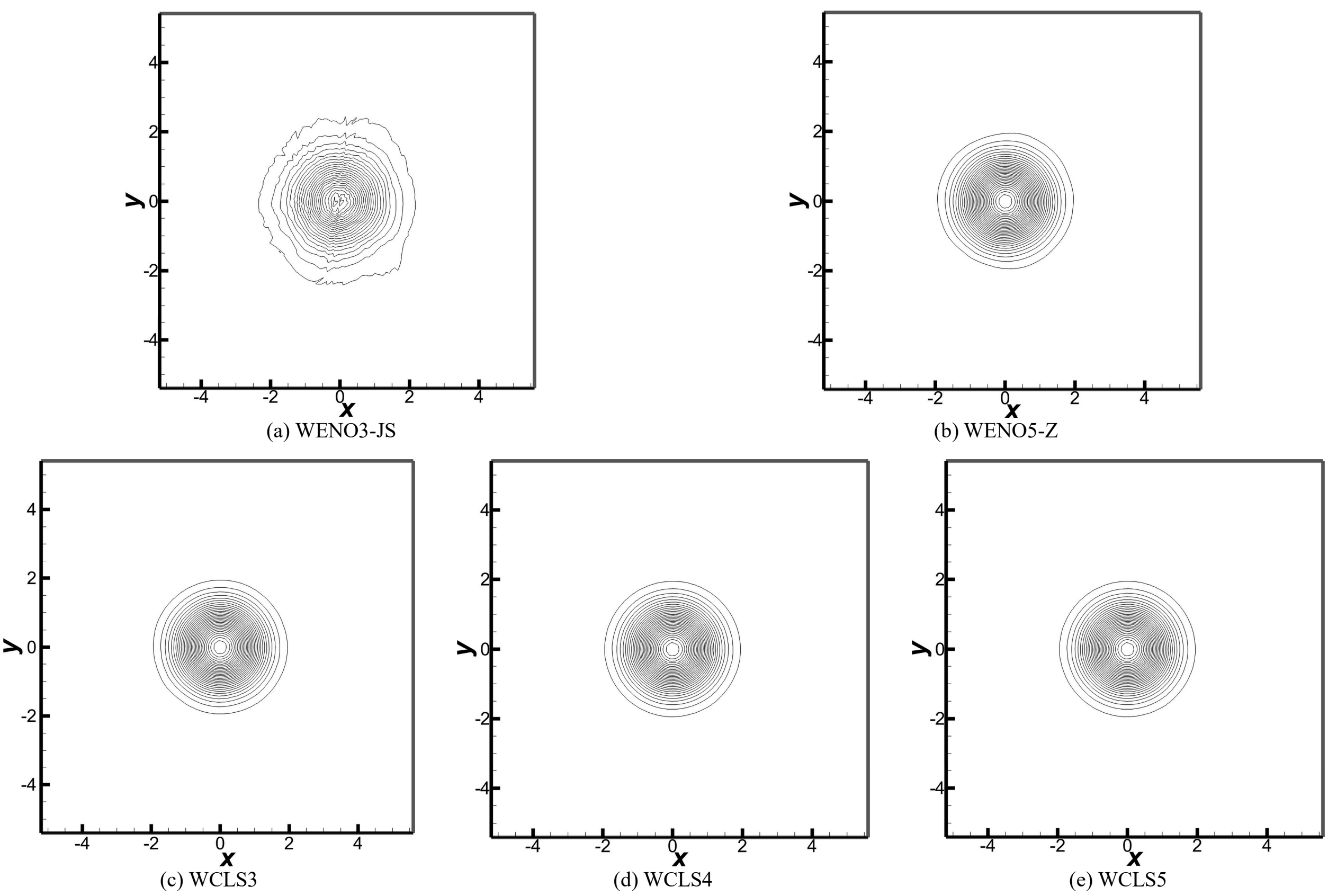

(a) WENO3-JS (b) WENO5-Z

(c) WCLS3 (d) WCLS4 (e) WCLS5

Figure 13: Density contours on distorted-rand grids at t = 10. 30 equally spaced contour lines from 0.5085 to 0.9860.

### 4.7 Shock Vortex Interaction on Non-Uniform Grids [20]

To further validate the proposed WCLS schemes on non-uniform grids, a shock vortex interaction problem is simulated in this section. A stationary shock with Mach number 1.1 is positioned at $x = 0.5$. Gas property before the shock is $[\rho, u_1, u_2, P] = [1.0, 1.1\sqrt{\gamma}, 0.0, 1.0]$. Initially, an isentropic vortex centered at $[x_c, y_c] = [0.25, 0.5]$ is added to the mean flow as

$$\begin{cases} \delta u_1 = -\dfrac{\varepsilon}{r_c} e^{\vartheta\left(1-\frac{r^2}{r_c^2}\right)}(y - y_c), \\ \delta u_2 = \dfrac{\varepsilon}{r_c} e^{\vartheta\left(1-\frac{r^2}{r_c^2}\right)}(x - x_c), \\ \delta T = -\dfrac{\varepsilon^2(\gamma-1)}{4\vartheta\gamma} e^{2\vartheta\left(1-\frac{r^2}{r_c^2}\right)}, \\ \delta S = 0, \end{cases} \tag{65}$$

where $\varepsilon = 0.3$, $\vartheta = 0.204$, $r_c = 0.05$, $S = p/\rho^\gamma$ and $r = \sqrt{(x-x_c)^2 + (y-y_c)^2}$. The computational domain is $[0,2] \times [0,1]$ discretized by $250 \times 125$ grids. The grids in $x$-direction are stretched with ratio 4.0 and concentrated around $x = 0.5$, while the grids in $y$-direction are disturbed by random distance as the ALT-RAND form in work of Wang [35]. The grids are illustrated as in Fig. 14. Inflow and outflow boundary conditions are applied on the left and right boundaries, respectively. Inviscid walls are applied on the top and bottom boundaries. The convergence criteria for both Euler equations and the WCLS reconstruction are $1 \times 10^{-6}$ and the maximum inner iteration step is 6.

The results at time $t = 0.6$ with CFL number 1.0 are shown in Fig. 15. As illustrated, the proposed WCLS schemes are successful in capturing the shocks on non-uniform grids robustly. After transmitting through shocks, the vortex is accurately captured. On the other hand, the uniform version of WENO3-JS scheme fails to resolve the smooth vortex after the shock.

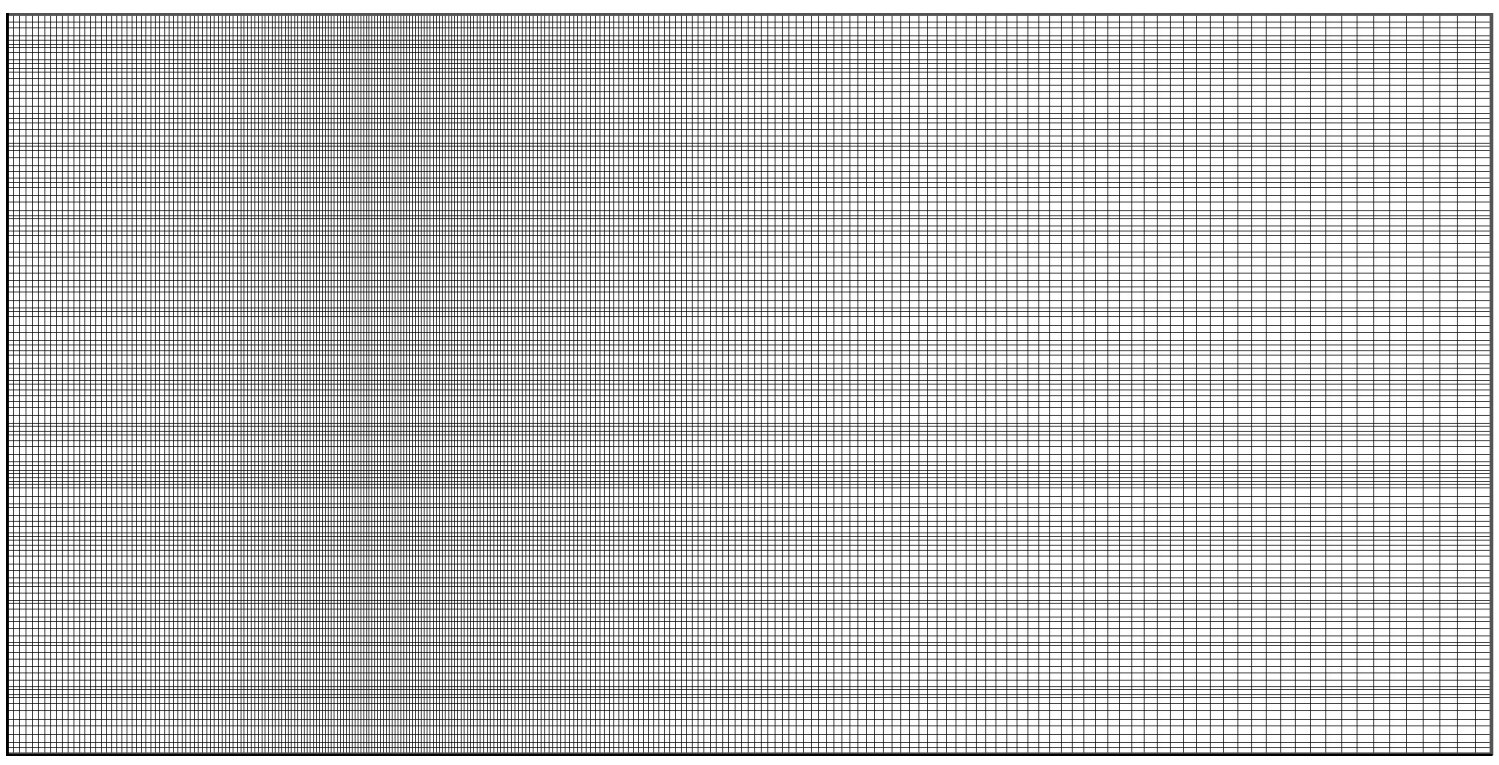

Figure 14: Non-uniform grids for shock-vortex interaction problem. In $x$, stretched grids; in $y$, ALT-RAND grids.

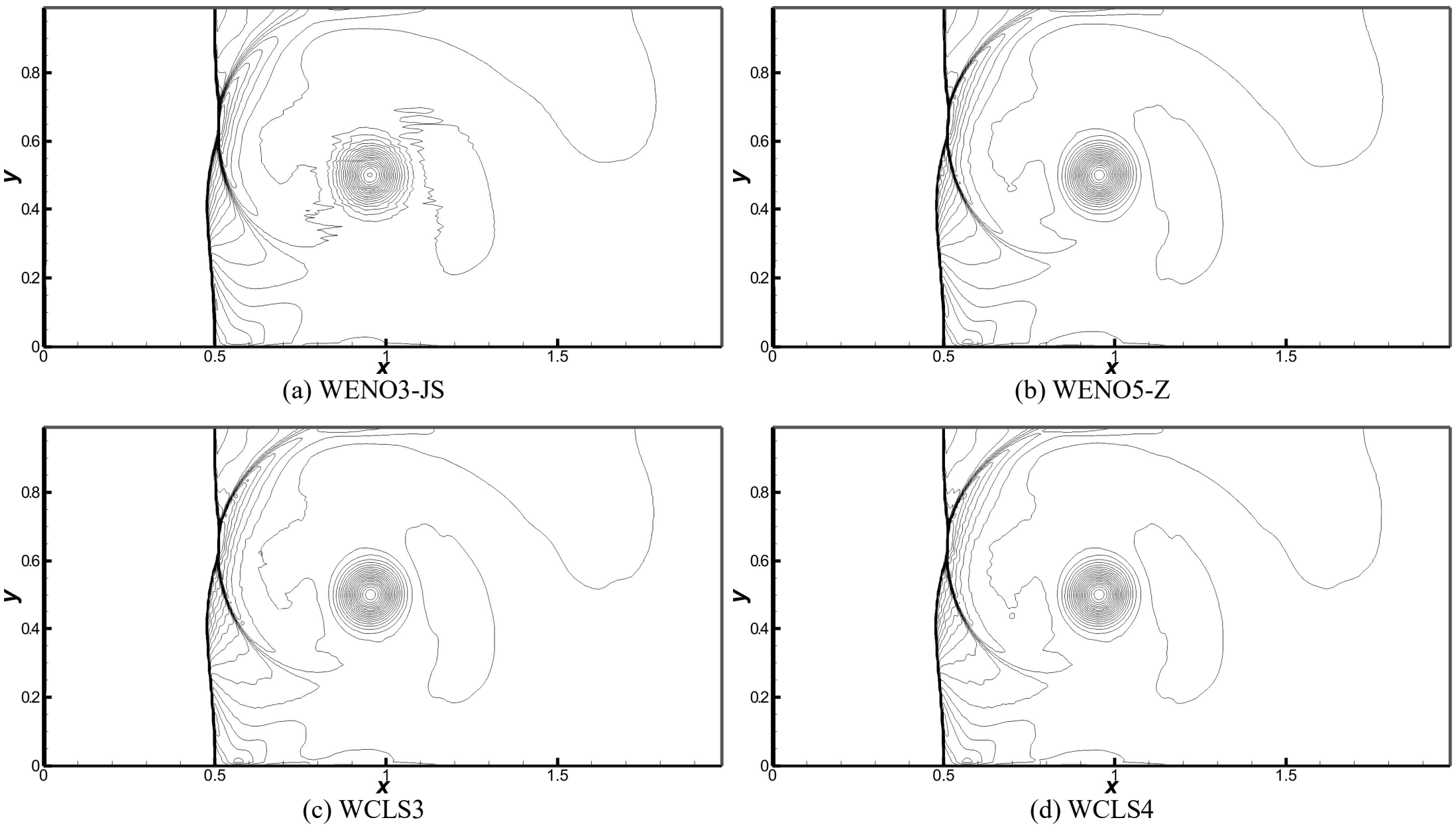

(a) WENO3-JS (b) WENO5-Z

(c) WCLS3 (d) WCLS4

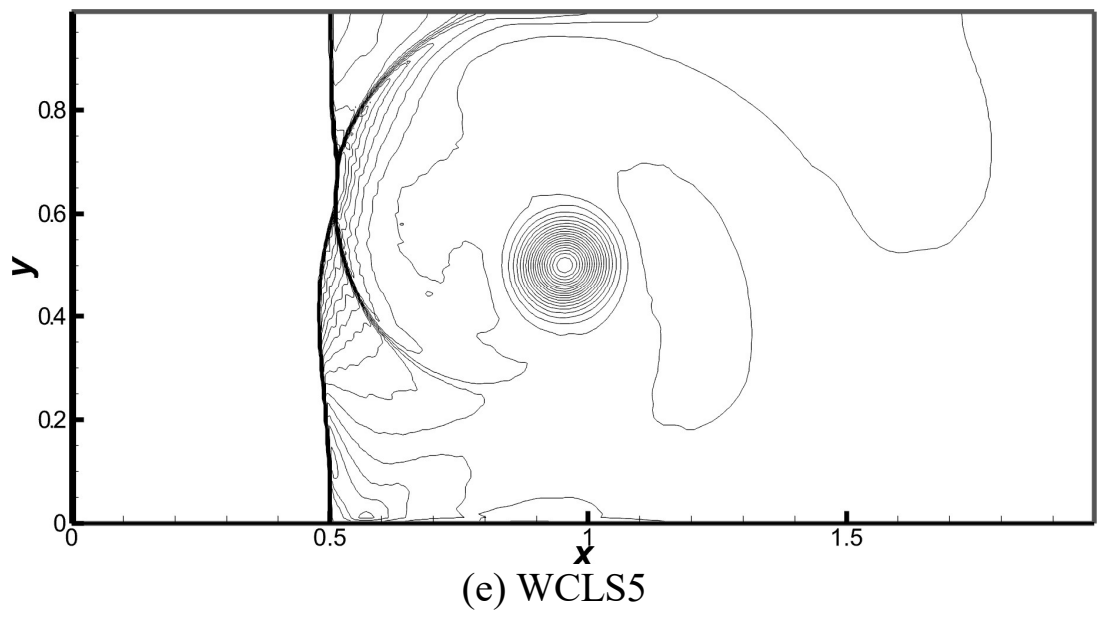

(e) WCLS5

Figure 15: Density contours for the shock vortex interaction with 30 lines ranging from 1.01 to 1.22.

## 4.8 2D Riemann Problem

A 2D Riemann problem with initial condition

$$[\rho, u_1, u_2, P] = \begin{cases} 1.5, & 0, & 0, & 1.5, & 0.8 \le x \le 1, 0.8 \le y \le 1, \\ 0.5323, & 1.206, & 0, & 0.3, & 0 \le x < 0.8, 0.8 \le y \le 1, \\ 0.138, & 1.206, & 1.206, & 0.029, & 0 \le x < 0.8, 0 \le y < 0.8, \\ 0.5323, & 0, & 1.206, & 0.3, & 0.8 \le x \le 1, 0 \le y < 0.8, \end{cases} \tag{66}$$

is utilized to check the performance of the proposed WCLS schemes. The simulation domain is $[0,1] \times [0,1]$ with $400 \times 400$ uniform control volumes. The simulation parameters are $\epsilon_{\text{Euler}} = \epsilon_{\text{Recon}} = 1 \times 10^{-4}$, $N_{\text{inner}} = 6$ and CFL = 1.0. To better preserve the symmetry of the flow structure, Eq. (55) of Euler equations is also iterated by GMRES method with $N_{\text{GMRES}} = 1$ in this case.

The density contours are presented in Fig. 16. As the shocks interact and move, Helmholtz instabilities along the slip lines are triggered. A jet with a mushroom-shaped head forms along the $y = x$ line. The WCLS schemes capture the shocks and contact discontinuities robustly and in high resolution. The resolution of contact discontinuities by the WCLS schemes are much higher than the WBAP-CLS schemes and the WENO schemes.

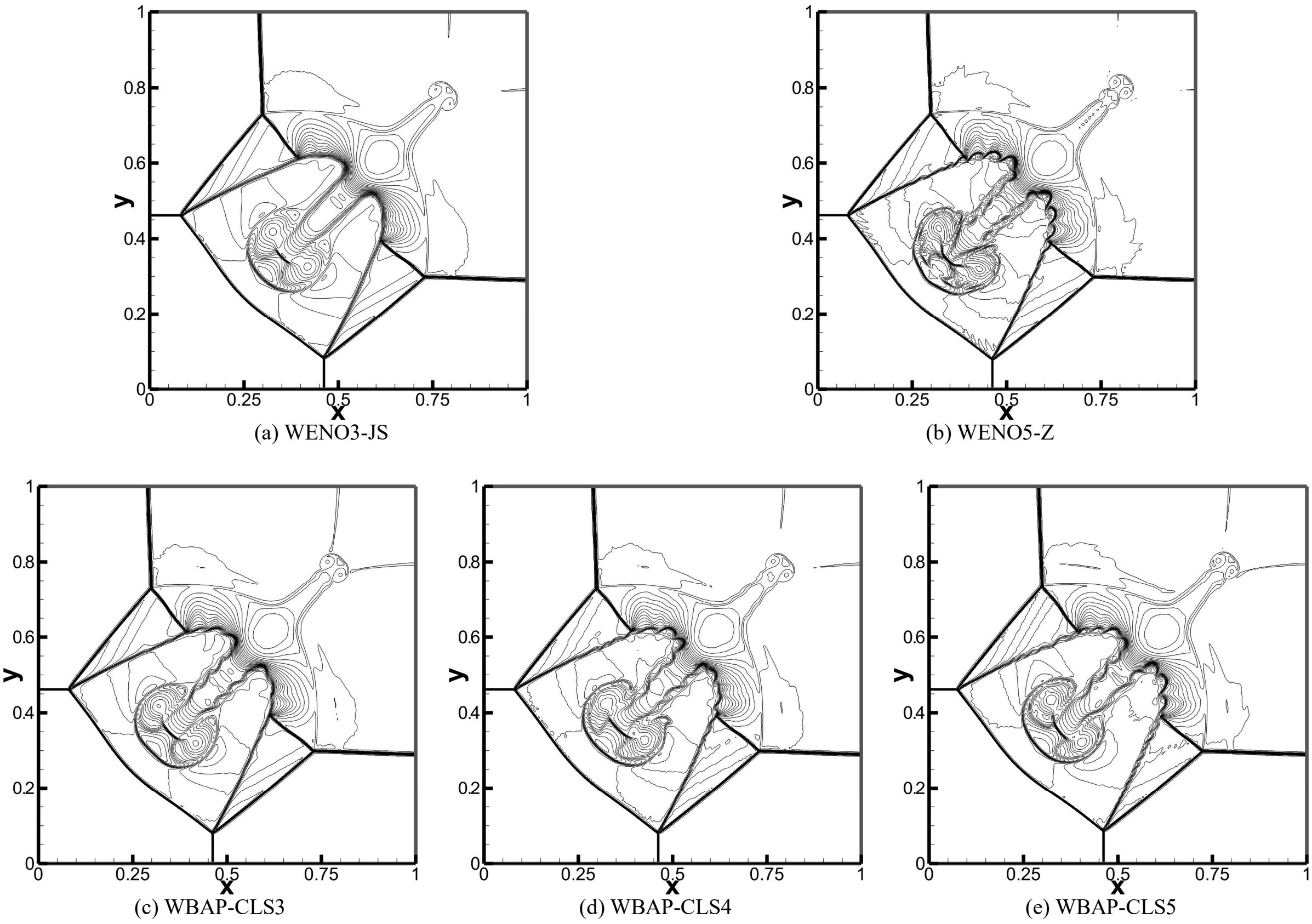

(a) WENO3-JS (b) WENO5-Z

(c) WBAP-CLS3 (d) WBAP-CLS4 (e) WBAP-CLS5

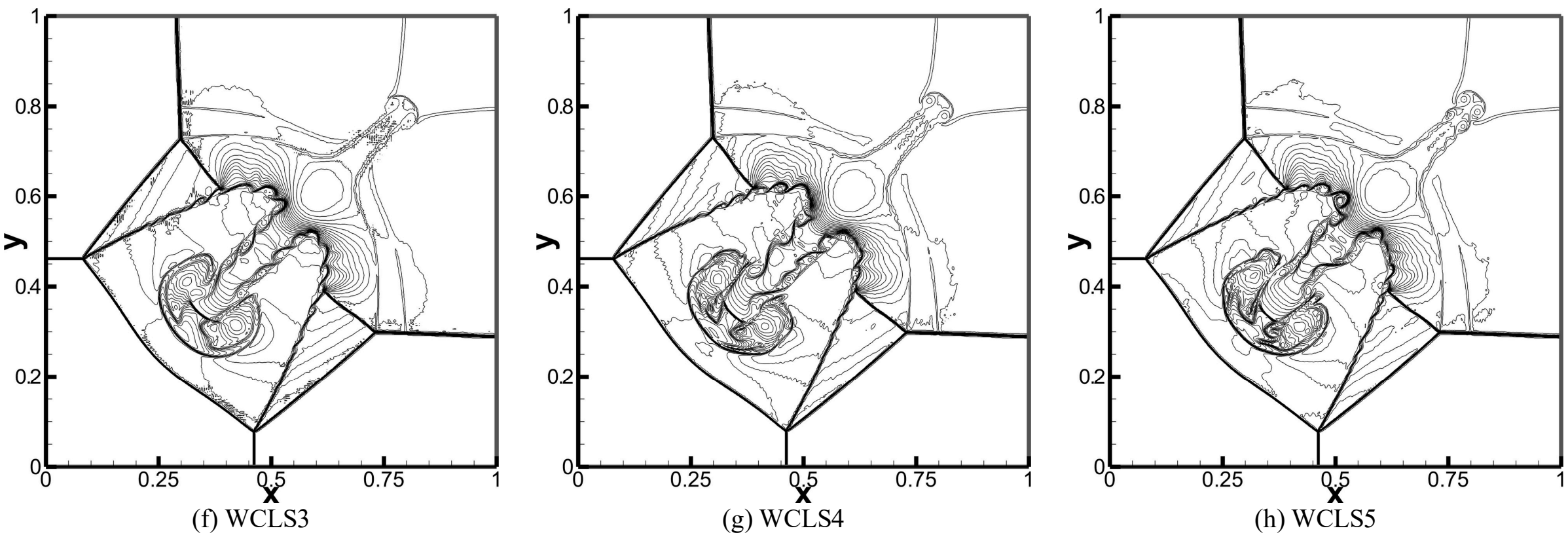

Figure 16: Density contours for the 2D Riemann problem with 32 lines ranging from 0.2 to 1.7. $N_x \times N_y = 400 \times 400$, $t$ = 0.8 and CFL = 1.0.

## 4.9 Double Mach Reflection

Double Mach reflection is another classic benchmark for validating shock-capturing schemes for compressible flows [39]. In this problem, a Mach 10 shock moves downward to the right and impinges on an inviscid wall, producing a bow shock, a Mach stem, and a slip line.

The simulation domain is $[0,4] \times [0,1]$ and the simulation end time is $t$ = 0.2. The inviscid wall locates at $x > 1/6$. The initial condition is

$$[\rho, u_1, u_2, P] = \begin{cases} 8.0, 7.1447, -4.125, 116.5, & x \leq 1/6 + \sqrt{3}y/3, \\ 1.4, 0, 0, 1, & x > 1/6 + \sqrt{3}y/3, \end{cases} \tag{67}$$

Other simulation parameters are $\epsilon_{\text{Euler}} = \epsilon_{\text{Recon}} = 1 \times 10^{-4}$, $N_{\text{inner}} = 6$ and $\text{CFL} = 1.0$.

The results on meshes with grid size $h = 1/240$ and $h = 1/480$ are shown in Figs. 17 and 19, respectively. The close views around the Mach stems and the slip lines are shown in Figs. 18 and 20. The results confirm the non-oscillatory property of the proposed WCLS schemes. Additionally, as demonstrated in the close views, the WCLS schemes resolve the contact discontinuities in much higher resolution than the traditional WENO schemes and the WBAP-CLS schemes. Even the third-order WCLS3 scheme can capture the Helmholtz instabilities in a comparable resolution with WENO5-Z.

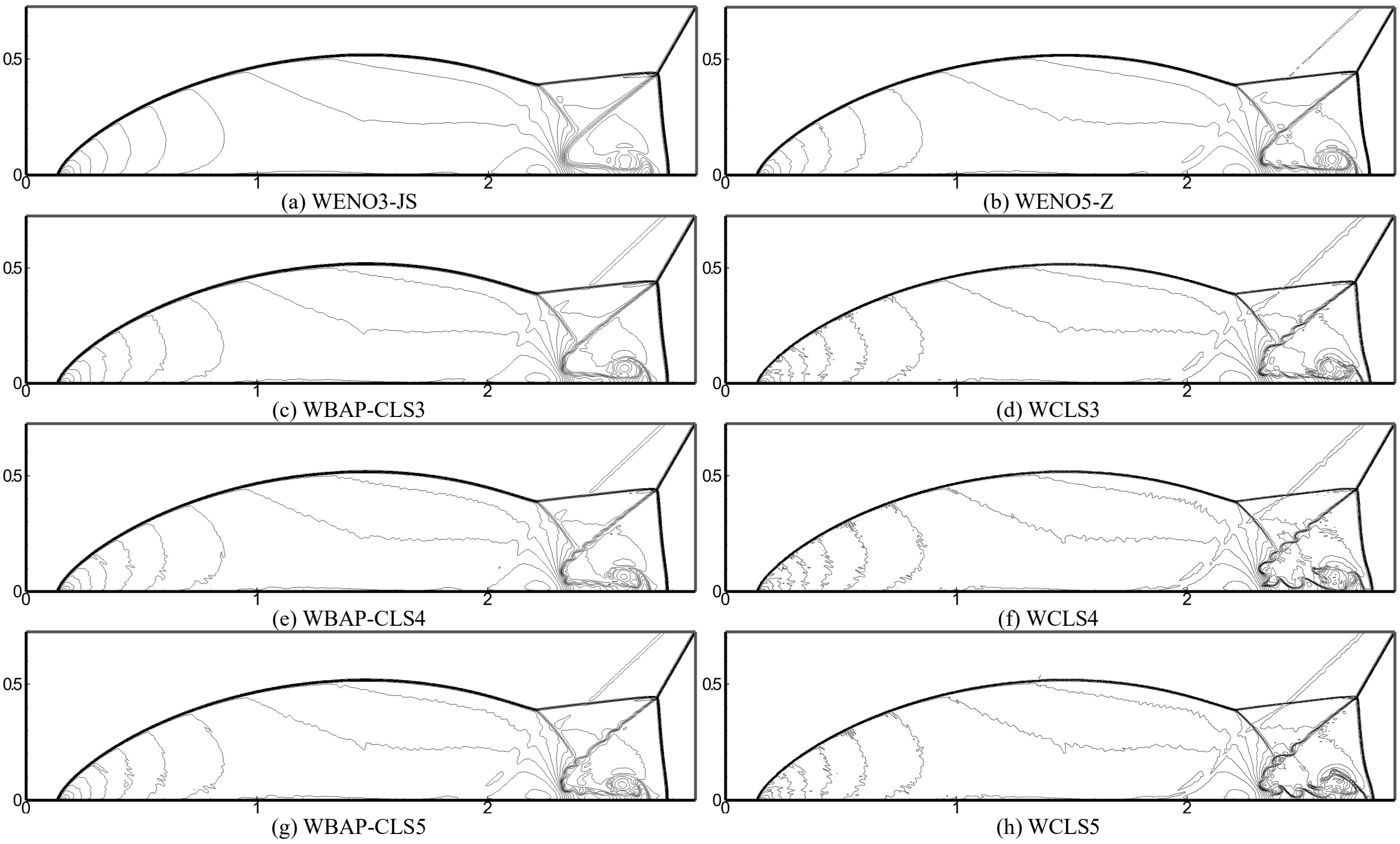

Figure 17: Density contours for the double Mach reflection with 32 lines ranging from 2.2 to 22. $N_x \times N_y = 960 \times 240$, $t$ = 0.2 and CFL =1.0.

(a) WENO3-JS (b) WENO5-Z

(c) WBAP-CLS3 (d) WBAP-CLS4 (e) WBAP-CLS5

(f) WCLS3 (g) WCLS4 (h) WCLS5

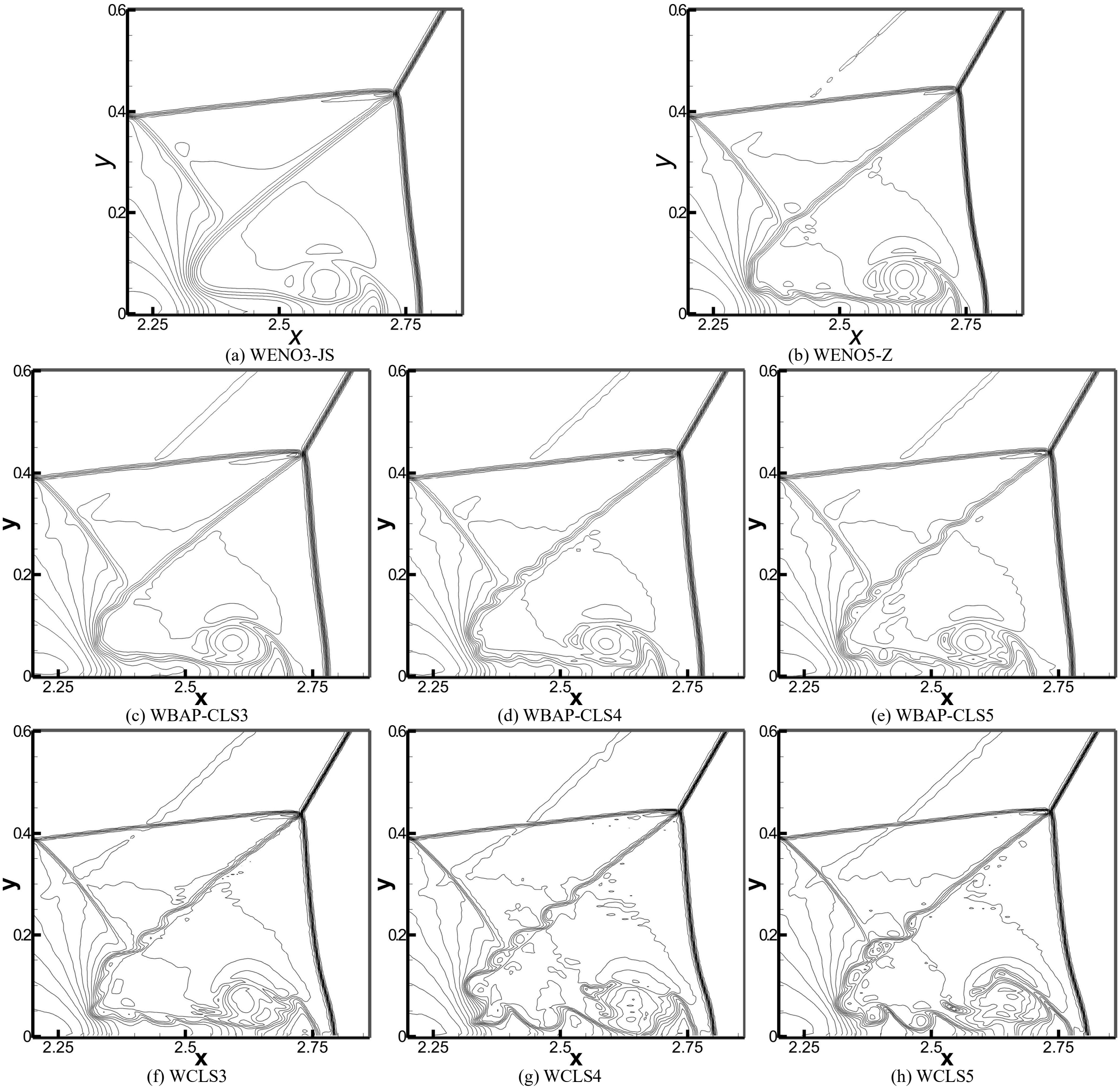

Figure 18: Close view of density contours for the double Mach reflection. $N_x \times N_y = 960 \times 240$.

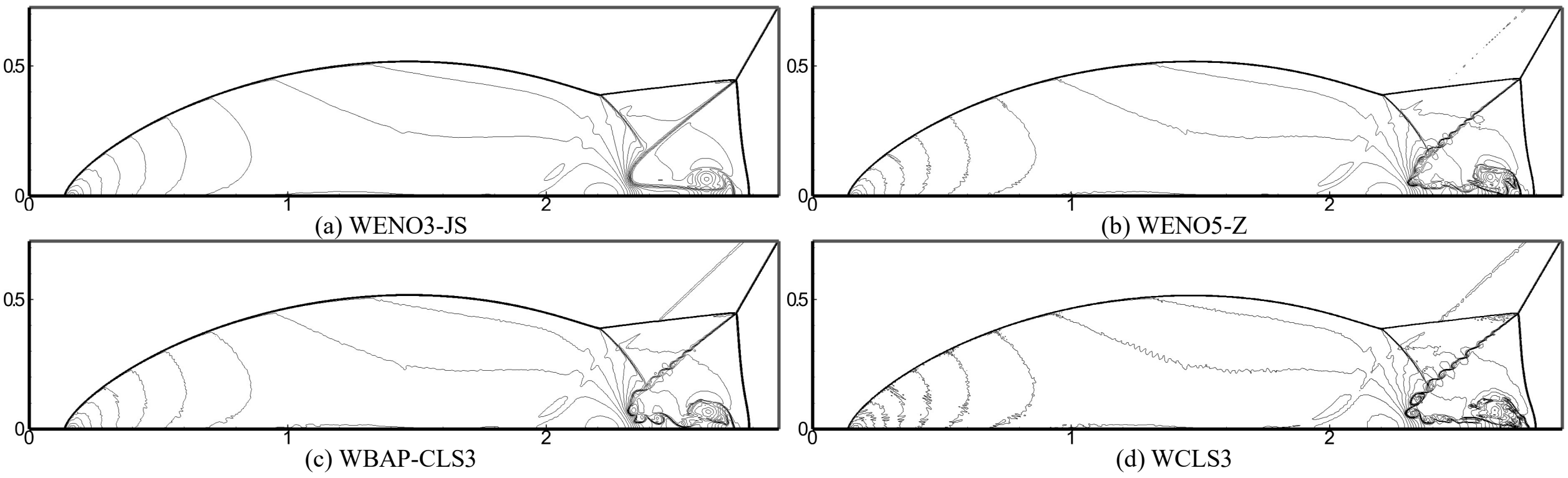

(a) WENO3-JS (b) WENO5-Z

(c) WBAP-CLS3 (d) WCLS3

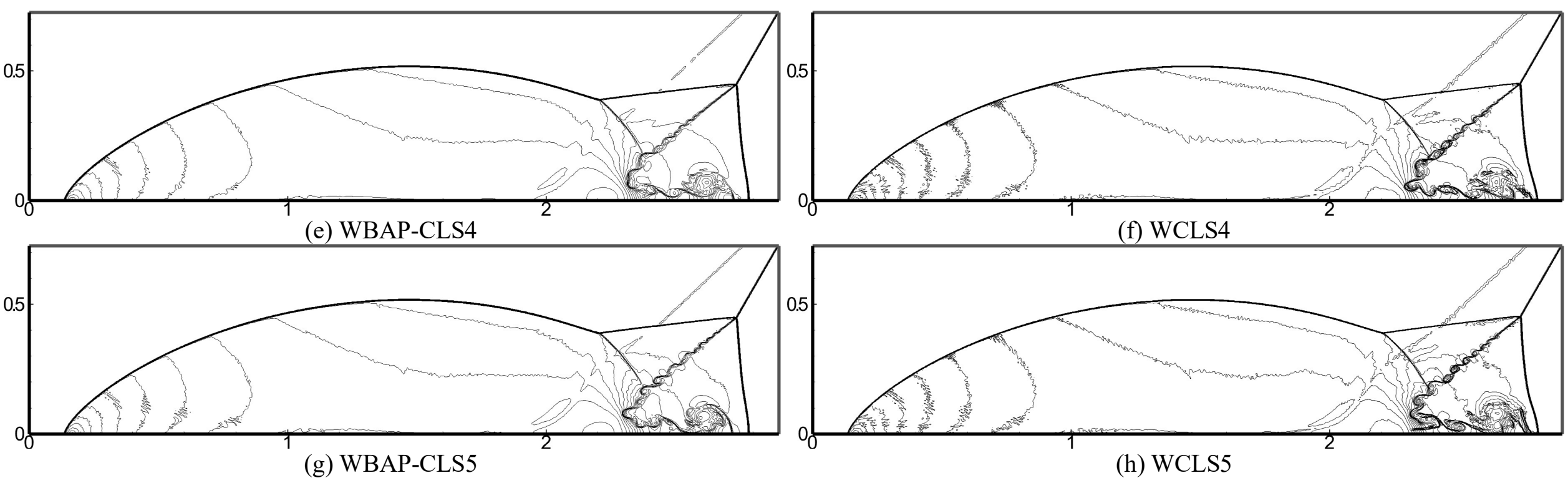

Figure 19: Density contours for the double Mach reflection with 32 lines ranging from 2.2 to 22. $N_x \times N_y = 1920 \times 480$, $t = 0.2$ and CFL =1.0.

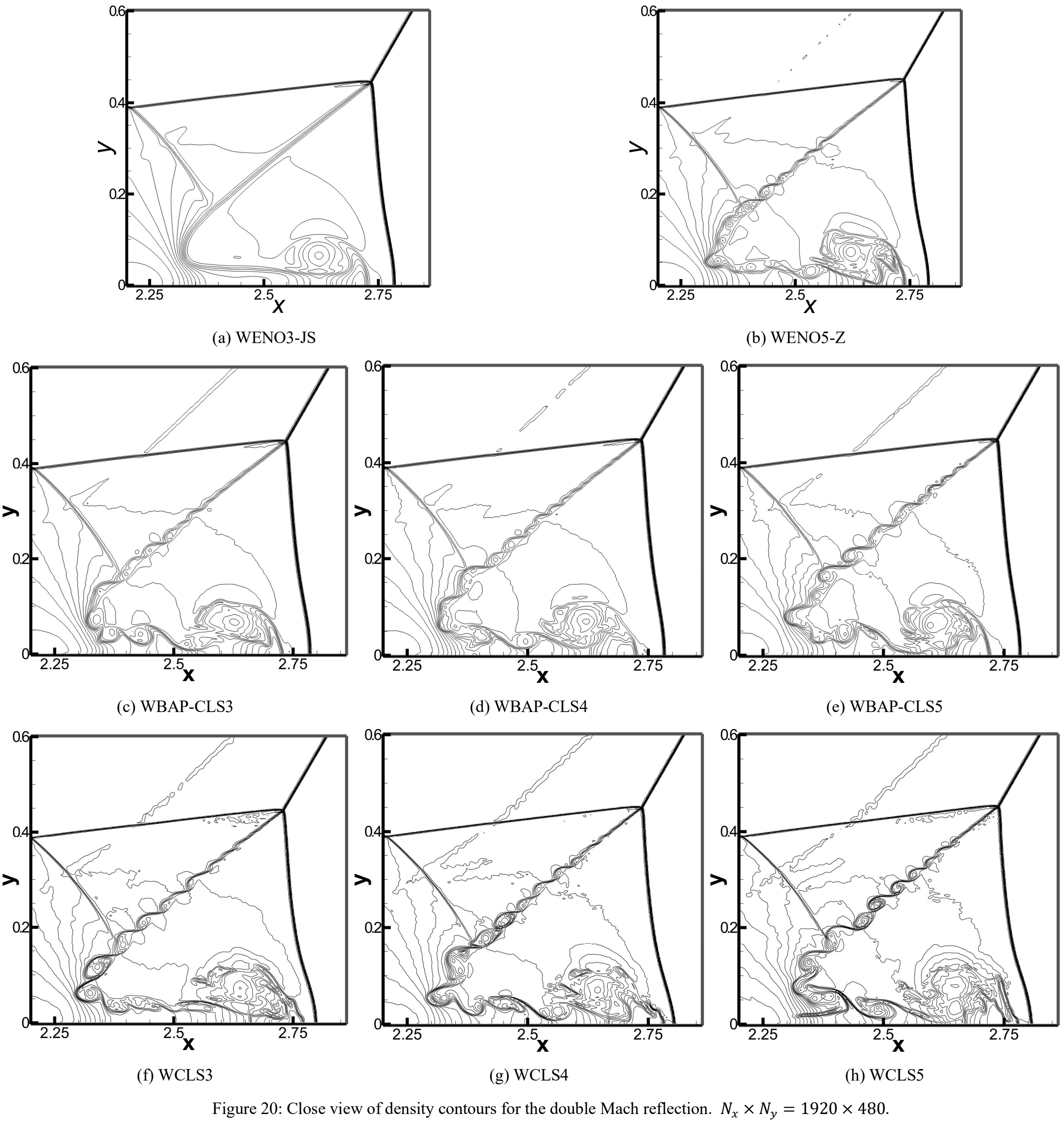

(a) WENO3-JS

(b) WENO5-Z

(c) WBAP-CLS3

(d) WBAP-CLS4

(e) WBAP-CLS5

(f) WCLS3

(g) WCLS4

(h) WCLS5

Figure 20: Close view of density contours for the double Mach reflection. $N_x \times N_y = 1920 \times 480$.

### 4.10 2D High-Mach-Number Astrophysical Jet [40]

To further assess the resolution and robustness of the proposed reconstruction methods for highly compressible flows with strong discontinuities, we consider the two-dimensional extremely high-Mach-number astrophysical jet problem. This test involves strong shocks, contact discontinuities, thin shear layers, and complex wave interactions.

For the Mach 80 case, the computational domain is $[0,\ 2]\times[-0.5,0.5]$.The ambient gas is initialized by $(\rho,u_1,u_2\,,p)=(0.5,0,0,0.4127)$ . At the left boundary, the jet inflow is prescribed as $(\rho,u_1,u_2,p)=(5,30,0,0.4127), y\in[-0.05,0.05]$, while the remaining part of the left boundary is set to the ambient state. The right, top, and bottom boundaries are treated as outflow boundaries. The ratio of specific heats is taken as $\gamma$=5/3. The mesh resolution is 448×224. The simulation parameters are $\epsilon_{\text{Euler}}=\epsilon_{\text{Recon}}=1\times10^{-4}$, $N_{\text{inner}}=6$ and CFL = 1.0.

For the Mach 2000 case, the computational domain is $[0,1]\times[-0.25,\ 0.25]$. The ambient state remains $(\rho,u_1,u_2,p)=(0.5,0,0,0.4127)$ , and the inflow jet is given by $(\rho,u_1,u_2,p)=(5,800,0,0.4127), y\in[-0.05,0.05]$, with the rest of the left boundary set to the ambient condition. The right, top, and bottom boundaries are again imposed as outflow boundaries. The mesh resolution is 800×400. The simulation parameters are $\epsilon_{\text{Euler}}=\epsilon_{\text{Recon}}=1\times10^{-4}$, $N_{\text{inner}}=6$ and CFL = 1.0.

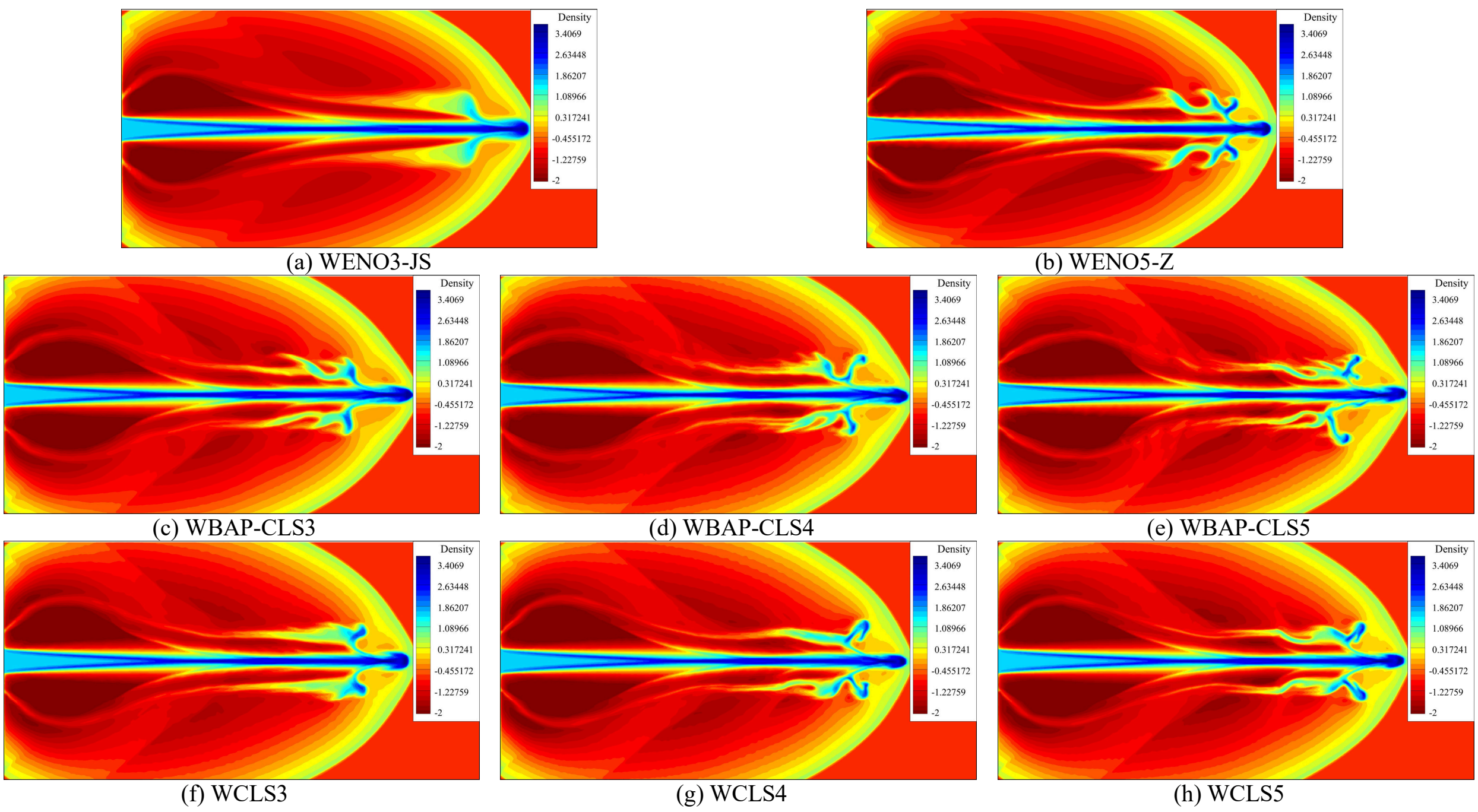

(a) WENO3-JS (b) WENO5-Z
(c) WBAP-CLS3 (d) WBAP-CLS4 (e) WBAP-CLS5
(f) WCLS3 (g) WCLS4 (h) WCLS5

Figure 21: Density contours of the 2D high-Mach-number astrophysical jet computed by different schemes for Mach 80.

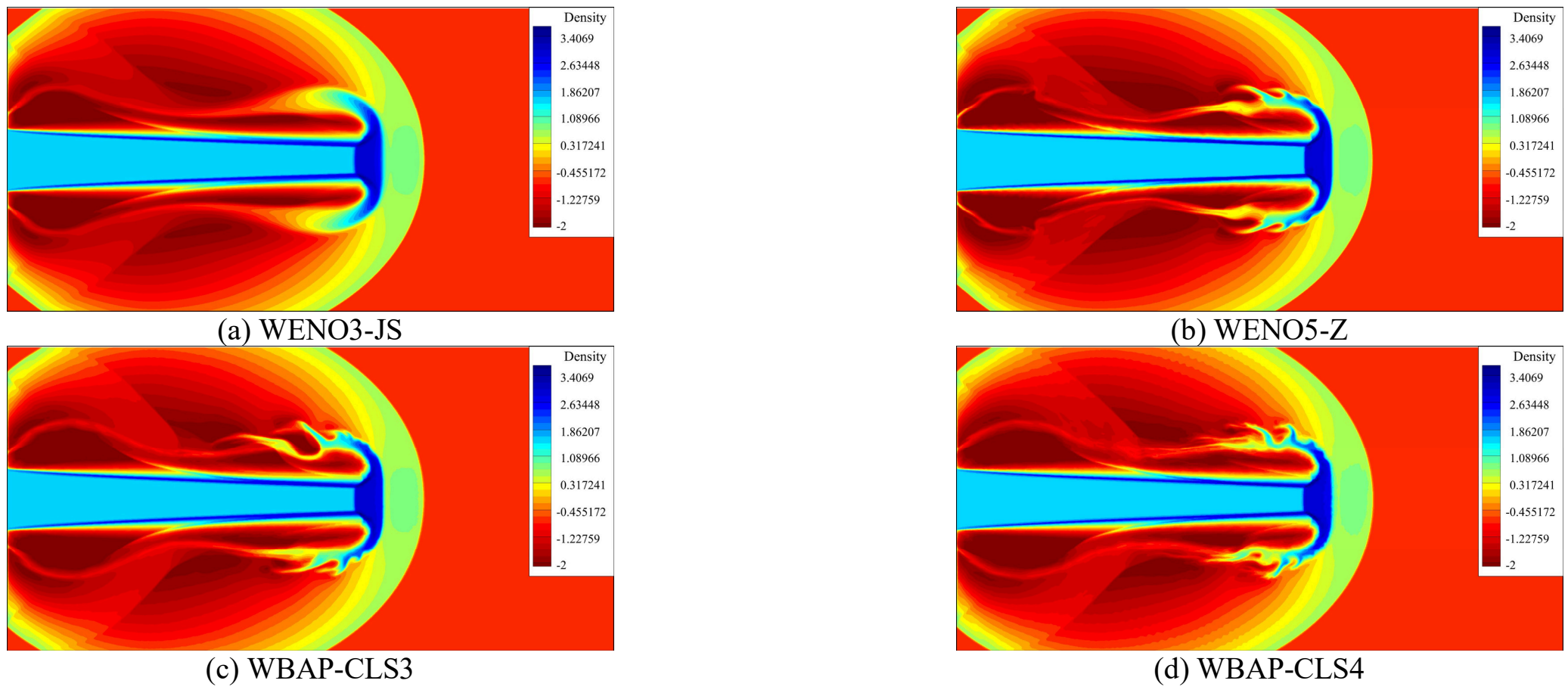

(a) WENO3-JS (b) WENO5-Z
(c) WBAP-CLS3 (d) WBAP-CLS4

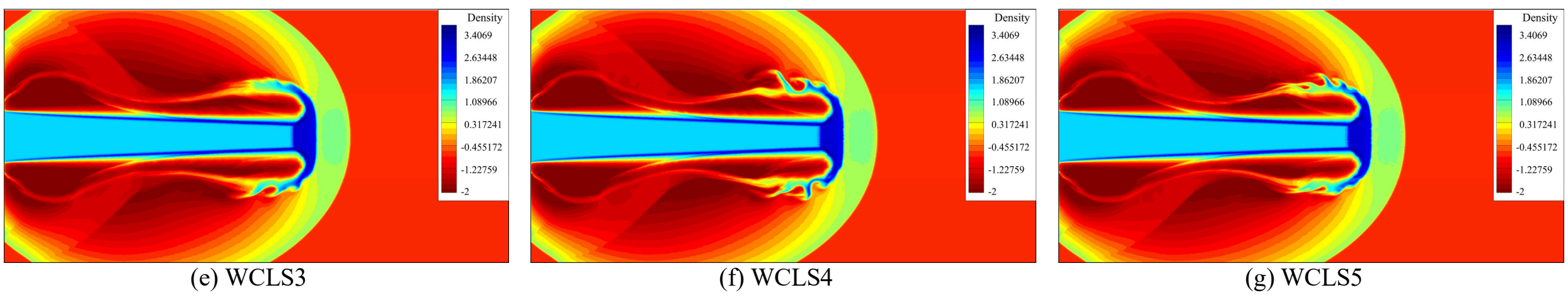

(e) WCLS3 (f) WCLS4 (g) WCLS5

Figure 22: Density contours of the 2D high-Mach-number astrophysical jet computed by different schemes for Mach 2000.

As shown in Figs. 21 and 22, the astrophysical jet test demonstrates that the WCLS schemes provide sharper resolution than the reference WENO methods, with finer structures near the jet head and along the shear layers better preserved. For the much more demanding Mach 2000 case, WBAP-CLS5 scheme fails because the WBAP limiter cannot suppress numerical oscillations under this extreme condition, whereas the proposed WCLS5 still captures the flow structures without oscillations.

### 4.11 Viscous Shock Tube

A viscous shock tube problem [41] with Reynolds number 200 is considered to further validate the WCLS schemes for viscous problems. The initial condition is

$$[\rho, u_1, u_2, P] = \begin{cases} 120, 0, 0, 120/\gamma, & x < 0.5, \\ 1.2, 0, 0, 1.2/\gamma, & x \geq 0.5. \end{cases} \tag{68}$$

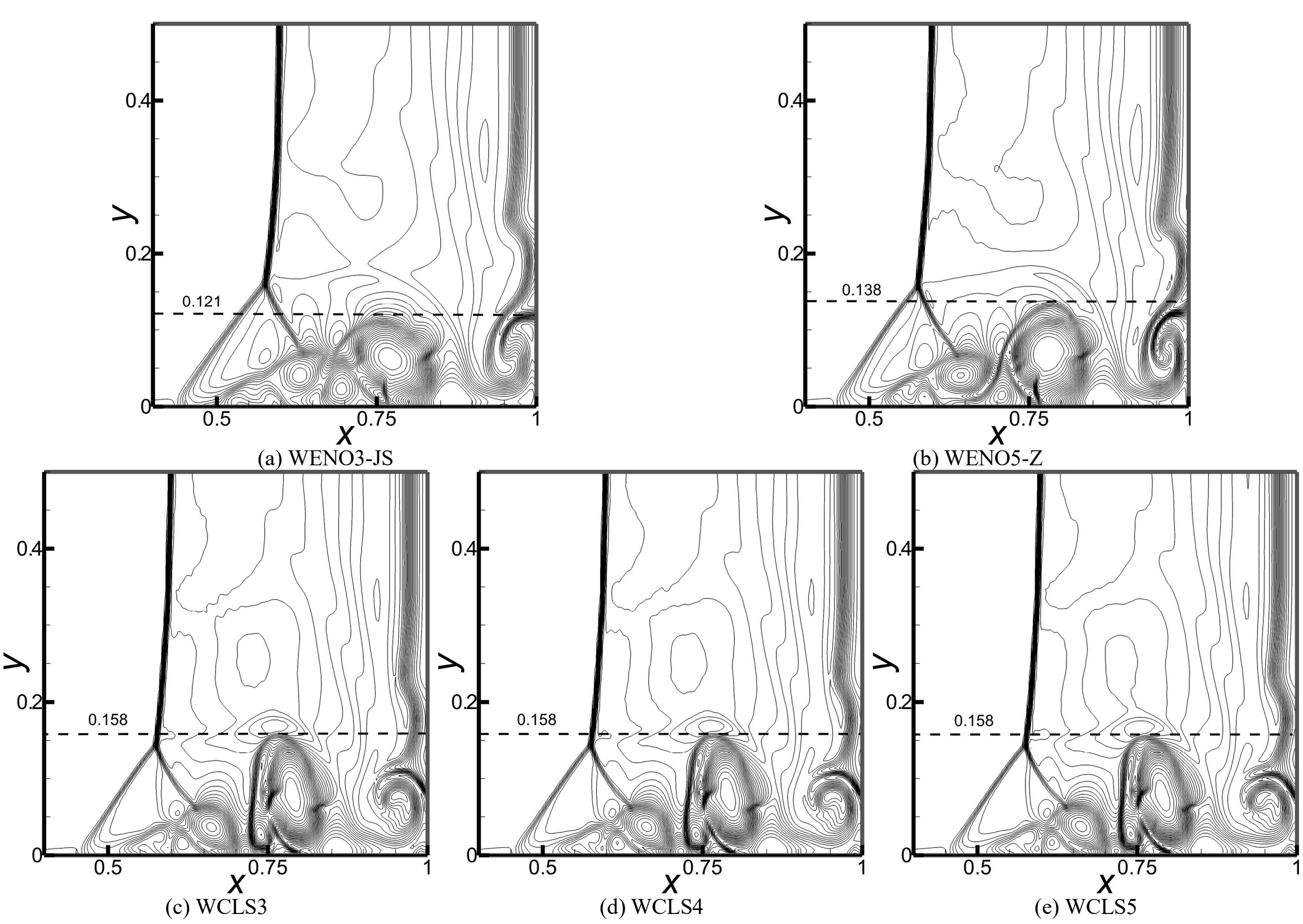

(a) WENO3-JS (b) WENO5-Z

(c) WCLS3 (d) WCLS4 (e) WCLS5

Figure 23: Density contours for the viscous shock tube problem. 30 equally spaced contour lines from 22.26 to 121.34. $N_x \times N_y = 300 \times 300$.

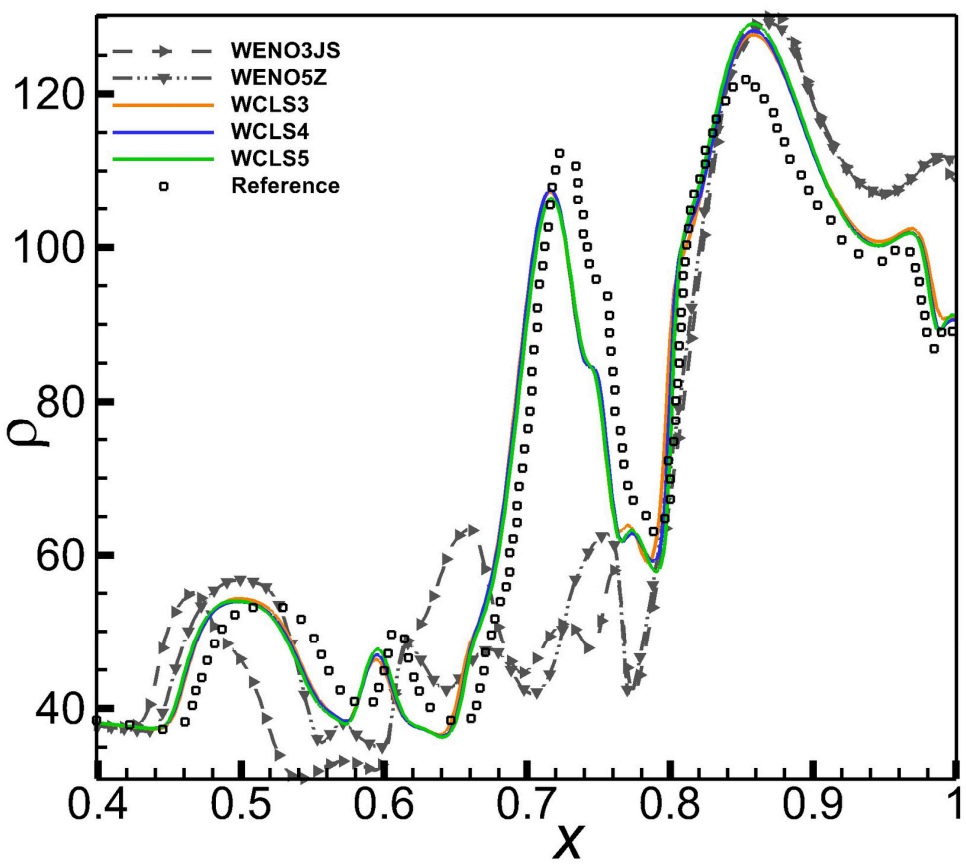

Figure 24: Distribution of the density along the bottom wall using different schemes.

The simulation domain is $[0,1] \times [0,1]$ with four adiabatic non-slip walls, consisting of $300 \times 300$ uniform control volumes. The simulation parameters are $\epsilon_{\mathrm{NS}} = \epsilon_{\mathrm{Recon}} = 1 \times 10^{-6}$, $N_{\mathrm{inner}} = 6$ and $\mathrm{CFL} = 2.0$. The simulation end time is $t = 1.0$.

Upon start, the shock moves rightward and reflects from the wall. The reflected shock interacts with the boundary layer and $\lambda$ shaped vortex is formed. The density contours for different schemes are shown in Fig. 23. As demonstrated, the vortex height by the proposed WCLS schemes is closer to the reference value of 0.165, illustrating a more accurate numerical viscous flux due to the derivative correction of the WCLS schemes along reconstruction axes. The density distribution along the bottom of the wall is shown in Fig. 24, where the reference result is from the work of Huang et al. [41], confirming again the superior resolution of the proposed WCLS schemes.

### 4.12 Shock–Mixing-Layer Interaction

In this section, a shock mixing layer interaction with Reynolds number 500 is simulated to validate the performance of the proposed WCLS schemes on slightly non-uniform grids. The Reynolds number is defined as $\mathrm{Re} = \rho_{\mathrm{mean}} U_{\mathrm{mean}} \delta/\mu$, where $\rho_{\mathrm{mean}} = (\rho^1 + \rho^2)/2$ and $U_{\mathrm{mean}} = (u_1^1 + u_1^2)/2$ denotes the averaged density and streamwise velocity, respectively, $\delta = 1$ is the initial thickness of the mixing layer. Superscripts 1 and 2 denote that the corresponding variable is defined for the upper and lower inflow of the mixing layer, respectively. In this problem, $[\rho^1, u_1^1, u_2^1, P^1] = [1.6374,3.0,0.0,0.3327]$ and $[\rho^2, u_1^2, u_2^2, P^2] = [0.3626,2.0,0.0,0.3327]$. The initial flow variables are set by the profile $\phi = (\phi^1 + \phi^2)/2 + (\phi^1 - \phi^2)\tanh(2y/\delta)/2$ where $\phi$ denotes $\rho$, $u_1$, $u_2$ or $P$.

An oblique shock originating from the upper-left corner interacts with the mixing layer and Helmholtz instabilities form. The flow condition after the shocks is $[\rho^3, u_1^3, u_2^3, P^3] = [2.1101, 2.9709, -0.1367, 0.4754]$. Supersonic inflow and outflow conditions are applied for the left and right boundaries; inviscid wall condition is applied for the bottom boundary and supersonic inflow with after-shock condition is applied on the upper boundary. Fluctuations of $y$-component of velocity are applied on the inflow boundary, the detailed expression of which can be found in Eq. (3.18) from the work of Yee et al. [42].

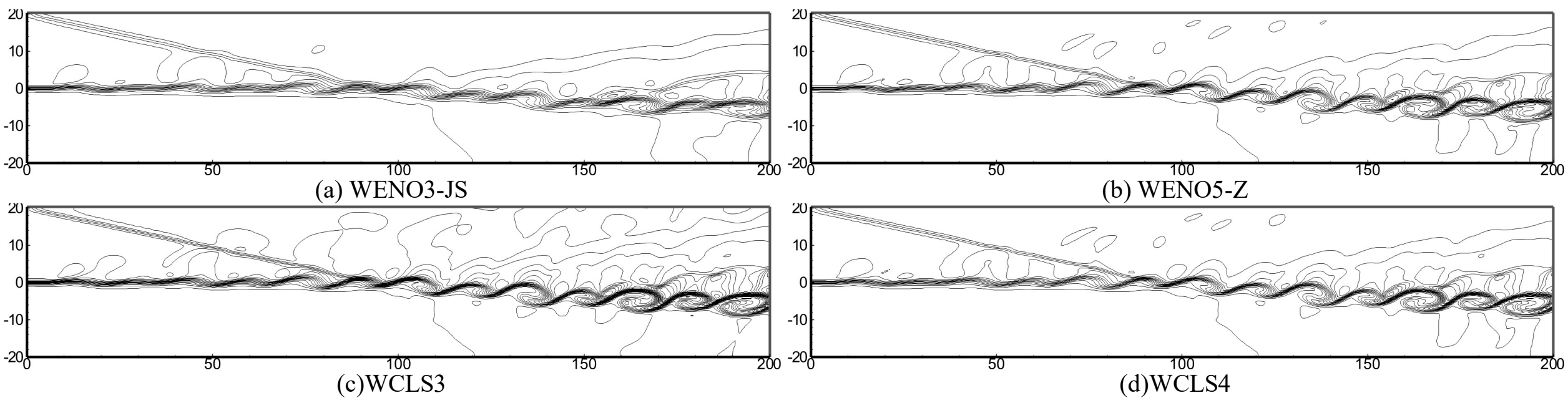

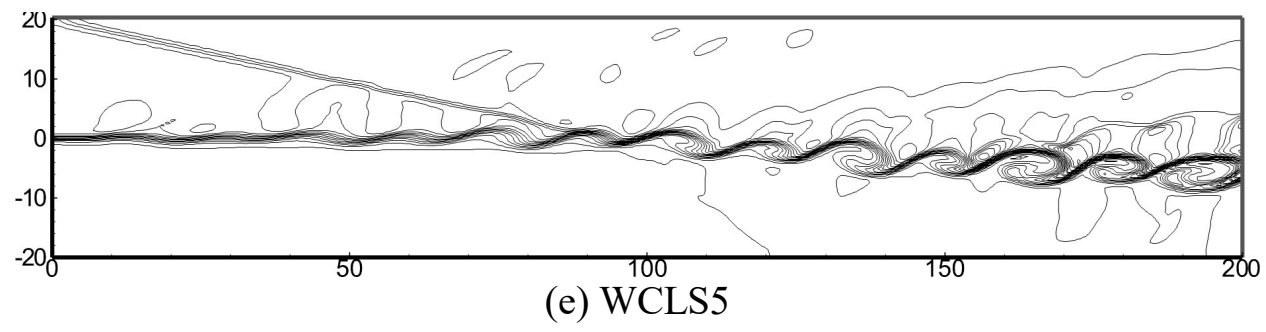

(e) WCLS5

Figure 25: Density contours by the WENO and WCLS schemes for shock–mixing-layer interaction. 20 equally spaced 0.4 to 2.8.

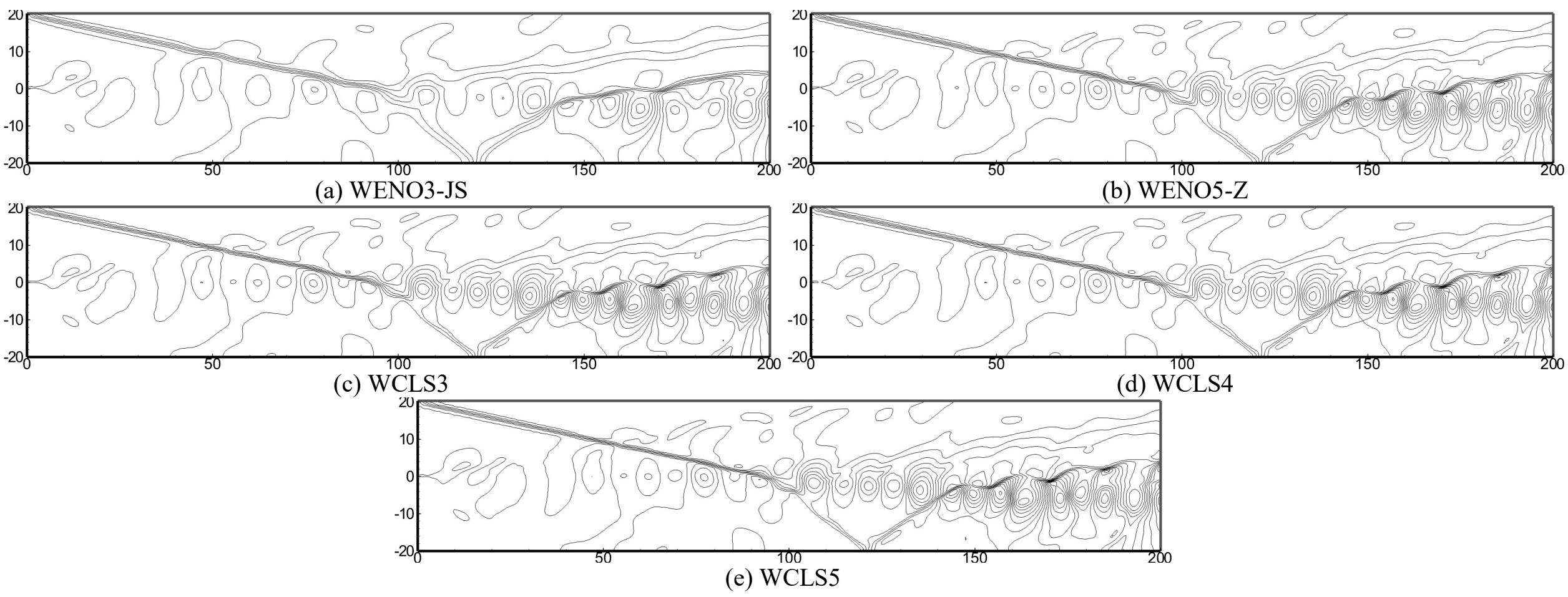

(a) WENO3-JS (b) WENO5-Z (c) WCLS3 (d) WCLS4 (e) WCLS5

Figure 26: Pressure contours by the WENO and WCLS schemes for shock–mixing-layer interaction. 20 equally spaced 0.23 to 0.73.

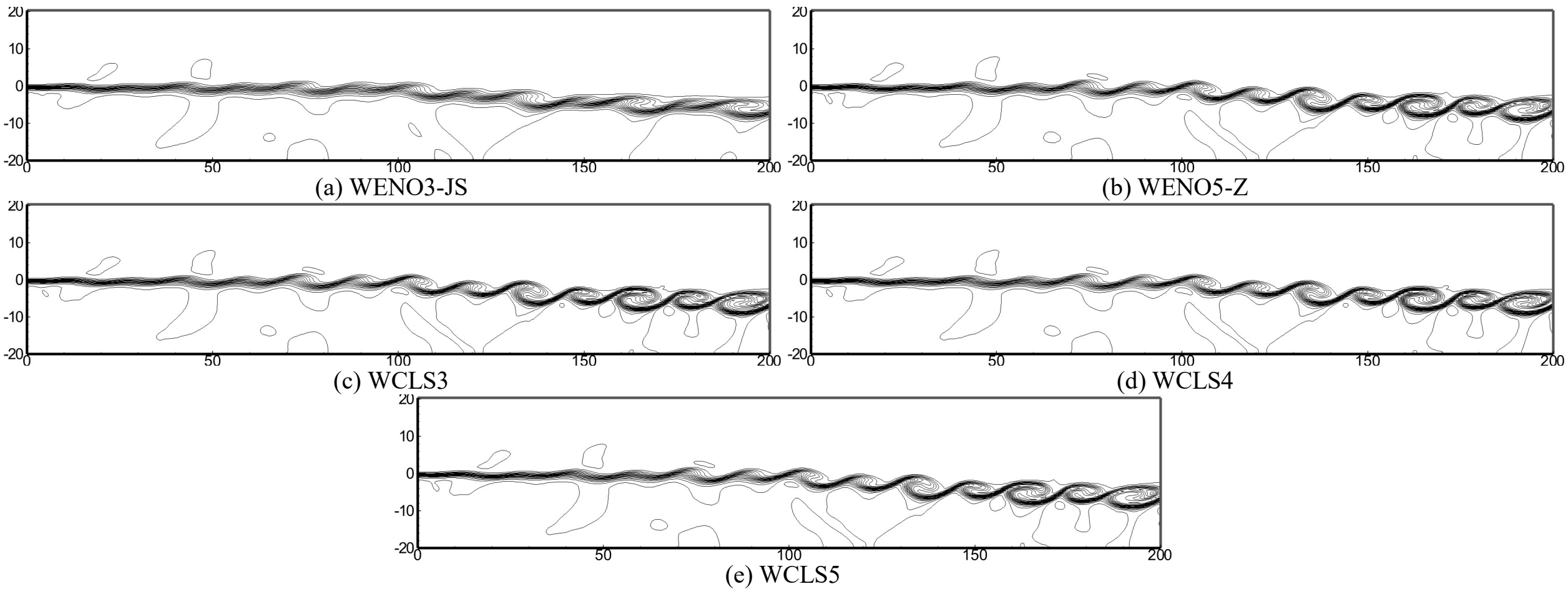

(a) WENO3-JS (b) WENO5-Z (c) WCLS3 (d) WCLS4 (e) WCLS5

Figure 27: Temperature contours by the WENO and WCLS schemes for shock–mixing-layer interaction. 20 equally spaced 0.2 to 1.1.

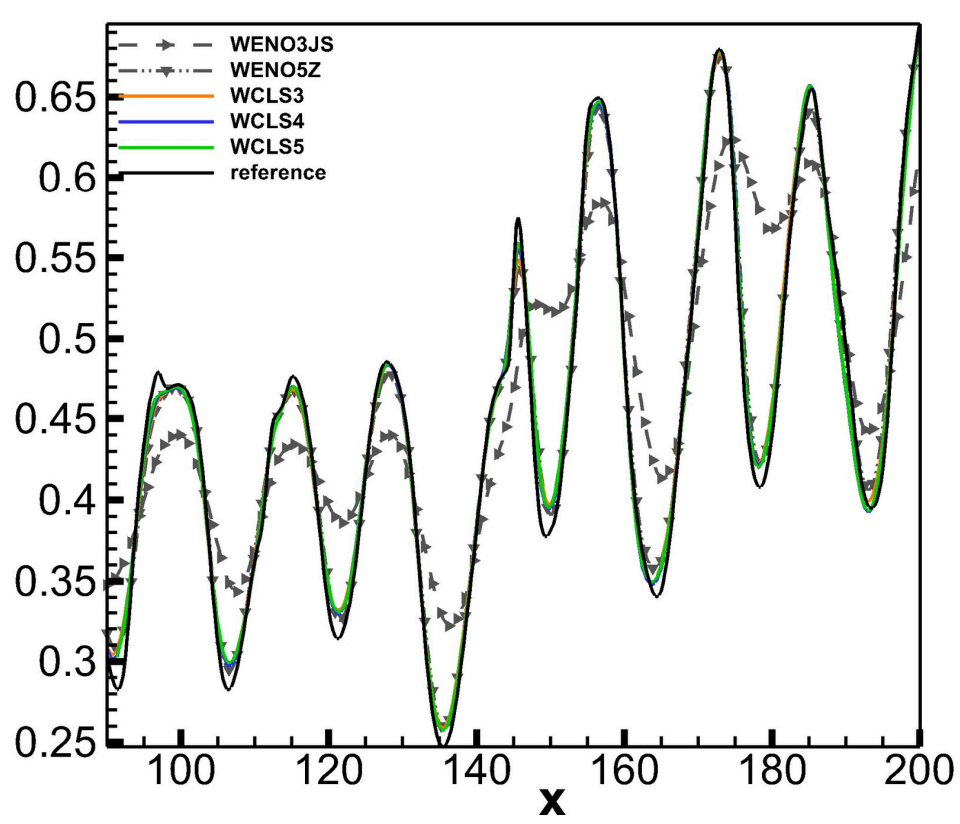

Figure 28: Pressure distributions along the line from point (90, 0) to point (200, −6).

The simulation domain is $[0, 200] \times [-20, 20]$ discretized by $320 \times 80$ control volumes. The grids are uniform in $x$-direction and nonuniform in $y$-direction. The $y$ coordinates are mapped in a form $20\tanh(\eta)/\tanh 1$ from uniform computational grids $\eta \in [-1,1]$. The simulation parameters are $\epsilon_{\text{NS}} = \epsilon_{\text{Recon}} = 1 \times 10^{-6}$, $N_{\text{inner}} = 6$ and $\text{CFL} = 2.0$. The simulation end time is $t = 160.0$.

Density contours are presented in Fig. 25. The WENO3-JS scheme is too dissipative on this mesh to capture the fine vortex structures. WENO5-Z and the proposed WCLS schemes all resolve the Helmholtz instabilities in a good resolution. As the last vortex around $x = 200$ shows, the WCLS schemes can resolve finer flow structures than the WENO5-Z scheme. Figures 26 and 27 present the contours of pressure and temperature, respectively, confirming the non-oscillatory property of the proposed WCLS schemes. Eventually, the pressure profiles along the line from point (90,0) to point (200, −6) are shown in Fig. 28, where the reference result is from WENO5-Z scheme on mesh with $640 \times 160$ control volumes. The result in Fig. 28 further confirms the high resolution of the proposed WCLS schemes for viscous problems and on non-uniform grids.

### 4.13 Steady Subsonic Viscous Flow around NACA0012

For viscous problems with boundary layers, the implicit scheme can speed up the simulation with a relatively higher CFL number. This case further validates the proposed WCLS schemes on non-uniform curvilinear grids. The chord length of the NACA0012 airfoil is 1. Subsonic fluid flows past the airfoil at a Mach number of 0.5 and an attack angle of $0°$. The Reynolds number is 5000. The distance of the first layer of grids from the viscous adiabatic wall is $1.24 \times 10^{-4}$. The simulation domain is a circle with radius of 20. There are $145$ grid points along the surface of airfoil and 107 points along the radius of computational domain. The computational grids are shown in Fig. 29. The CFL number is 1000 and simulation is continued until converged.

The skin-friction and pressure coefficients are presented in Fig. 30. The reference result is taken from Li [43]. The proposed WCLS schemes show clear order convergence and yield more accurate skin-friction and pressure coefficients than the WENO3-JS and WENO5-Z schemes.

Finally, Fig. 31 presents the density contours around the NACA0012 airfoil. On this mesh, the proposed WCLS schemes produce highly symmetric results.

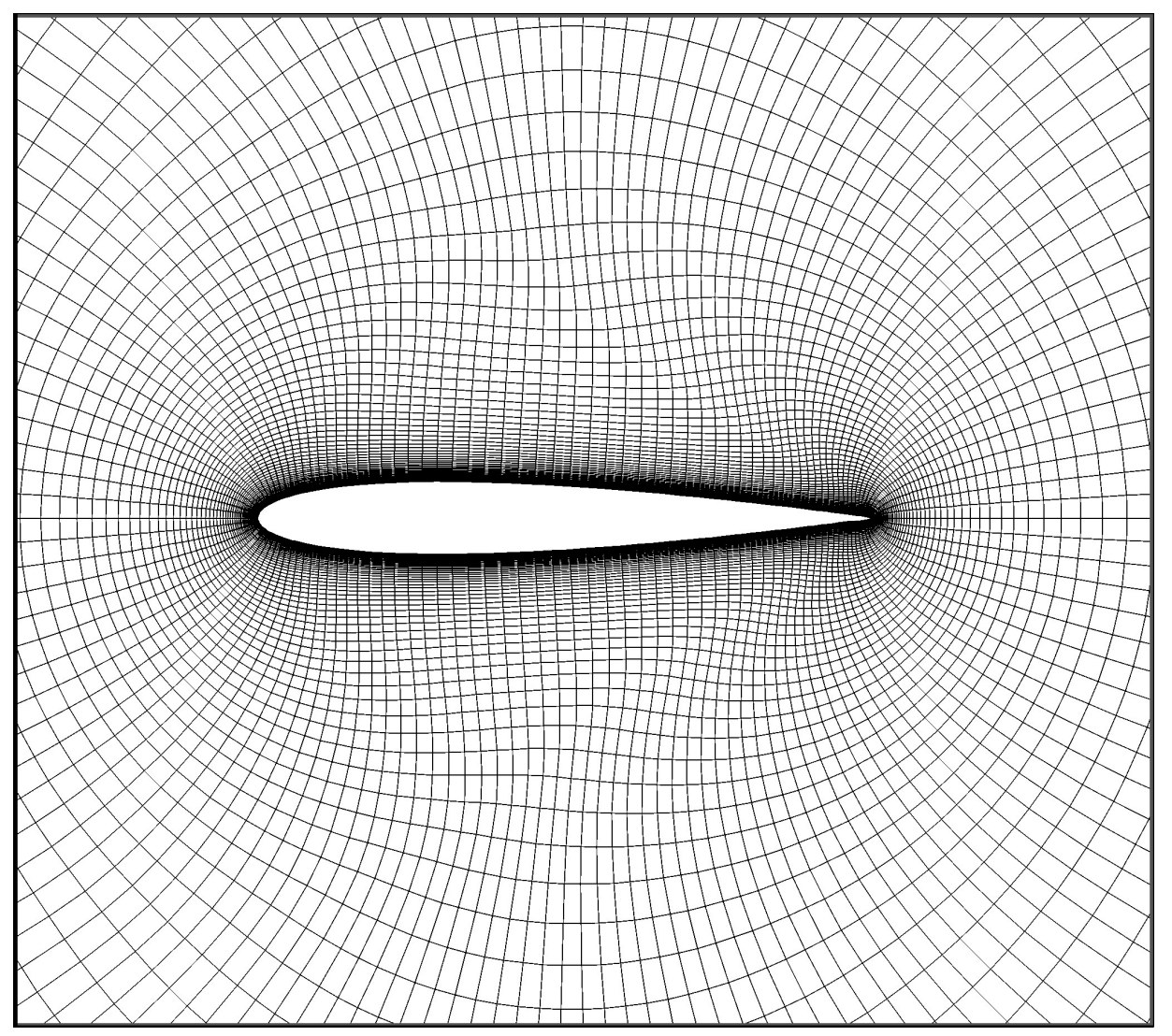

Figure 29: Structured non-uniform grids.

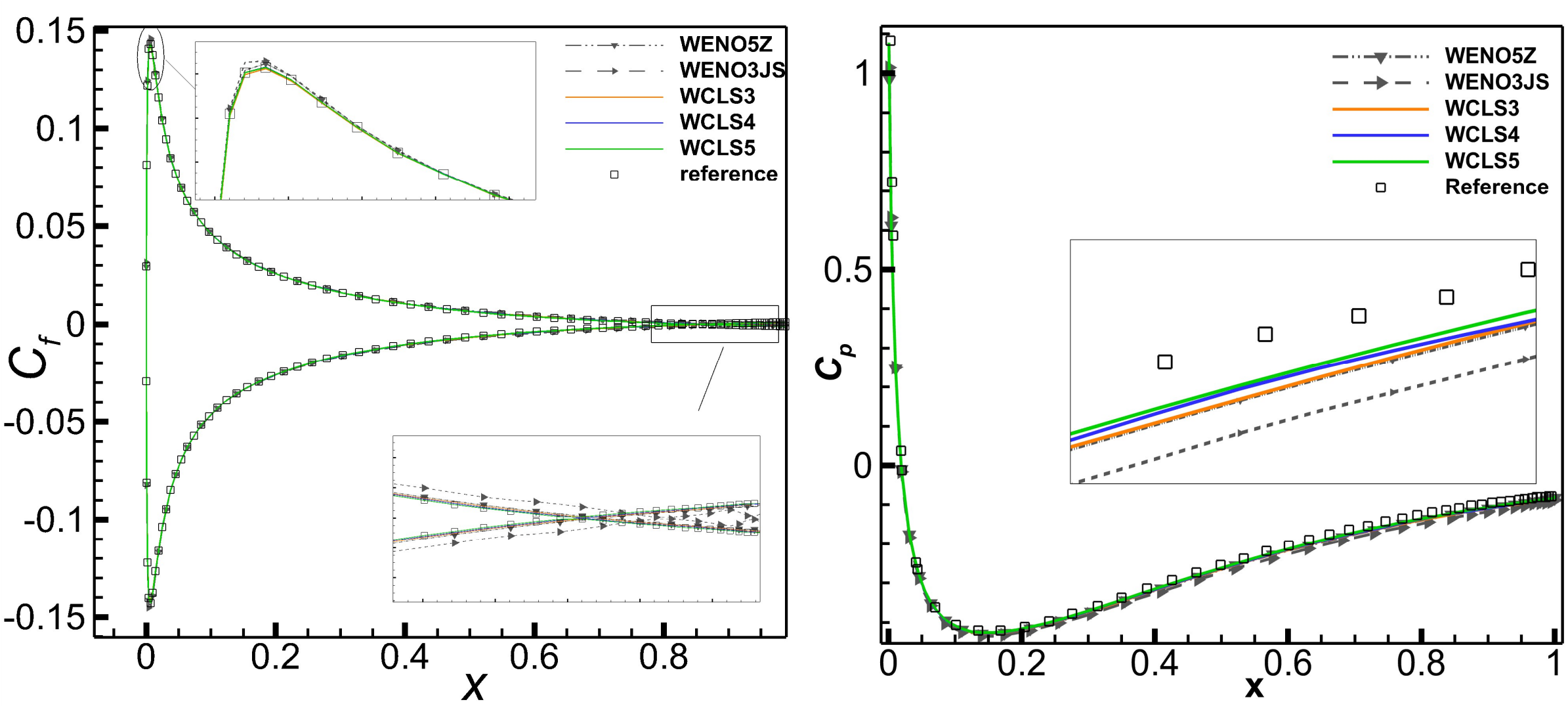

Figure 30: Skin friction coefficient (left) and pressure coefficients (right) on the airfoil obtained by different schemes.

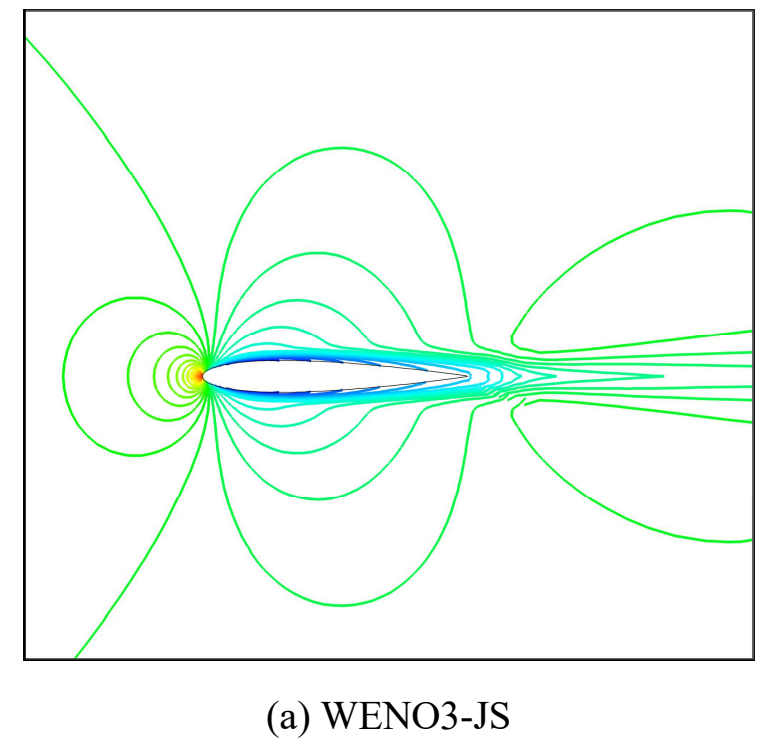

(a) WENO3-JS

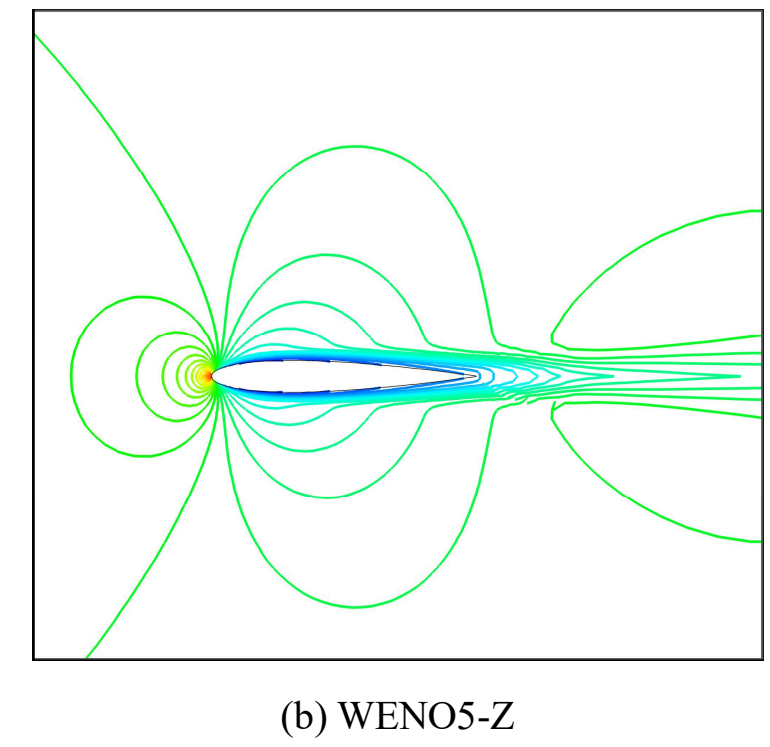

(b) WENO5-Z

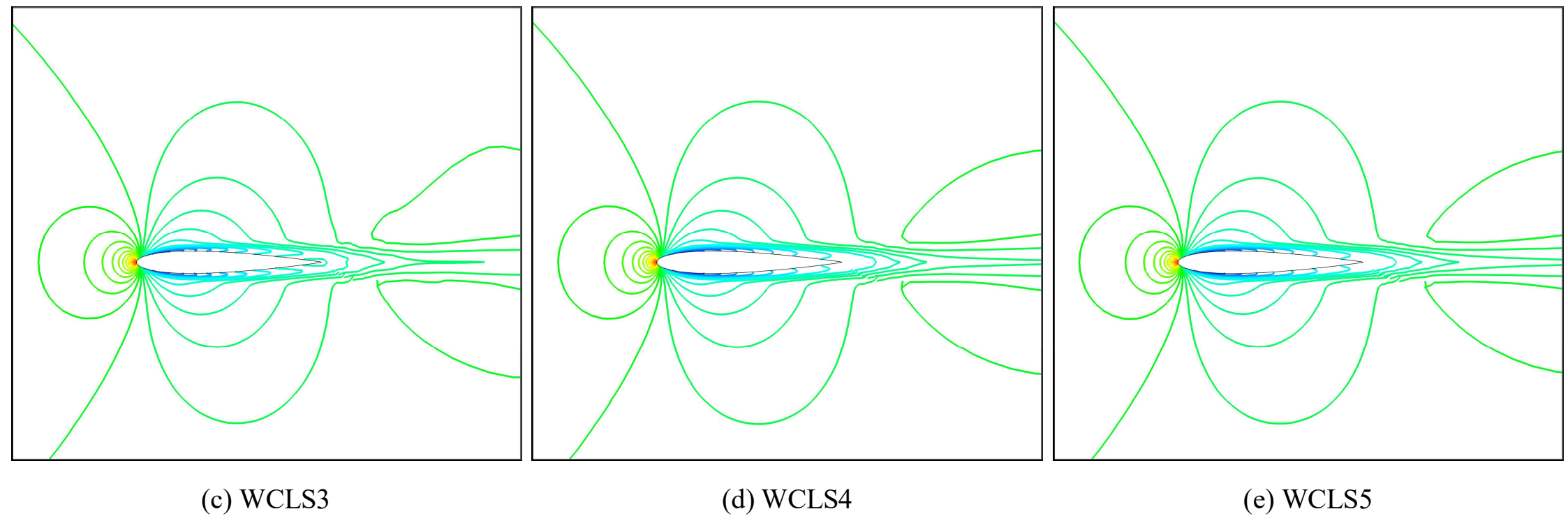

(c) WCLS3 (d) WCLS4 (e) WCLS5

Figure 31: Density contours using different schemes. 30 equally spaced lines from 0.88 to 1.13.

## 4.14 Unsteady Subsonic Viscous Flow around SD7003

An unsteady subsonic viscous flow with Reynolds number 10,000 is simulated. The Mach number is 0.2 and the angle of attack is $4°$. The mesh is C-type with grid number being $148 \times 1422$ and is refined around the wake region shown as in Fig. 32. The chord length is 1. The height of the first layer from the viscous adiabatic wall is 0.00166. The CFL number is 25 with $\epsilon_{\text{NS}} = \epsilon_{\text{Recon}} = 1 \times 10^{-3}$, $N_{\text{inner}} = 25$. The simulation is continued for 20 periods.

The contours of Mach number and $z$-vorticities are shown in Figs. 33 and 34, respectively. The WENO3-JS scheme is the most dissipative scheme, followed by the WENO5-Z scheme. The proposed WCLS schemes are successful in capturing the vortices originating from the separation bubble on the upper-back side of SD7003 airfoil. This test confirms the high resolution of the proposed WCLS schemes on highly non-uniform grids with very large aspect ratio near boundary layer.

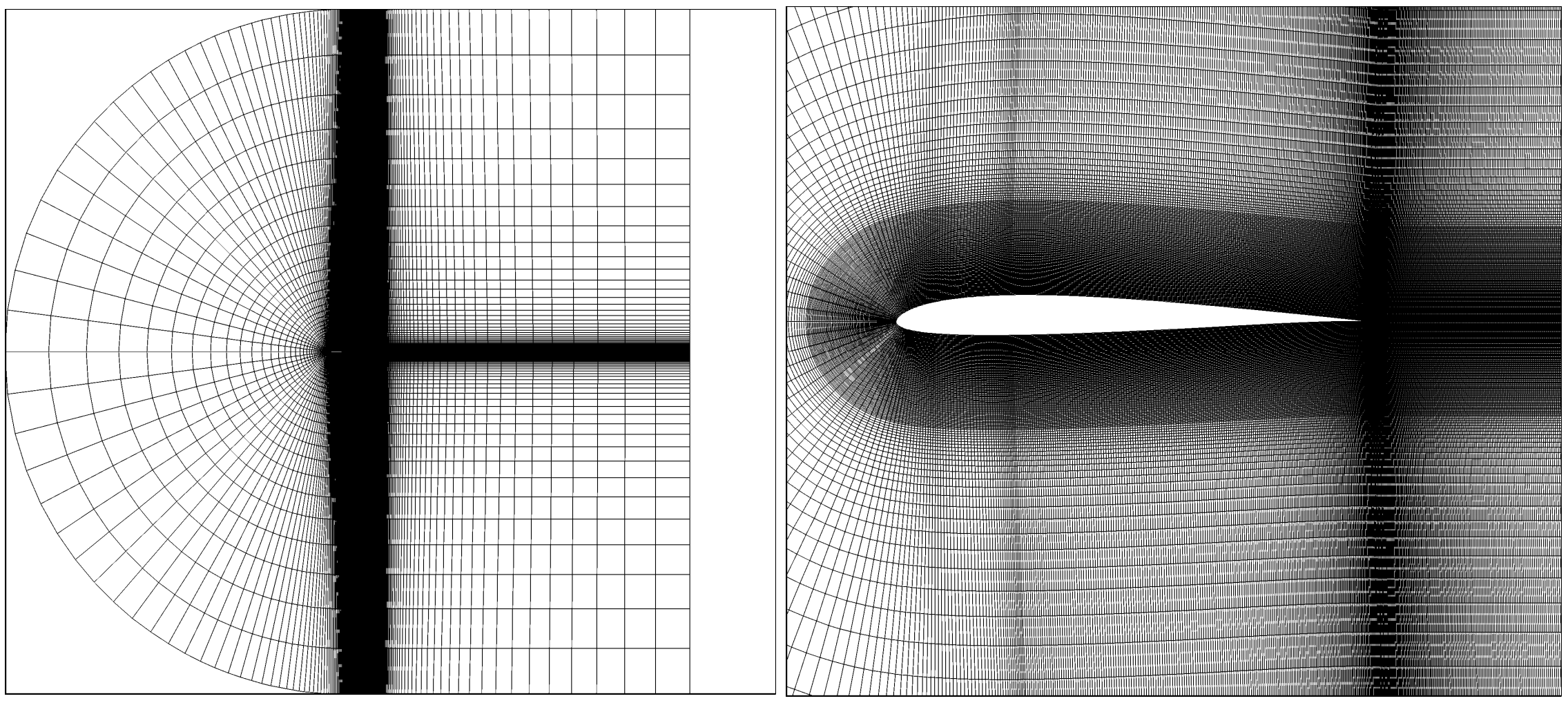

Figure 32: Structured non-uniform mesh.

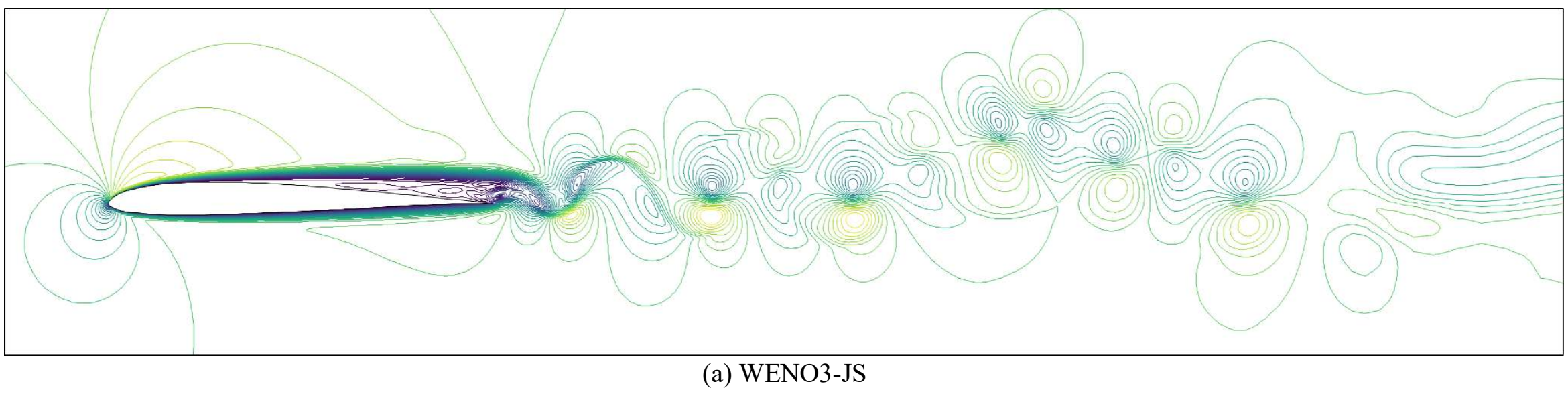

(a) WENO3-JS

(b) WENO5-Z

(c) WCLS3

(d) WCLS4

(e) WCLS5

Figure 33: Mach contours for the flow with Re 10,000 around the SD7003 airfoil using different schemes. 32 equally spaced 0.008 to 0.26.

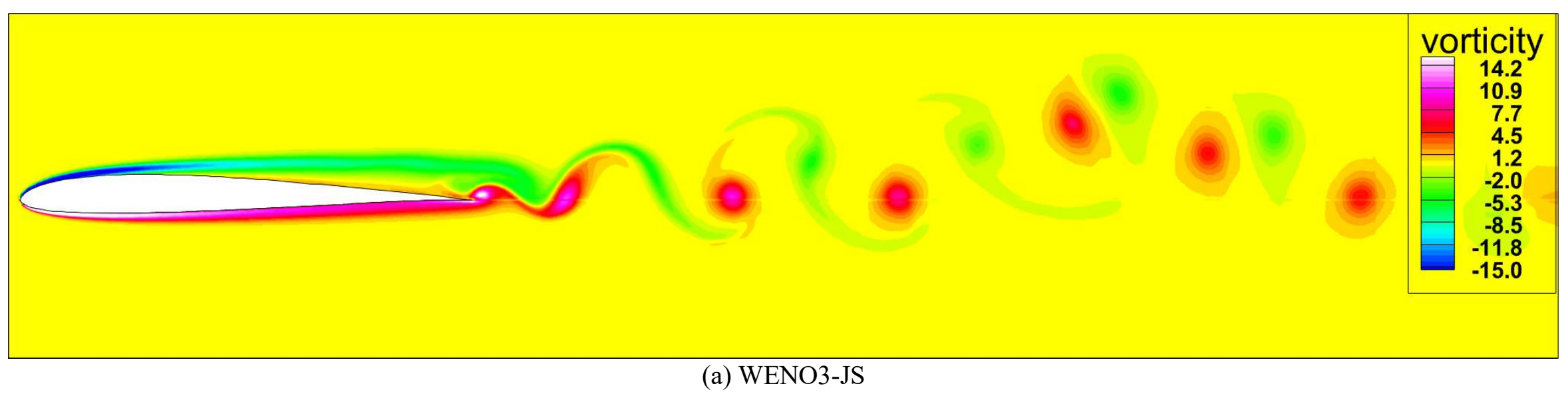

(a) WENO3-JS

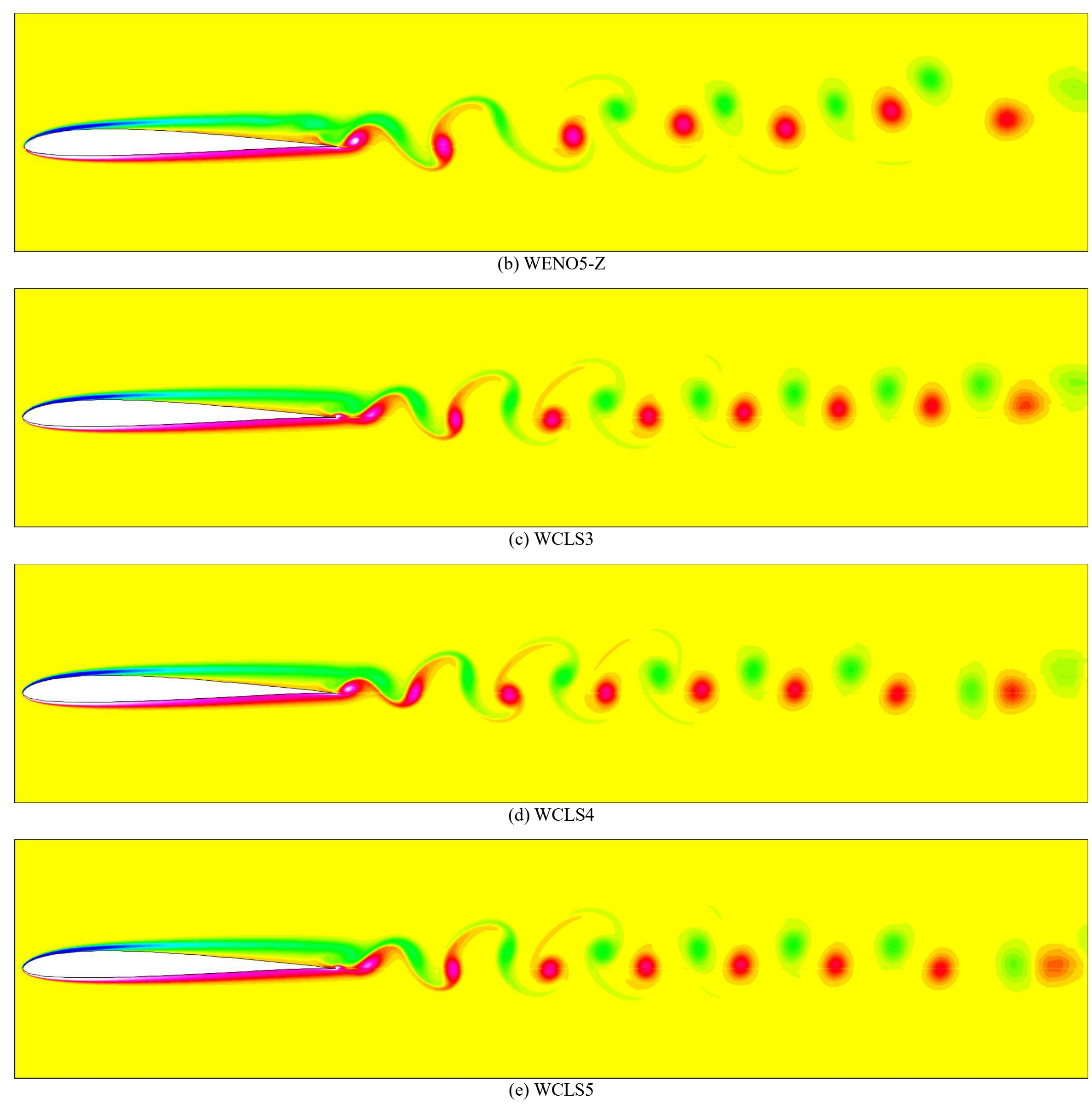

(b) WENO5-Z

(c) WCLS3

(d) WCLS4

(e) WCLS5

Figure 34: Vorticity contours for the flow with Re 10,000 around the SD7003 airfoil with different schemes.

### 4.15 Kelvin-Helmholtz Instability

The Kelvin-Helmholtz instability (KHI) problem is a classic benchmark widely used to evaluate the numerical dissipation of a scheme and its capability to resolve small-scale vortical structures during long-term evolution. Following the setup in [44,45], the initial conditions are given as follows:

$$[\rho, u, v, p] = \begin{cases} \left[1, -0.5 + 0.5e^{\frac{y+0.25}{\delta}}, 0.01\sin(4\pi x), 1.5\right], y \in [-0.5, -0.25), \\ \left[2, 0.5 - 0.5e^{\frac{-y-0.25}{\delta}}, 0.01\sin(4\pi x), 1.5\right], y \in [-0.25, 0), \\ \left[2, 0.5 - 0.5e^{\frac{y-0.25}{\delta}}, 0.01\sin(4\pi x), 1.5\right], y \in [0, 0.25), \\ \left[1, -0.5 + 0.5e^{\frac{0.25-y}{\delta}}, 0.01\sin(4\pi x), 1.5\right], y \in [0.25, 0.5], \end{cases} \tag{69}$$

where $\delta = 0.00625$ corresponds to the shear layer thickness in the simulation. The computational domain is $[-0.5,0.5] \times [-0.5,0.5]$. Periodic boundary conditions are applied in the $x$-direction to allow the spatial development of shear layer, while reflective boundary conditions are imposed in the $y$-direction. The simulations are performed on a uniform grid of $N_x \times N_y = 400 \times 400$ up to a final time of $t = 4.0$. The simulation parameters are $\epsilon_{\text{Euler}} =$

$\epsilon_{\mathrm{Recon}} = 1 \times 10^{-4}$, $N_{\mathrm{inner}} = 6$ and CFL = 1.0.

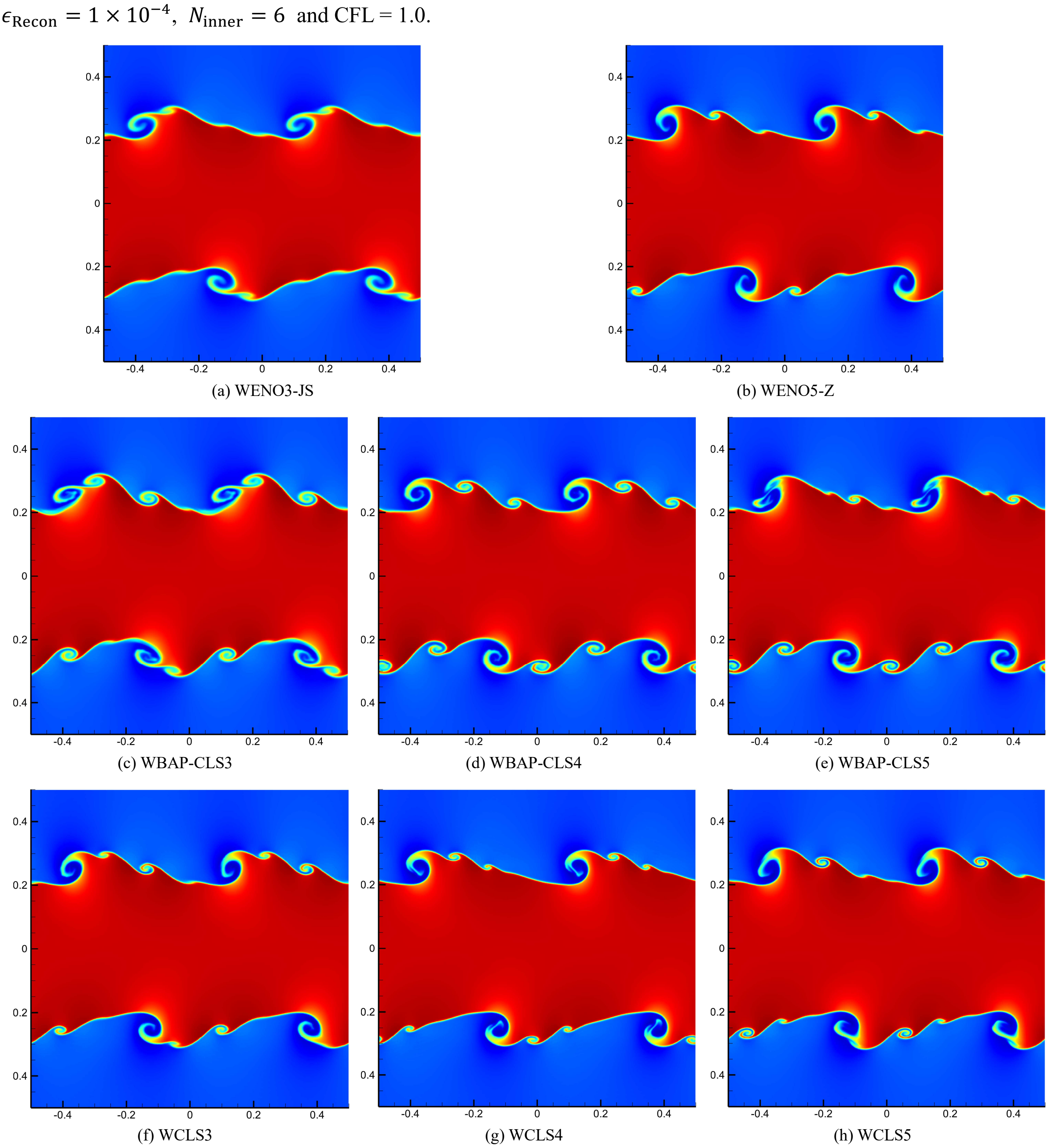

Figure 35: Density contours of the 2D Kelvin-Helmholtz instability computed by different schemes $t$ = 1.0.

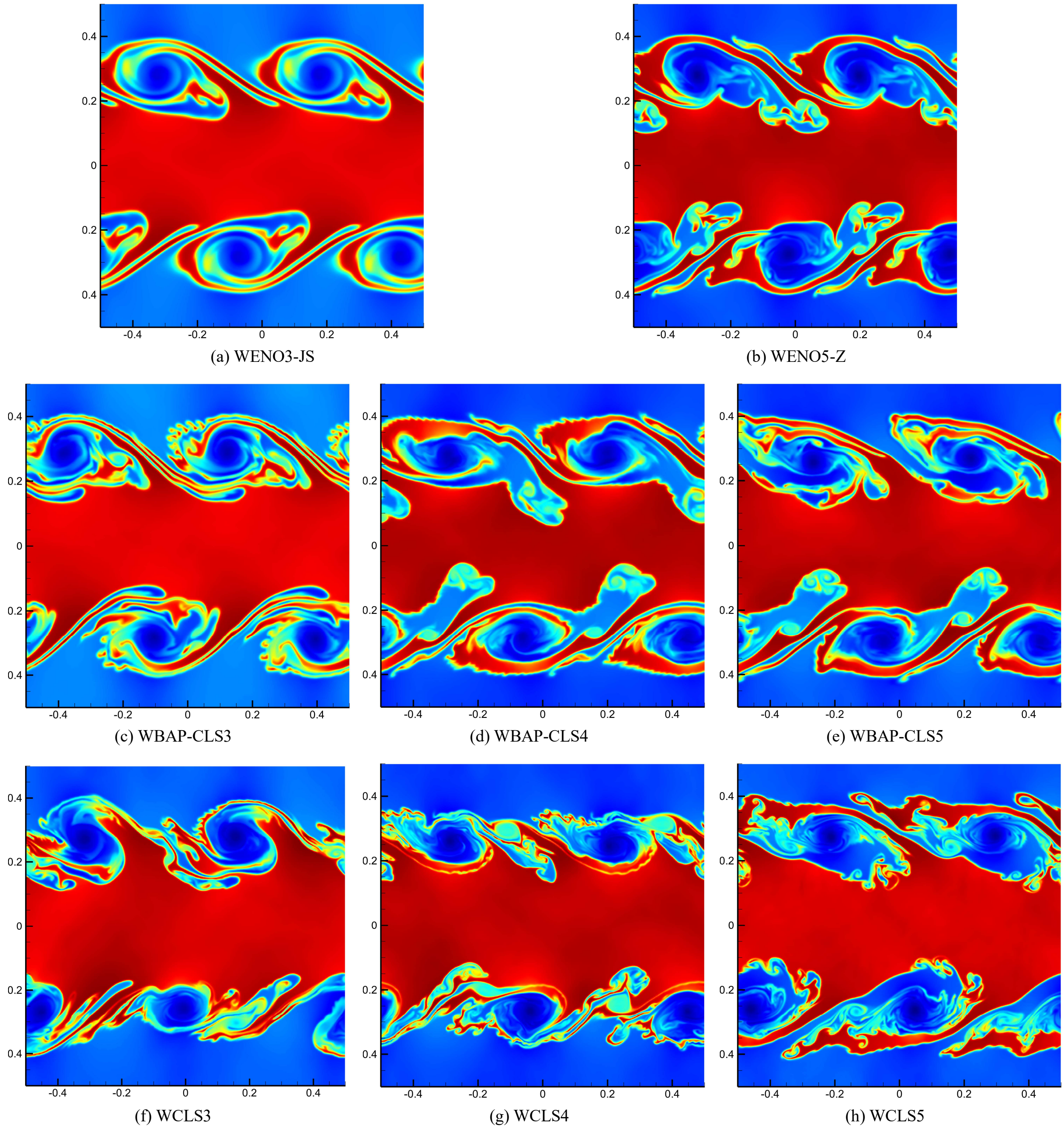

(a) WENO3-JS (b) WENO5-Z

(c) WBAP-CLS3 (d) WBAP-CLS4 (e) WBAP-CLS5

(f) WCLS3 (g) WCLS4 (h) WCLS5

Figure 36: Density contours of the 2D Kelvin-Helmholtz instability computed by different schemes $t = 2.5$.

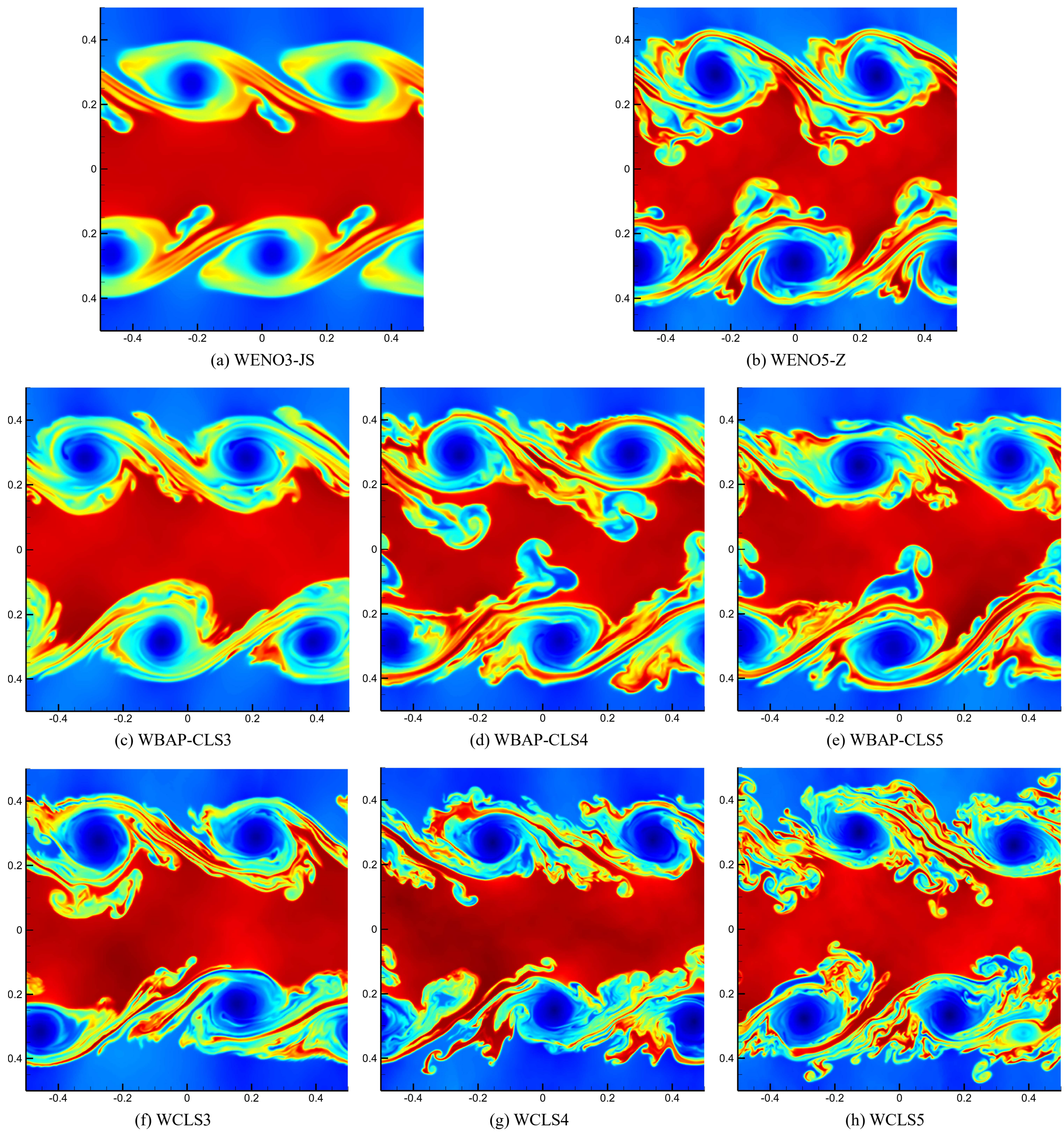

Figure 37: Density contours of the 2D Kelvin-Helmholtz instability computed by different schemes $t = 4.0$.

Figures 35–37 show the density contours at $t = 1.0$, $t = 2.5$, and $t = 4.0$, respectively.

Figure 35 at $t = 1.0$ shows that the proposed WCLS schemes capture the primary roll-up vortices and thin shear layers more sharply than the WENO3-JS and WENO5-Z schemes. The WCLS schemes also resolve the subtle interfaces of the Kelvin-Helmholtz billows more clearly than the corresponding WBAP-CLS schemes, even at this early stage.

As the flow evolves to $t = 2.5$ and $t = 4.0$, the differences in numerical dissipation among these schemes become more evident. Because of its excessive numerical dissipation, the traditional WENO3-JS scheme completely smears out the secondary instabilities. In contrast, the proposed WCLS schemes reveal finer structures in the complex wave-breaking patterns. As shown in Figs. 36 and 37, the WCLS schemes produce more intricate small-scale vortical structures than the WBAP-CLS and WENO schemes. These results indicate that the proposed WCLS method resolves small-scale instability-related phenomena more effectively while maintaining low numerical dissipation during the long-term evolution of compressible turbulent flows.

### 4.16 3D Compressible Inviscid Taylor-Green Vortex

To further evaluate the capability of the proposed high-order schemes in capturing small-scale instable structures during long-time evolution, the three-dimensional compressible inviscid Taylor-Green vortex (TGV) flow is considered. The TGV problem is widely used for assessing the dissipation performance and resolving efficiency of numerical schemes. The simulations are conducted in a cubic domain of $\Omega = [0,2\pi] \times [0,2\pi] \times [0,2\pi]$ with periodic boundary conditions in all three directions. The simulation parameters are $\epsilon_{\text{Euler}} = \epsilon_{\text{Recon}} = 1 \times 10^{-4}$, $N_{\text{inner}} = 6$ and CFL = 2.0. The initial conditions for this purely inviscid flow are given by [46,47] as

$$\begin{cases} \rho_0(x,y,z) = 1, \\ u_0(x,y,z) = \sin(x)\cos(y)\cos(z), \\ v_0(x,y,z) = -\cos(x)\sin(y)\cos(z), \\ w_0(x,y,z) = 0, \\ p_0(x,y,z) = 100 + \dfrac{1}{16}[(\cos(2z)+2)(\cos(2x)+\cos(2y))-2], \end{cases} \tag{70}$$

To quantify the low-dissipation property of the proposed WCLS and CLS methods, we monitor the temporal evolution of the volume-averaged kinetic energy, defined as [48]

$$E_k(t) = \frac{1}{\rho_0|\Omega|}\int_\Omega \frac{1}{2}\rho(u^2+v^2+w^2)\mathrm{d}\Omega, \tag{71}$$

where $|\Omega|$ is the volume of the computational domain.

For an inviscid flow, the kinetic energy should theoretically be conserved with an ideal constant value of $E_k = 0.125$. The numerical decay of $E_k(t)$ serves as an indicator of the inherent numerical dissipation introduced by different schemes. The simulation parameters are $\epsilon_{\text{Euler}} = \epsilon_{\text{Recon}} = 1 \times 10^{-4}$, $N_{\text{inner}} = 6$ and CFL = 2.0.

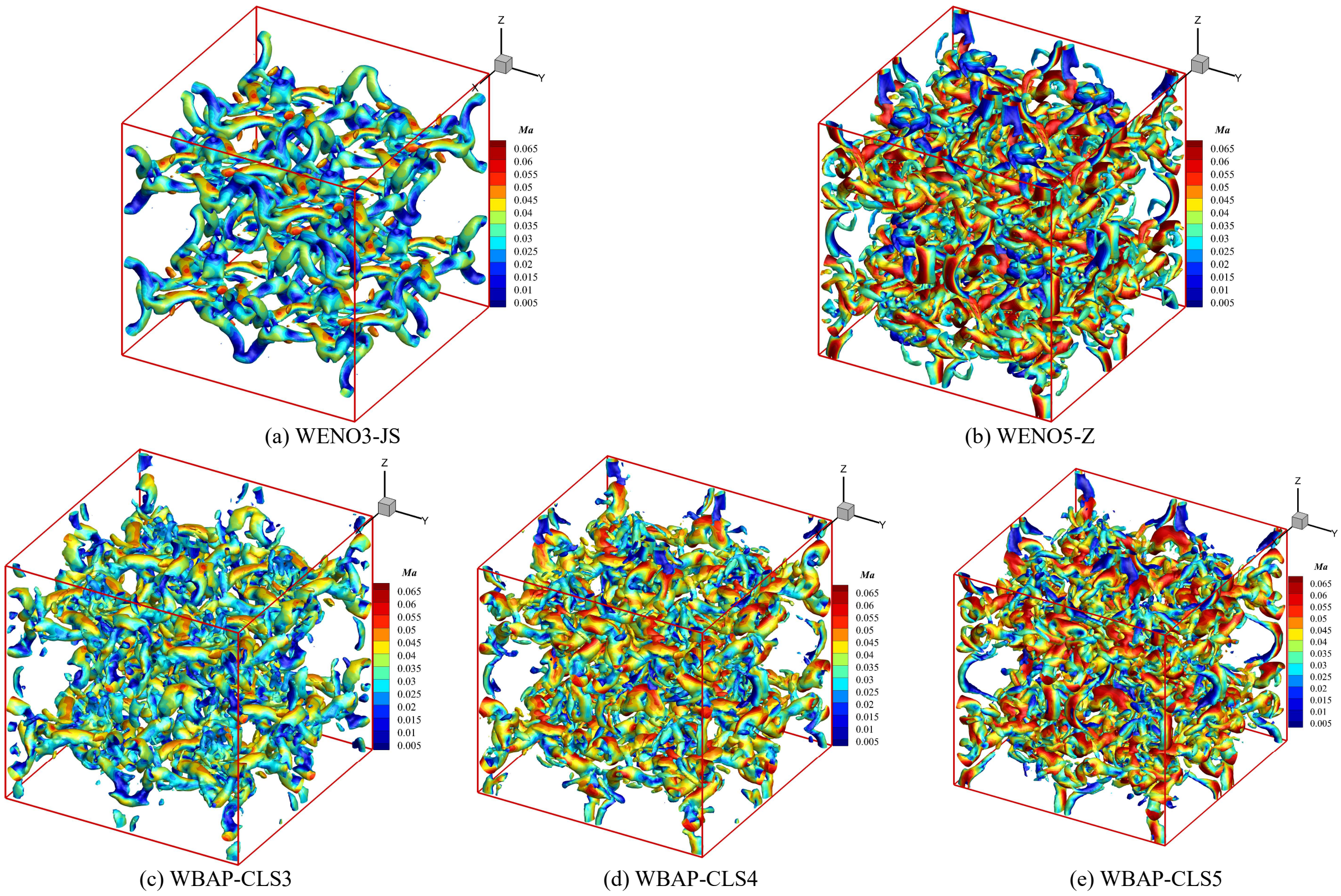

(a) WENO3-JS (b) WENO5-Z

(c) WBAP-CLS3 (d) WBAP-CLS4 (e) WBAP-CLS5

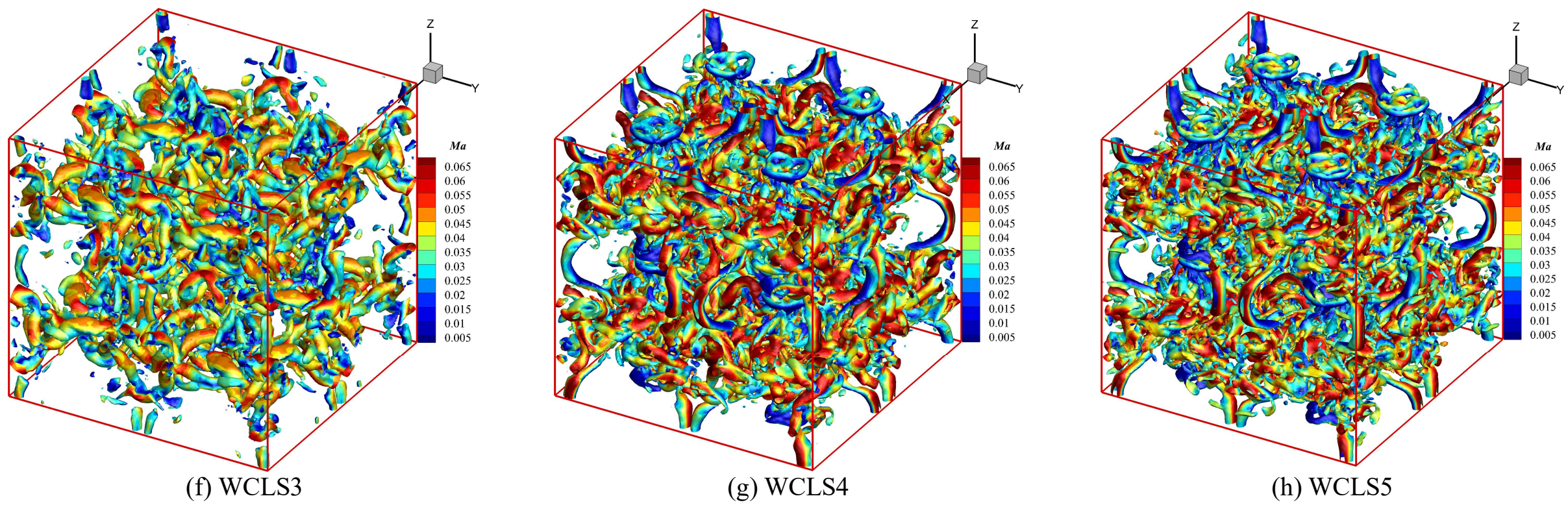

Figure 38: Iso-contours (Q = 1.5) obtained for the 3D inviscid TGV problem using different schemes at t = 10.

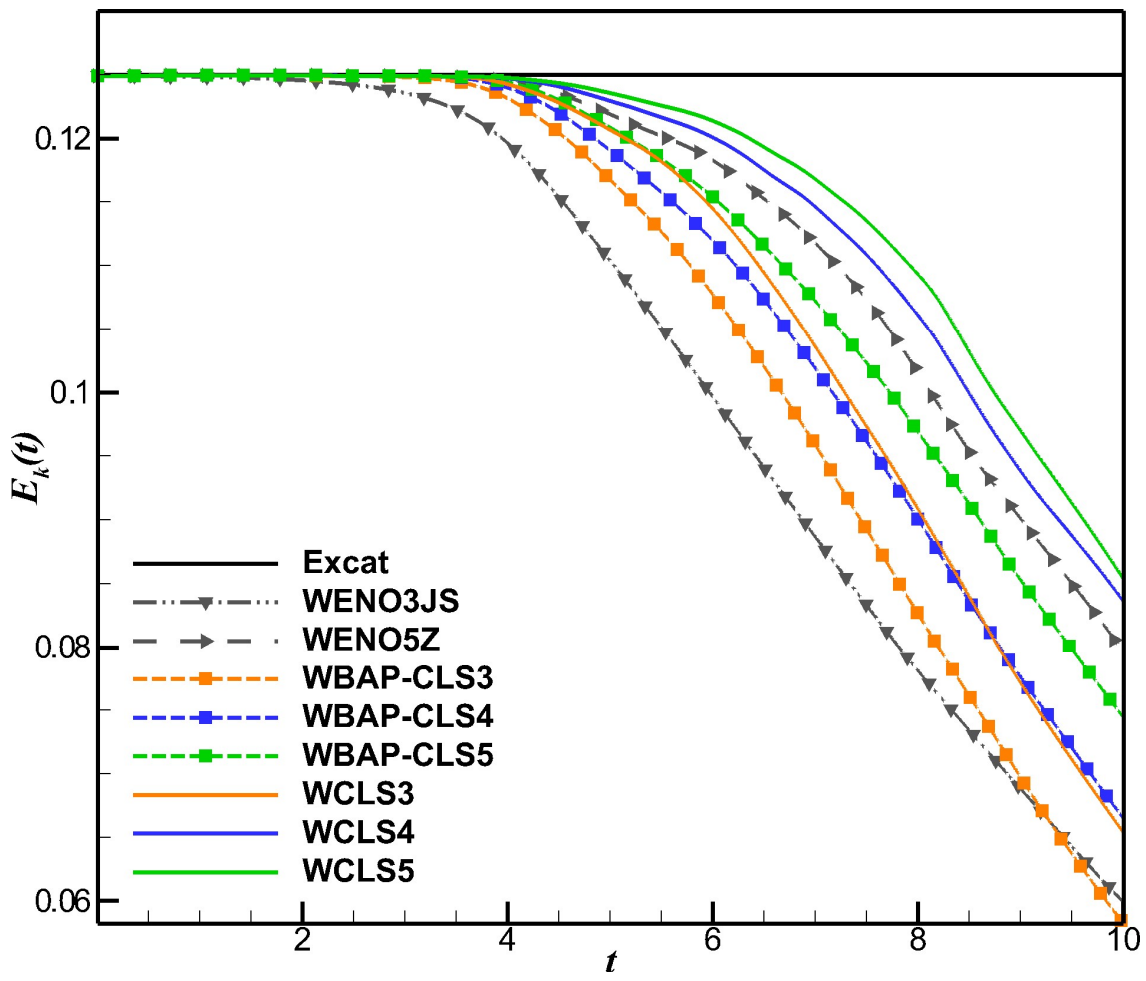

Figure 39: The evolution of kinetic energy versus dimensionless time obtained by different schemes for 3D inviscid TGV flows.

Figure 38 displays the Q-criterion iso-surfaces (Q = 1.5, colored by Mach number) for assessing the vortex-capturing capability of the tested schemes. While all methods capture the primary vortical structures, the WENO3-JS scheme introduces excessive numerical dissipation and smears the fine-scale details. The WBAP-CLS schemes exhibit improved resolution; however, the proposed WCLS schemes resolve more intricate vortical structures. In particular, WCLS5 achieves the highest fidelity, capturing the finest small-scale vortices and thereby demonstrating superior spatial accuracy for this complex flow field.

To further assess the dissipative properties, Fig. 39 plots the temporal evolution of the kinetic energy. The kinetic-energy decay rates associated with the original WENO and WBAP-CLS schemes are higher than those of the proposed WCLS schemes. Specifically, WCLS5 exhibits the lowest dissipation rate and retains the highest kinetic energy at the final simulation time. This confirms that the WCLS formulation reduces numerical dissipation relative to the reference schemes and therefore provides better energy preservation.

Overall, the 3D TGV simulations confirm the high accuracy and low dissipation of the proposed WCLS schemes. Compared with the traditional WENO and WBAP-CLS methods, the WCLS approach offers better overall performance in resolving small-scale vortical structures and preserving kinetic energy.

### 4.17 3D Compressible Viscous Taylor-Green Vortex [49,50]

To further illustrate the extended application capability of the proposed WCLS schemes for viscous turbulent problems, the 3D compressible viscous Taylor-Green vortex (vis-TGV) flow is considered. The evolution of small-scale vortical structures in the vis-TGV flow differs from the inviscid case, as its turbulent features are governed by two dimensionless parameters: the Reynolds number (Re) and the reference Mach number ($\mathrm{Ma_{ref}}$). In this study, two representative cases are simulated: Case I with $\mathrm{Re} = 1600$ and $\mathrm{Ma_{ref}} = 1.25$, and Case II with $\mathrm{Re} = 400$ and $\mathrm{Ma_{ref}} = 2.0$. The simulation parameters are $\epsilon_{\mathrm{Euler}} = \epsilon_{\mathrm{Recon}} = 1 \times 10^{-4}$, $N_{\mathrm{inner}} = 6$ and CFL = 2.0.

Case I combines a moderately high Reynolds number with significant compressibility, leading to complex interactions between fine-scale turbulence and localized shocklets. It therefore provides a stringent test of the scheme's ability to balance high-resolution vortex capture with robust shock capturing. By contrast, Case II is more strongly influenced by viscosity but also features stronger compressibility and stronger shocks, making it a demanding test of robustness at higher Mach numbers.

Computations are carried out on a uniform grid of $128^3$ control volumes. The computational domain, periodic boundary conditions, and the initial velocity field $[u_0, v_0, w_0]$ are identical to those used in the 3D inviscid TGV simulations. For the thermodynamic variables, the initial temperature is set to $T_0 = 1$, and the initial density and pressure are given by

$$\begin{cases} \rho_0(x,y,z) = 1 + \dfrac{c_f}{16}[\cos(2x) + \cos(2y)][\cos(2z) + 2], \\ p_0(x,y,z) = \dfrac{\rho_0 T_0}{\gamma \mathrm{Ma}_{\mathrm{ref}}^2}, \end{cases} \tag{72}$$

where $c_f = \gamma \mathrm{Ma}_{\mathrm{ref}}^2$ for Case I, and $c_f = 1$ for Case II. Temporal evolution of kinetic energy and the total viscous dissipation rate are monitored. The total viscous dissipation rate, denoted as $\varepsilon^T$, consists of both solenoidal and dilatational dissipation components, which correspond to the vortex dynamics and compressibility effects, respectively [51]:

$$\varepsilon^T = \varepsilon^S + \varepsilon^D = \frac{1}{\rho_{\mathrm{ref}} \mathrm{Re}|\Omega|} \int_\Omega \mu(\nabla \times \boldsymbol{u})^2 \, \mathrm{d}\Omega + \frac{4}{3\rho_{\mathrm{ref}} \mathrm{Re}|\Omega|} \int_\Omega \mu(\nabla \cdot \mathbf{u})^2 \, \mathrm{d}\Omega. \tag{73}$$

**Case I:** Re = 1600, $\mathrm{Ma}_{\mathrm{ref}} = 1.25$. Under the effects of relatively small physical viscosity and a supersonic reference Mach number, the interactions between localized shocklets and turbulent vortices occur, eventually breaking down into intricate small-scale turbulent structures. The 3D iso-contours of Q-criterion colored by velocity $u$ at time $t = 18$ are displayed in Fig. 40. The time evolution of kinetic energy $E_k$, total viscous dissipation rate $\varepsilon^T$, and the kinetic energy spectrum $\mathcal{E}(k)$ at $t = 20s$ are shown in Fig. 41.

As shown in Fig. 40, the WENO3-JS and WENO5-Z schemes exhibit excessive numerical dissipation, leaving only sparse, large-scale vortex rings with fine details smeared out. The proposed WCLS schemes capture substantially richer and more intricate fine-scale turbulent structures than the WENO and WBAP-CLS schemes. Among them, WCLS5 shows the richest small-scale structures. This visual evidence indicates that the proposed WCLS schemes are particularly effective in resolving compressible turbulent vortices.

The quantitative assessments in Fig. 41 further corroborate these qualitative observations. The evolutions of kinetic energy $E_k$ exhibit relatively small but visible differences among the schemes. However, the curves of the total viscous dissipation rate ($\varepsilon^T$), which are highly sensitive to numerical dissipation, reveal profound disparities. The WENO schemes underpredict the dissipation peak. In contrast, the proposed WCLS schemes yield a higher $\varepsilon^T$ peak compared to both the WBAP-CLS and WENO schemes, indicating lower numerical dissipation and a stronger ability to preserve un-damped physical viscous behaviors. Furthermore, as depicted in the kinetic energy spectrum $\mathcal{E}(k)$, the WCLS scheme maintains the turbulent energy cascade much closer to the theoretical Kolmogorov $-5/3$ inertial subrange at high wavenumbers, whereas the WENO3-JS scheme dissipates energy earlier than the other tested schemes.

In summary, the proposed WCLS schemes combine low dissipation with high accuracy for compressible viscous turbulent flows.

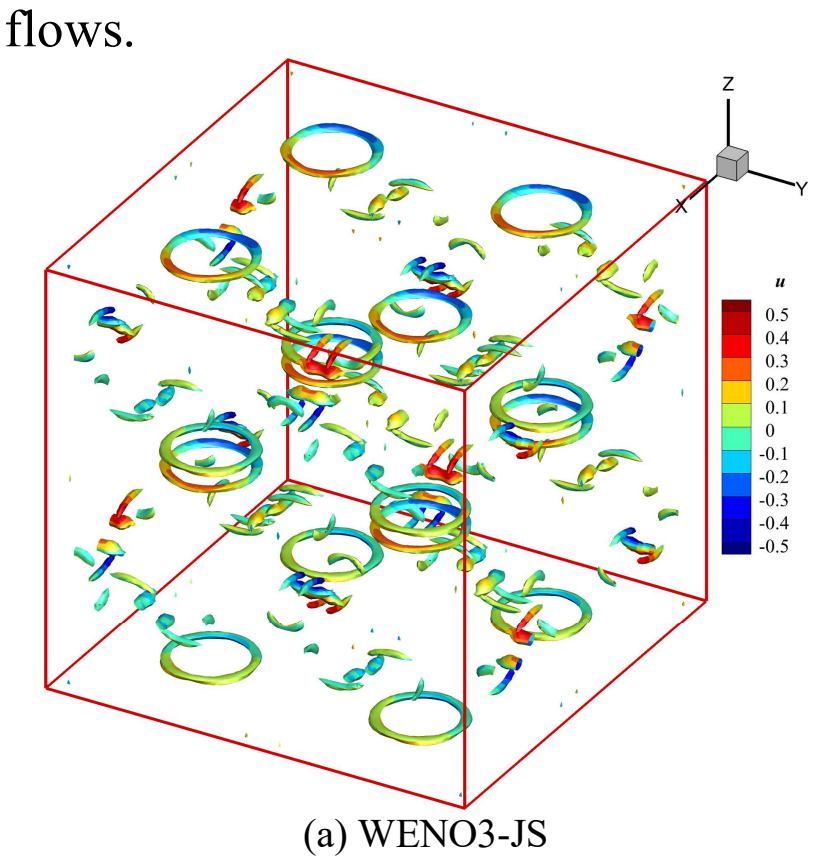

(a) WENO3-JS

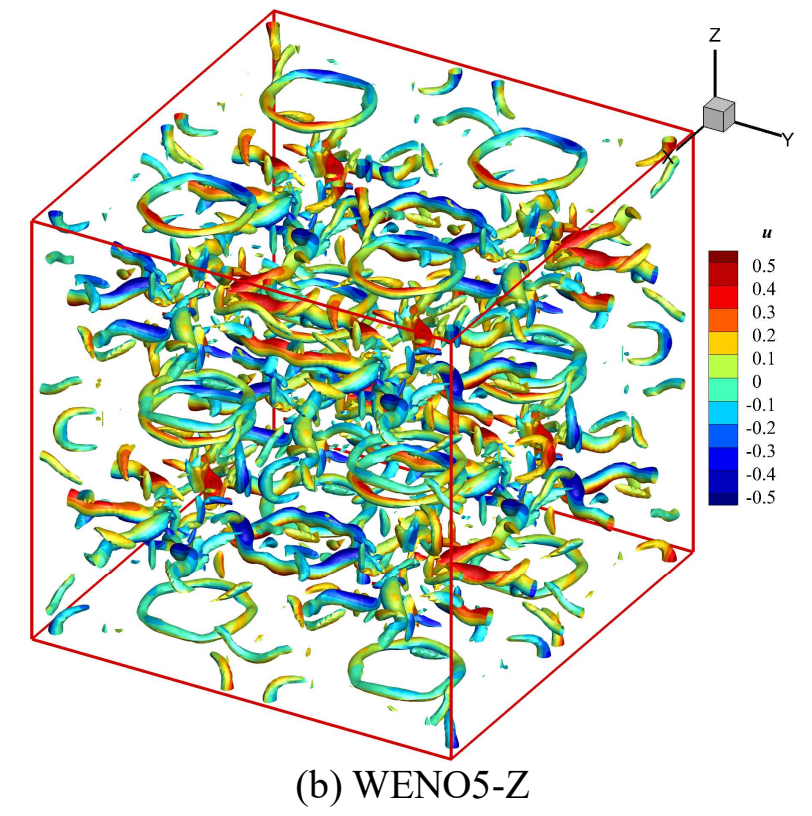

(b) WENO5-Z

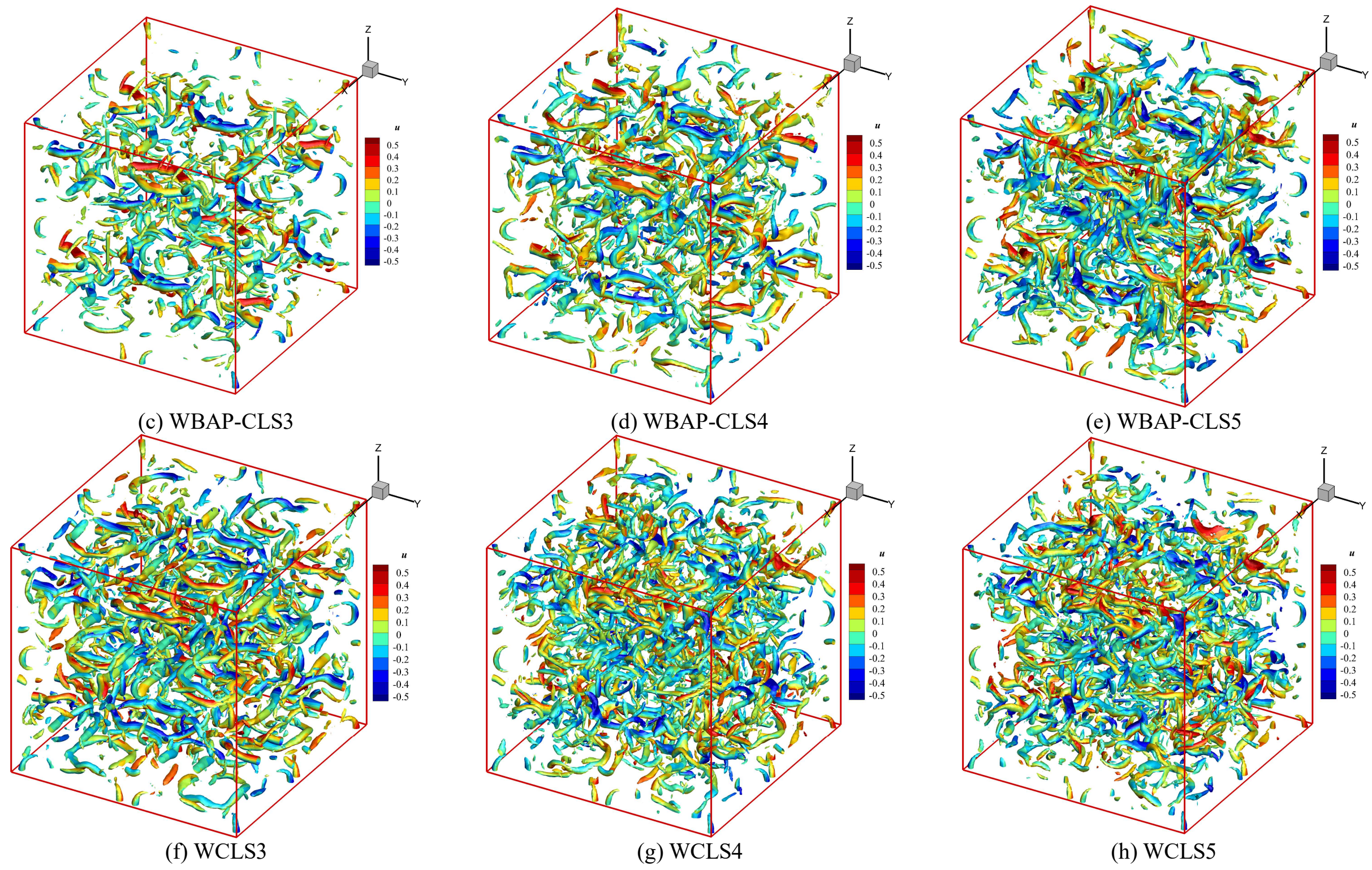

Figure 40: Iso-surfaces of Q = 2.5, colored by the u-velocity, obtained by different schemes for the viscous TGV problem.

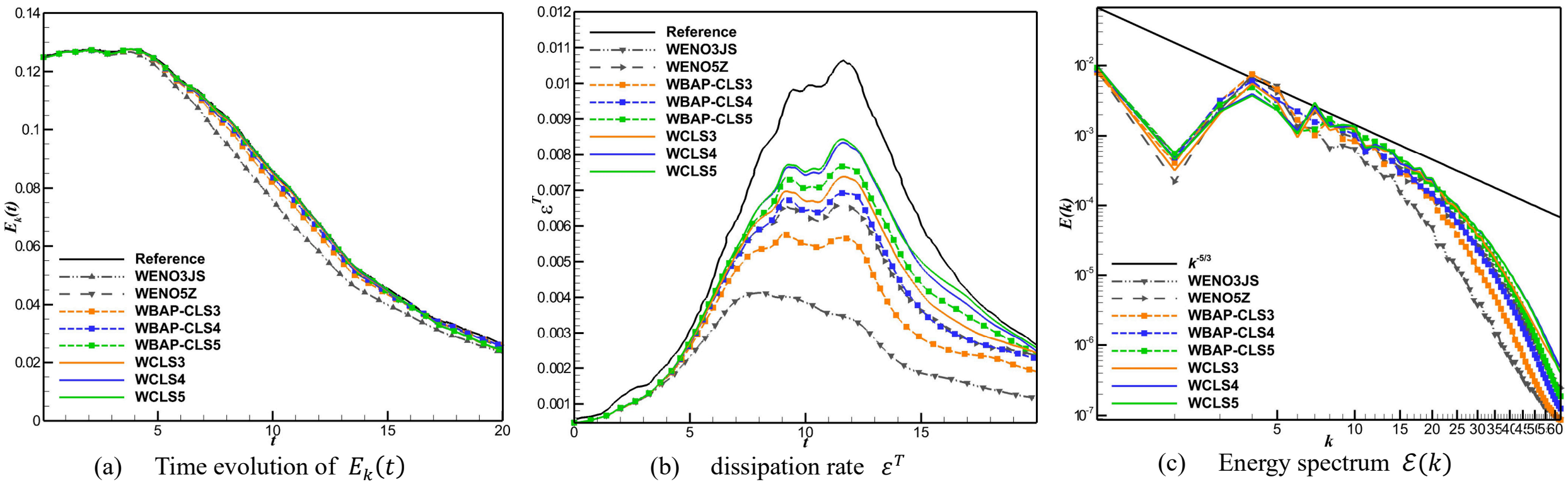

Figure 41: Comparison of different schemes for the viscous TGV problem. (a) Temporal evolution of turbulent kinetic energy; (b) temporal evolution of the turbulent dissipation rate; (c) energy spectrum at t = 20 s.

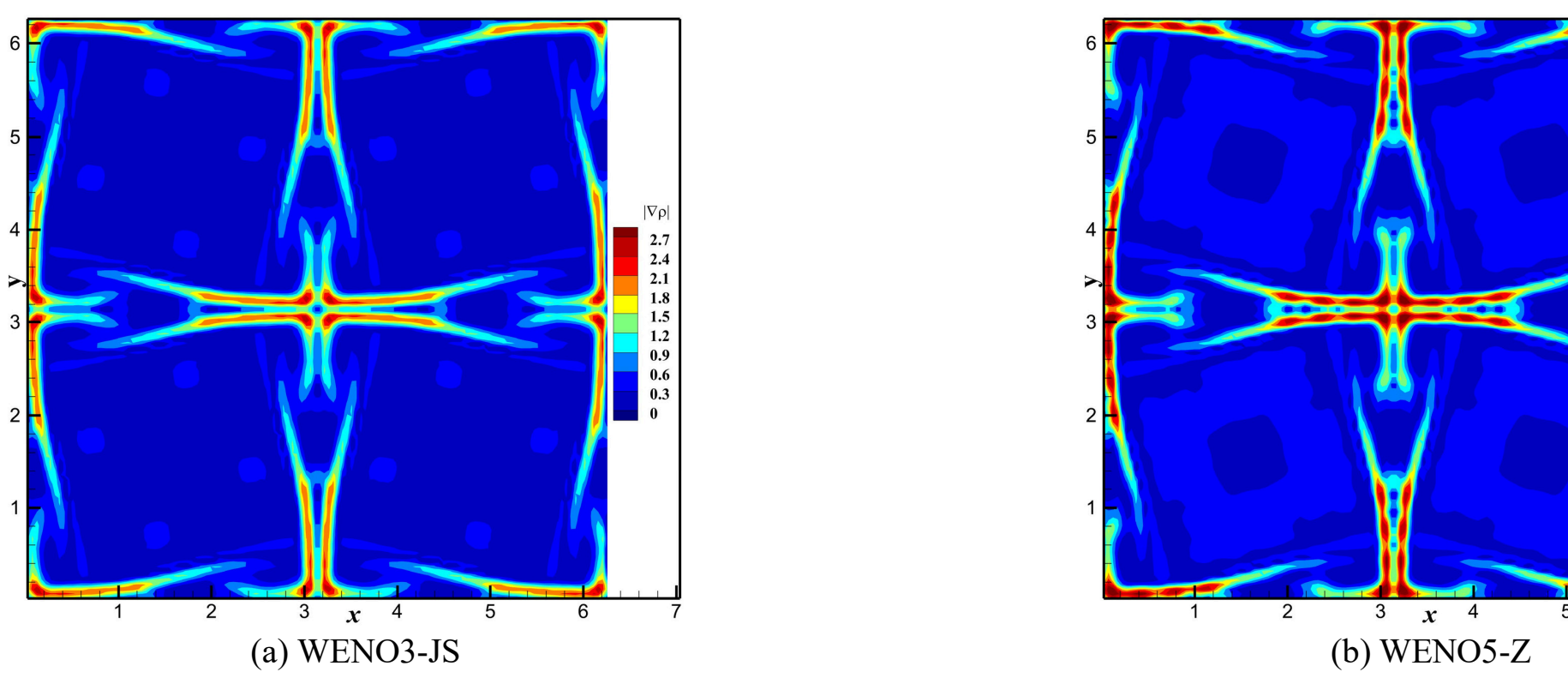

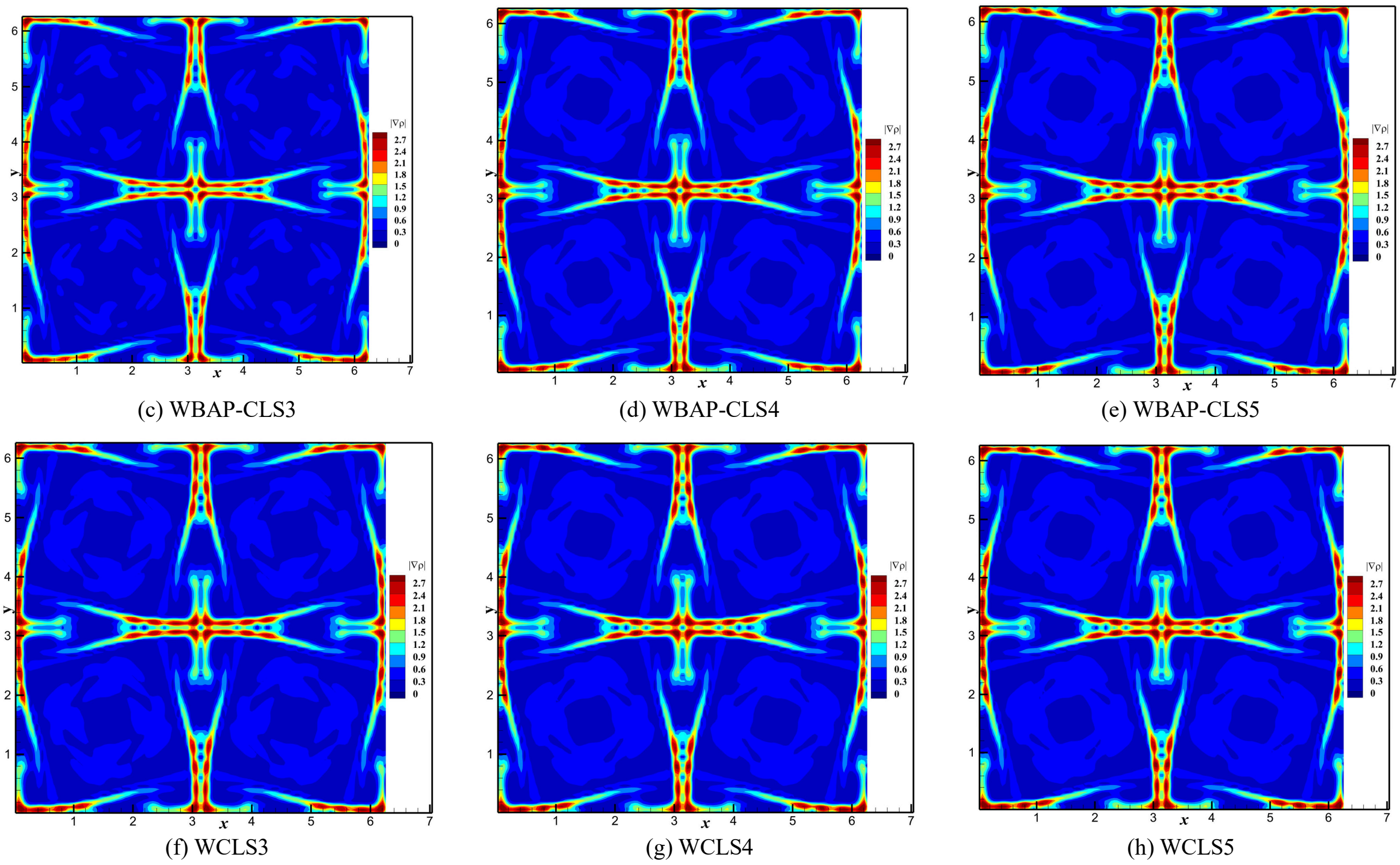

(c) WBAP-CLS3 (d) WBAP-CLS4 (e) WBAP-CLS5

(f) WCLS3 (g) WCLS4 (h) WCLS5

Figure 42: Cross-sectional contours of the density-gradient norm at $z = \pi$ obtained by different schemes for the viscous TGV problem at $t$ = 7.0.

To further investigate the capability of the schemes in resolving compressibility effects and shock-vortex interactions, the cross-sectional contours of the density gradient norm ($|\nabla\rho|$) at $z = \pi$ and an intermediate time $t = 7.0$ are compared in Fig. 42. At this stage, the supersonic flow ($\mathrm{Ma}_{\mathrm{ref}} = 1.25$) generates localized shocklets induced by the strong vortical motions. As observed in Fig. 42, the proposed WCLS schemes provide improved resolution among the tested schemes, not only maintain sharp and thin shock interfaces but also successfully resolve the faint acoustic waves and density fluctuations trapped within the vortex cores.

**Case II:** Re = 400, $\mathrm{Ma}_{\mathrm{ref}}$ = 2.0. This case represents a highly compressible and viscous-dominated flow regime. The 3D iso-contours of the vortical structures at $t$ = 10$s$ is presented in Fig. 43. The time evolution of kinetic energy $E_k$, total viscous dissipation rate $\varepsilon^T$, and the kinetic energy spectrum $\mathcal{E}(k)$ are presented in Fig. 44.

From the 3D visualizations in Fig. 43, the macroscopic vortical structures obtained by different schemes appear almost identical, featuring relatively thick and coherent vortex tubes. This visual similarity is physically reasonable: at a low Reynolds number of $\mathrm{Re} = 400$, the strong physical viscosity dominates the flow field, inherently smoothing out the ultra-fine turbulent structures and masking the visual differences in numerical dissipation among the schemes.

Due to the high physical viscosity, the evolution of turbulent kinetic energy $E_k$ displays a faster overall decay rate and all schemes perform similarly on this macroscopic scale. However, the total viscous dissipation rate $\varepsilon^T$ reveals clear differences. Driven by strong compressibility at $\mathrm{Ma}_{\mathrm{ref}} = 2.0$, the $\varepsilon^T$ curves exhibit multiple distinct peaks. The traditional WENO schemes fail to resolve these peaks, excessively suppressing the physical dissipation due to the numerical damping in their shock-capturing mechanisms. Conversely, the proposed WCLS schemes capture the highest $\varepsilon^T$ peaks among the tested schemes. Furthermore, the kinetic energy spectrum $\mathcal{E}(k)$ confirms that the WCLS scheme maintains slightly higher turbulent kinetic energy across the wavenumber range before the viscous cut-off. These quantitative results suggest that the proposed WCLS schemes maintain high robustness and low numerical dissipation in high-Mach shock-dominated viscous flows.

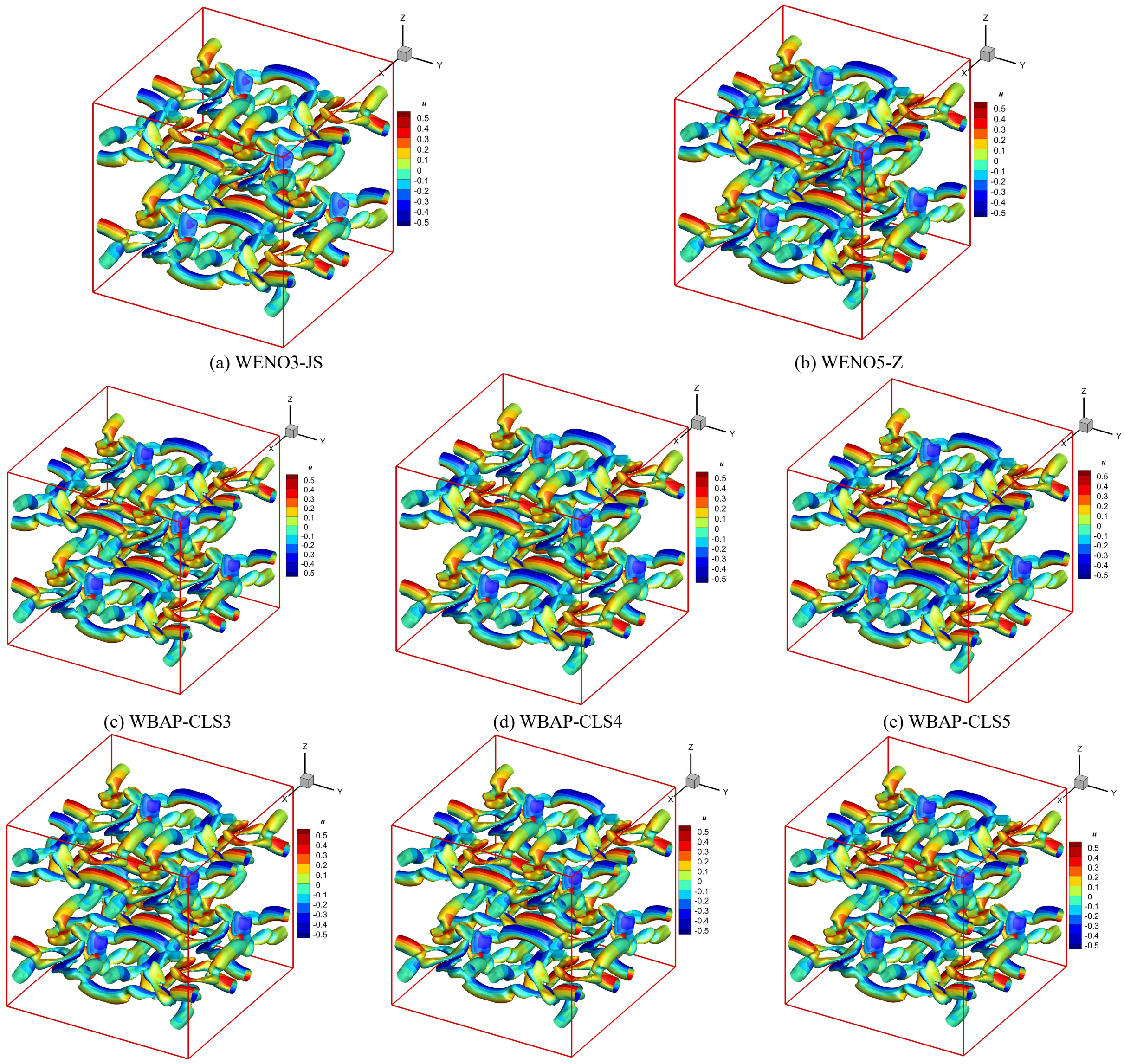

(a) WENO3-JS

(b) WENO5-Z

(c) WBAP-CLS3

(d) WBAP-CLS4

(e) WBAP-CLS5

(f) WCLS3

(g) WCLS4

(h) WCLS5

Figure 43: Iso-surfaces of Q = 0.75, colored by the x-velocity, obtained by different schemes for the viscous TGV problem of Case II.

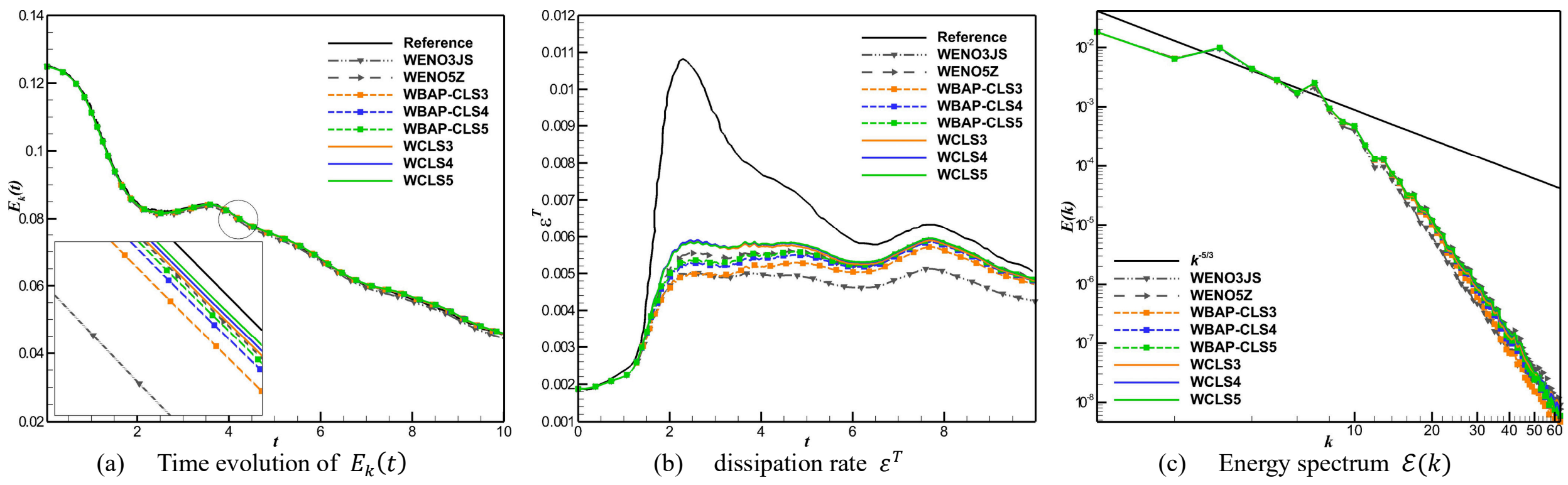

(a) Time evolution of $E_k(t)$ (b) dissipation rate $\varepsilon^T$ (c) Energy spectrum $\mathcal{E}(k)$

Figure 44: Comparison of different schemes for the viscous TGV problem of Case II at $t$ = 10. (a) Temporal evolution of $E_k(t)$; (b) temporal evolution of $\varepsilon^T$; (c) energy spectrum $\mathcal{E}(k)$.

Figure 45 shows the normalized computational time for various unsteady cases, further demonstrating the high efficiency of the proposed WCLS schemes within the implicit time-integration framework.

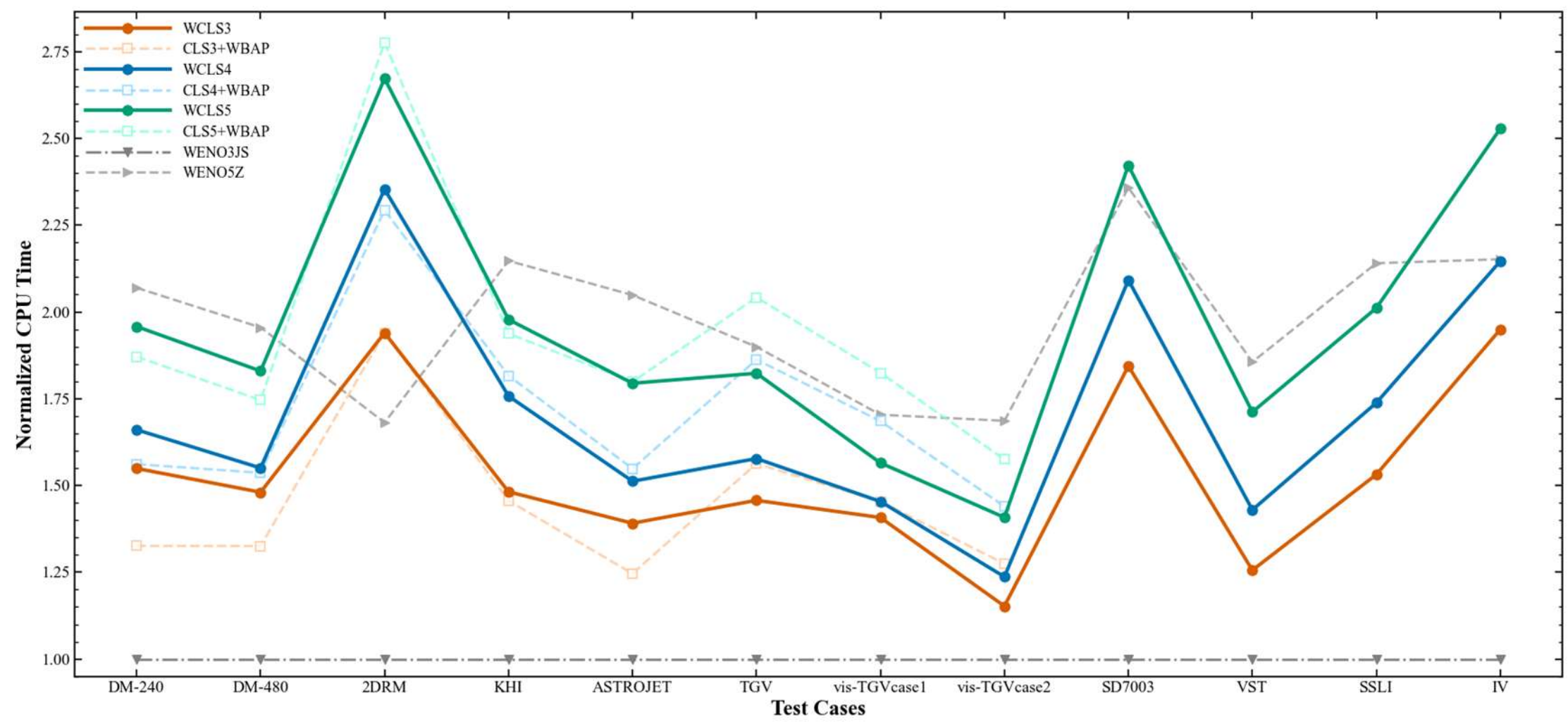

Figure 45: Comparison of normalized computational cost for different unsteady test cases. DM-240 and DM-480 denote the double Mach reflection case with 1/240 and 1/480 grid resolutions, respectively; 2DRM denotes the 2D Riemann problem; KHI denotes Kelvin–Helmholtz instability; ASTROJET denotes the 2D high-Mach-number astrophysical jet; TGV denotes the 3D compressible inviscid Taylor–Green vortex; vis-TGV case 1 and vis-TGV case 2 denote the two viscous Taylor–Green vortex cases; SD7003 denotes the unsteady subsonic viscous flow around the SD7003 airfoil; VST denotes the viscous shock tube; SSLI denotes the shock mixing-layer interaction; and IV denotes the isentropic vortex problem on non-uniform grids.

## 5. Conclusions

In this work, we developed a family of third- to fifth-order weighted compact least-squares schemes for the compressible Navier-Stokes equations on curvilinear grids within an implicit time-integration framework. The proposed formulation reconstructs compact polynomials only along smooth reconstruction lines, employs a compact smoothness indicator based on von Neumann neighbors, and introduces a boundary-value-difference-based penalty mechanism to control the solution in under-resolved regions. The reconstruction is further coupled with the implicit flow solver through generalized minimal residual iterations, thereby avoiding the direct solution of block tridiagonal systems and reducing computational overhead. As a result, the method retains high-order accuracy in smooth regions while maintaining stable, non-oscillatory behavior near discontinuities.

Numerical results for one-, two-, and three-dimensional inviscid and viscous test problems show that the proposed schemes recover the designed orders of accuracy, resolve shocks and contact discontinuities sharply, and remain effective on non-uniform curvilinear grids. Relative to the reference weighted essentially non-oscillatory schemes and weighted-biased-average-procedure-limited compact least-squares schemes, the present formulation exhibits lower numerical dissipation and better preservation of small-scale flow structures in demanding viscous cases, while retaining competitive computational efficiency within the implicit framework. Overall, the proposed weighted compact least-squares approach offers a practical balance among resolution, robustness, and efficiency for compressible flows containing both smooth multiscale structures and strong discontinuities.

### Declaration of competing interests

The authors declare that they have no competing interests relevant to this work.

### CRediT authorship contribution statement

Yongzhi Luo: Methodology, Software, Validation, Formal analysis, Visualization, Writing – original draft.
Huiheng Fan: Methodology, Software, Validation, Visualization, Writing – original draft.
Wei-Gang Zeng: Methodology, Formal analysis, Funding acquisition, Writing – review & editing.
Yu-Xin Ren: Methodology, Writing – review & editing.
Jianhua Pan: Supervision, Conceptualization, Methodology, Software, Validation, Formal analysis, Funding ac quisition, Writing – original draft, Writing – review & editing.

### Funding

This work was supported by the National Natural Science Foundation of China [12102211, 12292980, 1229 2982], the Yongjiang Youth Innovation Talent Program of Ningbo, China, and the project at Northeast Nor mal University [GFPY202505].

**Data availability statement**

The data, input files, and post-processing materials supporting the findings of this study are available from the corresponding author upon reasonable request.